\documentclass[a4paper,12pt]{article}
\begin{document}

\title{Multidimensional Laplace
transforms over Cayley-Dickson algebras and partial differential
equations.}
\author{Ludkovsky S.V.}
\date{25 January 2010}
\maketitle
\begin{abstract}
Multidimensional noncommutative Laplace transforms over octonions
are studied. Theorems about direct and inverse transforms and other
properties of the Laplace transforms over the Cayley-Dickson
algebras are proved. Applications to partial differential equations
including that of elliptic, parabolic and hyperbolic type are
investigated. Moreover, partial differential equations of higher
order with real and complex coefficients and with variable
coefficients with or without boundary conditions are considered.
\end{abstract}

\section{Introduction.}
The Laplace transform over the complex field is already classical
and plays very important role in mathematics including complex
analysis and differential equations \cite{vladumf,lavrsch,polbremm}.
The classical Laplace transform is used frequently for ordinary
differential equations and also for partial differential equations
sufficiently simple to be resolved, for example, of two variables.
But it meets substantial difficulties or does not work for general
partial differential equations even with constant coefficients
especially for that of hyperbolic type. \par To overcome these
drawbacks of the classical Laplace transform in the present paper
more general noncommutative multiparameter transforms over
Cayley-Dickson algebras are investigated. In the preceding paper a
noncommutative analog of the classical Laplace transform over the
Cayley-Dickson algebras was defined and investigated
\cite{lutsltjms}. This paper is devoted to its generalizations for
several real parameters and also variables in the Cayley-Dickson
algebras. For this the preceding results of the author on
holomorphic, that is (super)differentiable functions, and
meromorphic functions of the Cayley-Dickson numbers are used
\cite{ludoyst,ludfov}. The super-differentiability of functions of
Cayley-Dickson variables is stronger than the Fr\'echet's
differentiability. In those works also a noncommutative line
integration was investigated.
\par We remind that quaternions and
operations over them had been first defined and investigated by W.R.
Hamilton in 1843 \cite{hamilt}. Several years later on Cayley and
Dickson had introduced generalizations of quaternions known now as
the Cayley-Dickson algebras \cite{baez,kansol,kurosh,rothe}. These
algebras, especially quaternions and octonions, have found
applications in physics. They were used by Maxwell, Yang and Mills
while derivation of their equations, which they then have rewritten
in the real form because of the insufficient development of
mathematical analysis over such algebras in their time
\cite{emch,guetze,lawmich}. This is important, because
noncommutative gauge fields are widely used in theoretical physics
\cite{solov}.
\par Each Cayley-Dickson algebra ${\cal A}_r$ over the real
field $\bf R$ has $2^r$ generators $\{ i_0,i_1,...,i_{2^r-1} \} $
such that $i_0=1$, $i_j^2=-1$ for each $j=1,2,...,2^r-1$,
$i_ji_k=-i_ki_j$ for every $1\le k\ne j \le 2^r-1$, where $r\ge 1$.
The algebra ${\cal A}_{r+1}$ is formed from the preceding algebra
${\cal A}_r$ with the help of the so-called doubling procedure by
generator $i_{2^r}$. In particular, ${\cal A}_1=\bf C$  coincides
with the field of complex numbers, ${\cal A}_2=\bf H$ is the skew
field of quaternions, ${\cal A}_3$ is the algebra of octonions,
${\cal A}_4$ is the algebra of sedenions. This means that a sequence
of embeddings $...\hookrightarrow {\cal A}_r\hookrightarrow {\cal
A}_{r+1}\hookrightarrow ...$ exists. \par Generators of the
Cayley-Dickson algebras have a natural physical meaning as
generating operators of fermions. The skew field of quaternions is
associative, and the algebra of octonions is alternative. The
Cayley-Dickson algebra ${\cal A}_r$ is power associative, that is,
$z^{n+m}=z^nz^m$ for each $n, m \in \bf N$ and $z\in {\cal A}_r$. It
is non-associative and non-alternative for each $r\ge 4$. A
conjugation $z^*={\tilde z}$ of Cayley-Dickson numbers $z\in {\cal
A}_r$ is associated with the norm $|z|^2 = zz^* = z^*z$. The
octonion algebra has the multiplicative norm and is the division
algebra. Cayley-Dickson algebras ${\cal A}_r$ with $r\ge 4$ are not
division algebras and have not multiplicative norms. The conjugate
of any Cayley-Dickson number $z$ is given by the formula:
\par $(M1)$ $z^* := \xi ^* - \eta {\bf l}$. \\
The multiplication in ${\cal A}_{r+1}$ is defined by the following
equation:
\par $(M2)$ $(\xi + \eta {\bf l})(\gamma +\delta {\bf l})=(\xi \gamma
-{\tilde {\delta }}\eta )+(\delta \xi +\eta {\tilde {\gamma }}){\bf l}$ \\
for each $\xi $, $\eta $, $\gamma $, $\delta \in {\cal A}_r$, $z :=
\xi +\eta {\bf l}\in {\cal A}_{r+1}$, $\zeta :=\gamma +\delta {\bf
l} \in {\cal A}_{r+1}$.
\par At the beginning of this article a multiparameter noncommutative transform is
defined. Then new types of the direct and inverse noncommutative
multiparameter transforms over the general Cayley-Dickson algebras
are investigated, particularly, also over the quaternion skew field
and the algebra of octonions. The transforms are considered in
${\cal A}_r$ spherical and ${\cal A}_r$ Cartesian coordinates. At
the same time specific features of the noncommutative multiparameter
transforms are elucidated, for example, related with the fact that
in the Cayley-Dickson algebra ${\cal A}_r$ there are $2^r-1$
imaginary generators $\{ i_1,...,i_{2^r-1} \} $ apart from one in
the field of complex numbers such that the imaginary space in ${\cal
A}_r$ has the dimension $2^r-1$. Theorems about properties of images
and originals in conjunction with the operations of linear
combinations, differentiation, integration, shift and homothety are
proved. An extension of the noncommutative multiparameter transforms
for generalized functions is given. Formulas for noncommutative
transforms of products and convolutions of functions are deduced.
\par Thus this solves the problem of non-commutative mathematical
analysis to develop the multiparameter Laplace transform over the
Cayley-Dickson algebras. Moreover, an application of the
noncommutative integral transforms for solutions of partial
differential equations is described. It can serve as an effective
means (tool) to solve partial differential equations with real or
complex coefficients with or without boundary conditions and their
systems of different types. An algorithm is described which permits
to write fundamental solutions and functions of Green's type. A
moving boundary problem and partial differential equations with
discontinuous coefficients are also studied with the use of the
noncommutative transform.
\par Moreover, a decomposition theorem of linear partial
differential operators over the Cayley-Dickson algebras is proved. A
relation between fundamental solutions of an initial and component
operators is demonstrated. In conjunction with a line integration
over the Cayley-Dickson algebras and the decomposition theorem of
partial differential operators it permits to solve partial
differential equations linear with constant and variable
coefficients and non-linear as well as boundary problems (see also
\cite{ludifeqcdla}). Certainly, this approach effectively
encompasses systems of partial differential equations, because each
function $f$ with values in the Cayley-Dickson algebra is the sum of
functions $f_ji_j$, where each function $f_j$ is real-valued.
\par All results of this paper are obtained for the first time.

\section{Multidimensional noncommutative integral transforms.}
\par {\bf 1. Definitions. Transforms in ${\cal A}_r$ Cartesian coordinates.}
\par Denote by ${\cal A}_r$ the Cayley-Dickson algebra, $0\le r$,
which may be, in particular, ${\bf H} = {\cal A}_2$ the quaternion
skew field or ${\bf O} = {\cal A}_3$ the octonion algebra. For
unification of the notation we put ${\cal A}_0 = {\bf R}$, ${\cal
A}_1 = {\bf C}$. A function $f: {\bf R}^n\to {\cal A}_r$ we call a
function-original, where $2\le r$, $n\in \bf N$, if it fulfills the
following conditions $(1-5)$.
\par $(1).$ The function $f(t)$ is almost everywhere continuous on
${\bf R}^n$ relative to the Lebesgue measure $\lambda _n$ on ${\bf
R}^n$.
\par $(2).$ On each finite interval in $\bf R$ each function
$g_j(t_j)= f(t_1,...,t_n)$ by $t_j$ with marked all other variables
may have only a finite number of points of discontinuity of the
first kind, where $t=(t_1,...,t_n)\in {\bf R}^n$, $t_j\in \bf R$,
$j=1,...,n$. Recall that a point $u_0\in \bf R$ is called a point of
discontinuity of the first type, if there exist finite left and
right limits $\lim_{u\to u_0, u<u_0} g(u) =: g(u_0-0)\in {\cal A}_r$
and $\lim_{u\to u_0, u>u_0} g(u) =: g(u_0+0)\in {\cal A}_r$.
\par $(3).$ Every partial function $g_j(t_j)=f(t_1,...,
t_n)$ satisfies the H\"older condition: $|g_j(t_j+h_j)-g_j(t_j)| \le
A_j |h_j|^{\alpha _j}$ for each $|h_j|<\delta $, where $0<\alpha
_j\le 1$, $A_j=const
>0$, $\delta _j>0$ are constants for a given $t=(t_1,...,
t_n)\in {\bf R}^n$, $j=1,...,n$, everywhere on ${\bf R}^n$ may be
besides points of discontinuity of the first type.
\par $(4).$ The function $f(t)$ increases not faster, than
the exponential function, that is there exist constants $C_v = const
>0$, $v= (v_1,...,v_n)$, $a_{-1}, a_1 \in \bf R$, where $v_j\in \{
-1, 1 \}$ for every $j=1,...,n$, such that
\par $|f(t)|<C_v \exp ((q_v,t))$ for each $t\in {\bf R}^n$ with
$t_j v_j\ge 0$ for each $j=1,...,n$, $q_v =
(v_1a_{v_1},...,v_na_{v_n})$; where
\par $(5)$ $(x,y) := \sum_{j=1}^n x_j y_j$ denotes
the standard scalar product in ${\bf R}^n$.
\par Certainly for a bounded original $f$ it is possible to take
$a_{-1} = a_1 = 0$.
\par Each Cayley-Dickson number $p\in {\cal A}_r$ we write in the
form \par $(6)$ $p = \sum_{j=0}^{2^r-1} p_j i_j$, where $ \{ i_0,
i_1, ...,i_{2^r-1} \} $ is the standard basis of generators of
${\cal A}_r$ so that $i_0=1$, $i_j^2=-1$ and $i_0i_j=i_j=i_ji_0$ for
each $j>0$, $i_ji_k = - i_ki_j$ for each $j>0$ and $k>0$ with $k\ne
j$, $p_j\in {\bf R}$ for each $j$.
\par If there exists an integral
\par $(7)$ $F^n(p) := F^n(p;\zeta ):=
\int_{{\bf R}^n} f(t) e^{- <p,t) - \zeta }dt$, \\
then $F^n(p)$ is called the noncommutative multiparameter (Laplace)
transform at a point $p\in {\cal A}_r$ of the function-original
$f(t)$, where $\zeta -\zeta _0 = \zeta _1i_1+...+\zeta
_{2^r-1}i_{2^r-1}\in {\cal A}_r$ is the parameter of an initial
phase, $\zeta _j\in \bf R$ for each $j=0,1,...,2^r-1$, $\zeta \in
{\cal A}_r$, $n=2^r-1$, $dt = \lambda _n(dt)$,
\par $(8)$ $<p,t) =p_0 (t_1+...+t_{2^r-1}) +
\sum_{j=1}^{2^r-1} p_j t_j i_j$, we also put \par $(8.1)$
$u(p,t;\zeta ) = <p,t) + \zeta $.
\par For vectors $v, w \in {\bf R}^n$ we shall consider a partial ordering
\par $(9)$ $v\prec w$ if and only if $v_j\le w_j$ for each $j=1,...,n$ and a
$k$ exists so that $v_k<w_k$, $1\le k \le n$.

\par {\bf 2. Transforms in ${\cal A}_r$ spherical coordinates.}
\par Now we consider also the non-linear function $u=u(p,t;\zeta )$ taking
into account non commutativity of the Cayley-Dickson algebra ${\cal
A}_r$.  Put \par $(1)$ $u(p,t) := u(p,t;\zeta ) := p_0 s_1 +
M(p,t)+\zeta _0$, where
\par $(2)$ $M(p,t)=M(p,t;\zeta ) = (p_1s_1+\zeta _1)[ i_1 \cos
(p_2s_2 +\zeta _2) + i_2 \sin (p_2s_2+\zeta _2)$ \\
$\cos (p_3s_3+\zeta _3) +...+ i_{2^r-2} \sin (p_2s_2+\zeta _2)
...\sin (p_{2^r-2}s_{2^r-2}+\zeta _{2^r-2}) \cos
(p_{2^r-1}s_{2^r-1}$
\\ $+\zeta _{2^r-1}) + i_{2^r-1}\sin (p_2s_2+\zeta _2)...\sin
(p_{2^r-2}s_{2^r-2}+\zeta _{2^r-2}) \sin
(p_{2^r-1}s_{2^r-1} + \zeta _{2^r-1})]$ \\
for the general Cayley-Dickson algebra with $2\le r<\infty $,
\par $(2.1) \quad s_j := s_j(n;t) := t_j +...+t_n$ for each $j=1,...,n$,
$n=2^r-1$, so that $s_1=t_1+...+t_n$, $s_n=t_n$.  More generally,
let
\par $(3)$ $u(p,t)=u(p,t;\zeta )=p_0 s_1 + w(p,t)+\zeta _0$, where
$w(p,t)$ is a locally analytic function, $Re (w(p,t))=0$ for each
$p\in {\cal A}_r$ and $t\in {\bf R}^{2^r-1}$, $Re (z) := (z+{\tilde
z})/2$, ${\tilde z}=z^*$ denotes the conjugated number for $z\in
{\cal A}_r$. Then the more general non-commutative multiparameter
transform over ${\cal A}_r$ is defined by the formula:
\par $(4)$ $F_u^n(p;\zeta ) := \int_{{\bf R}^n} f(t)
\exp (-u(p,t;\zeta ))dt$ \\
for each Cayley-Dickson numbers $p\in {\cal A}_r$ whenever this
integral exists as the principal value of either Riemann or Lebesgue
integral, $n=2^r-1$. This non-commutative multiparameter transform
is in ${\cal A}_r$ spherical coordinates, when $u(p,t;\zeta )$ is
given by Formulas $(1,2)$.
\par At the same time the components $p_j$ of the number $p$ and $\zeta
_j$ for $\zeta $ in $u(p,t;\zeta )$ we write in the $p$- and $\zeta
$-representations respectively such that
\par $(5)$ $h_j=(-hi_j+ i_j(2^r-2)^{-1} \{ -h
+\sum_{k=1}^{2^r-1}i_k(hi_k^*) \} )/2$ for each $j=1,2,...,2^r-1$, \\
\par $(6)$ $h_0=(h+ (2^r-2)^{-1} \{ -h +
\sum_{k=1}^{2^r-1}i_k(hi_k^*) \} )/2$, \\
where $2\le r\in \bf N$, $h=h_0i_0+...+h_{2^r-1}i_{2^r-1}\in {\cal
A}_r$, $h_j\in \bf R$ for each $j$, $i_k^* = {\tilde i}_k = - i_k$
for each $k>0$, $i_0=1$, $h\in {\cal A}_r$. Denote $F_u^n(p;\zeta )$
in more details by ${\cal F}^n(f,u;p;\zeta )$.
\par Henceforth, the functions $u(p,t; \zeta )$ given by 1$(8,8.1)$
or $(1,2, 2.1)$ are used, if another form $(3)$ is not specified. If
for $u(p,t; \zeta )$ concrete formulas are not mentioned, it will be
undermined, that the function $u(p,t; \zeta )$ is given in ${\cal
A}_r$ spherical coordinates by Expressions $(1,2, 2.1)$. If in
Formulas 1$(7)$ or $(4)$ the integral is not by all, but only by
$t_{j(1)},...,t_{j(k)}$ variables, where $1\le k<n$, $1\le
j(1)<...<j(k)\le n$, then we denote a noncommutative transform by
$F_u^{k; t_{j(1)},...,t_{j(k)}} (p;\zeta )$ or ${\cal F}^{k;
t_{j(1)},...,t_{j(k)}} (f,u;p;\zeta )$. If $j(1)=1,$...,$j(k)=k$,
then we denote it shortly by $F^k_u(p;\zeta )$ or ${\cal
F}^k(f,u;p;\zeta )$. Henceforth, we take $\zeta _m=0$ and $t_m=0$
and $p_m=0$ for each $1\le m\notin \{ j(1),...,j(k) \} $ if
something other is not specified.

\par {\bf 3. Remark.} The spherical ${\cal A}_r$ coordinates
appear naturally from the following consideration of iterated
exponents:
\par $(1)$ $\exp (i_1(p_1s_1+\zeta _1)\exp (-i_3(p_2s_2
+\zeta _2) \exp (-i_1(p_3s_3+\zeta _3)))= \exp (i_1(p_1s_1+\zeta
_1)\exp (- (p_2s_2+\zeta _2)(i_3\cos (p_3s_3+\zeta _3) - i_2\sin
(p_3s_3+\zeta _3))))$
\par $= \exp (i_1(p_1s_1+\zeta _1)(\cos (p_2s_2+\zeta _2) -
\sin (p_2s_2+\zeta _2)(i_3\cos (p_3s_3+\zeta _3) - i_2\sin
(p_3s_3+\zeta _3))))$
\par $= \exp ((p_1s_1+\zeta _1)(i_1\cos (p_2s_2+\zeta _2) + i_2
\sin (p_2s_2+\zeta _2)\cos (p_3s_3+\zeta _3) + i_3\sin (p_2s_2+\zeta
_2)\sin (p_3s_3+\zeta _3)))$. \par Consider $i_{2^r}$ the generator
of the doubling procedure of the Cayley-Dickson algebra ${\cal
A}_{r+1}$ from the Cayley-Dickson algebra ${\cal A}_r$, such that
$i_ji_{2r}=i_{2^r+j}$ for each $j=0,...,2^r-1$. We denote now the
function $M(p,t;\zeta )$ from Definition 2 over ${\cal A}_r$ in more
details by $\mbox{ }_rM$.
\par Then by induction we write:
$$(2)\quad \exp (\mbox{ }_{r+1}M(p,t;\zeta ))= \exp \{\mbox{
}_rM((i_1p_1+_...+i_{2^r-1}p_{2^r-1}),(t_1,...,t_{2^r-2},
(t_{2^r-1}+s_{2^r}));$$   $$ (i_1\zeta _1+...+i_{2^r-1} \zeta
_{2^r-1}) \exp (-i_{2^r+1}(p_{2^r}s_{2^r} +\zeta _{2^r})$$  $$ \exp
(-\mbox{ }_rM((i_1p_{2^r+1}+...+i_{2^r-1}p_{2^{r+1}-1}),
(t_{2^r+1},...,t_{2^{r+1}-1}); (i_1\zeta _{2^r+1}+...+i_{2^r-1}
\zeta _{2^{r+1}-1}))) \} ,$$  where $t=(t_1,...,t_n)$, $n=n(r+1) =
2^{r+1}-1$, $s_j=s_j(n(r+1);t)$ for each $j=1,...,n(r+1)$, since
$s_m(n(r+1);t) = t_m+...+t_{n(r+1)} = s_m(n(r);t) +
s_{2^r}(n(r+1);t)$ for each $m=1,...,2^r-1$.
\par An image function can be written in the form
\par $(3)$ $F_u^n(p;\zeta ):=
\sum_{j=0}^{2^r-1} i_j F_{u,j}^n(p;\zeta )$, \\ where a function $f$
is decomposed in the form \par $(3.1)$ $f(t)=\sum_{j=0}^{2^r-1} i_j
f_j(t)$, \\ $f_j: {\bf R}^n\to \bf R$ for each $j=0,1,...,2^r-1$,
$F_{u,j}^n(p;\zeta )$ denotes the image of the function-original
$f_j$.
\par If an automorphism of the Cayley-Dickson algebra ${\cal A}_r$
is taken and instead of the standard generators $ \{
i_0,...,i_{2^r-1} \} $ new generators $ \{ N_0,...,N_{2^r-1} \} $
are used, this provides also $M(p,t;\zeta )=M_N(p,t;\zeta )$
relative to new basic generators, where $2\le r\in \bf N$. In this
more general case we denote by $\mbox{ }_NF_u^n(p;\zeta )$ an image
for an original $f(t)$, or in more details we denote it by $\mbox{
}_N {\cal F}^n(f,u;p;\zeta )$. \par Formulas 1$(7)$ and 2$(4)$
define the right multiparameter transform. Symmetrically is defined
a left multiparameter transform. They are related by conjugation and
up to a sign of basic generators. For real valued originals they
certainly coincide. Henceforward, only the right multiparameter
transform is investigated.
\par Particularly, if $p=(p_0,p_1,0,...,0)$ and $t=(t_1,0,...,0)$,
then the multiparameter non-commutative Laplace transforms 1$(7)$
and 2$(4)$ reduce to the complex case, with parameters $a_1$,
$a_{-1}$. Thus, the given above definitions over quaternions,
octonions and general Cayley-Dickson algebras are justified.

\par {\bf 4. Theorem.} {\it If an original $f(t)$ satisfies
Conditions 1$(1-4)$ and $a_1<a_{-1}$, then its image ${\cal
F}^n(f,u;p;\zeta )$ is ${\cal A}_r$-holomorphic (that is locally
analytic) by $p$ in the domain $\{ z\in {\cal A}_r: a_1< Re
(z)<a_{-1} \} $, as well as by $\zeta \in {\cal A}_r$, where $1\le
r\in \bf N$, $2^{r-1}\le n \le 2^r-1$, the function $u(p,t; \zeta )$
is given by 1$(8,8.1)$ or 2$(1,2, 2.1)$.}
\par {\bf Proof.} At first consider the characteristic functions
$\chi _{U_v} (t)$, where $\chi _U(t) =1$ for each $t\in U$, while
$\chi _U(t)=0$ for every $t\in {\bf R}^n\setminus U$, $U_v := \{
t\in {\bf R}^n: v_jt_j\ge 0 ~ \forall j=1,...,n \} $ is the domain
in the Euclidean space ${\bf R}^n$ for any $v$ from \S 1. Therefore,
\par $(1)$ $F_u^n(p;\zeta ) := \sum_{[v=(v_1,...,v_n):
v_1,...,v_n \in \{ -1, 1 \} ]} \int_{U_v} f(t) \exp (-u(p,t;\zeta )) dt,$ \\
since $\lambda _n (U_v\cap U_w)= 0$ for each $v\ne w$. Each integral
$\int_{U_v} f(t) \exp (-u(p,t;\zeta )) dt$ is absolutely convergent
for each $p\in {\cal A}_r$ with the real part $ a_1< Re (p)
<a_{-1}$, since it is majorized by the converging integral
\par $(2)$ $|\int_{U_v} f(t) \exp (-u(p,t;\zeta )) dt| \le \int_0^{\infty }...
\int_0^{\infty } C_v \exp \{ - v_1(w - a_{v_1})y_1-...-v_n(w -
a_{v_n})y_n - \zeta _0 \} dy_1...dy_n =
C_v e^{-\zeta _0} \prod_{j=1}^n v_j(w-a_{v_j})^{-1} $, \\
where $w=Re (p)$, since $|e^z|=\exp (Re (z))$ for each $z\in {\cal
A}_r$ in view of Corollary 3.3 \cite{ludfov}. While an integral,
produced from the integral $(1)$ differentiating by $p$ converges
also uniformly:  $$(3)\quad |\int_{U_v}f(t)[\partial \exp
(-u(p,t;\zeta ))/ \partial p].hdt| $$
$$\le \int_0^{\infty }...\int_0^{\infty } C_v
|(h_0(v_1y_1+...+v_ny_n),h_1(v_1y_1+...+v_ny_n),...,h_{n-1}(v_{n-1}y_{n-1}+v_ny_n),
h_nv_ny_n)|$$ $\exp \{ - v_1(w - a_{v_1})y_1-...-v_n(w - a_{v_n})y_n
- \zeta _0 \} dy_1...dy_n$
$$\le |h|C_v e^{-\zeta _0} \prod_{j=1}^n (w-a_{v_j})^{-2} $$
for each $h\in {\cal A}_r$, since each $z\in {\cal A}_r$ can be
written in the form $z=|z|\exp (M)$, where $|z|^2=z{\tilde z}\in
[0,\infty )\subset \bf R$, $M\in {\cal A}_r$, $Re (M):= (M+{\tilde
M})/2=0$ in accordance with Proposition 3.2 \cite{ludfov}. In view
of Equations 2$(5,6)$:
\par $(4)$ $\partial (\int_{{\bf R}^n}f(t)\exp (- u(p,t;\zeta
))dt)/\partial {\tilde p}=0$ and \par $(5)$ $\partial (\int_{{\bf
R}^n}f(t)\exp (- u(p,t;\zeta ))dt)/\partial {\tilde \zeta }=0$,
while
\par $(6)$ $|\int_{U_v} f(t) [\partial \exp (- u(p,t;\zeta ))/\partial
\zeta ].hdt| \le |h| \int_0^{\infty }...\int_0^{\infty } C_v \exp \{
- v_1(w - a_{v_1})y_1-...-v_n(w - a_{v_n})y_n - \zeta _0 \}
dy_1...dy_n = |h|
C_v e^{-\zeta _0}\prod_{j=1}^n v_j(w-a_{v_j})^{-1}$ \\
for each $h\in {\cal A}_r$. In view of convergence of integrals
given above $(1-6)$ the multiparameter non-commutative transform
$F_u^n(p;\zeta )$ is (super)differentiable by $p$ and $\zeta $,
moreover, $\partial F_u^n(p;\zeta )/\partial {\tilde p}=0$ and
$\partial F_u^n(p;\zeta )/\partial {\tilde \zeta }=0$ in the
considered $(p,\zeta )$-representation. In accordance with
\cite{ludoyst,ludfov} a function $g(p)$ is locally analytic by $p$
in an open domain $U$ in the Cayley-Dickson algebra ${\cal A}_r$,
$2\le r$, if and only if it is (super)differentiable by $p$, in
another words ${\cal A}_r$-holomorphic. Thus, $F_u^n(p;\zeta )$ is
${\cal A}_r$-holomorphic  by $p\in {\cal A}_r$ with $a_1<Re (p)<
a_{-1}$ and $\zeta \in {\cal A}_r$ due to Theorem 2.6
\cite{lutsltjms}.
\par {\bf 4.1. Corollary.} {\it Let suppositions of Theorem 4 be
satisfied. Then the image ${\cal F}^n(f,u;p;\zeta )$ with
$u=u(p,t;\zeta )$ given by 2$(1,2)$ has the following periodicity
properties:
\par $(1)$ ${\cal F}^n(f,u;p;\zeta +\beta i_j) = {\cal
F}^n(f,u;p;\zeta )$ for each $j=1,...,n$ and $\beta \in 2\pi {\bf
Z}$;
\par $(2)$ ${\cal F}^n(f,u;p^1;\zeta ^1) = (-1)^{\kappa }
{\cal F}^n(f,u;p^2;\zeta ^2)$ for each $j=1,...,n-1$ so that $\zeta
_0^1 = \zeta _0^2$ and $\zeta _j^1 = - \zeta _j^2$, $\zeta _{j+1}^1
= \pi + \zeta _{j+1}^2$, $\zeta _s^1 = \zeta _s^2$ for each $s\ne j$
and $s\ne j+1$, while either $p_j^1 = - p_j^2$ and $p_l^1=p_l^2$ for
each $l\ne j$ with $\kappa =2$ or $p^1=p^2$ and $f(t)$ is an even
function with $\kappa =2$ by the $s_j=(t_j+...+t_n)$ variable or an
odd function by $s_j=(t_j+...+t_n)$ with $\kappa =1$;
\par $(3)$ ${\cal F}^n(f,u;p;\zeta + \pi i_1) = - {\cal
F}^n(f,u;p;\zeta )$.}
\par {\bf Proof.} In accordance with Theorem 4 the image ${\cal
F}^n(f,u;p;\zeta )$ exists for each $p\in W_f := \{ z\in {\cal A}_r:
~ a_1< Re (z) < a_{-1} \} $ and $\zeta \in {\cal A}_r$, where $1\le
r$. Then from the $2\pi $ periodicity of sine and cosine functions
the first statement follows. From $\sin (-\phi ) = - \sin (\phi )$,
$\cos ( \phi ) = \cos ( - \phi )$, $\sin (\pi +\phi ) = - \sin (\phi
)$, $\cos (\phi +\pi ) = - \cos (\phi )$ we get that $\cos ( p_j s_j
+ \zeta _j^1) = \cos ( - p_j s_j + \zeta _j^2)$, $\sin (p_j s_j +
\zeta _j^1) \cos (p_{j+1} s_{j+1} + \zeta _{j+1}^1) = ( - \sin ( -
p_j s_j + \zeta _j^2) ) ( - \cos (p_{j+1} s_{j+1} + \zeta _{j+1}^2)
)$ and $\sin (p_j s_j + \zeta _j^1) \sin (p_{j+1} s_{j+1} + \zeta
_{j+1}^1) = ( - \sin ( - p_j s_j + \zeta _j^2) ) ( - \sin ( p_{j+1}
s_{j+1} + \zeta _{j+1}^2) )$ for each $t\in {\bf R}^n$. On the other
hand, either $p_j^1 = - p_j^2$ and $p_l^1=p_l^2$ for each $l\ne j\ge
1$ with $\kappa =2$ or $p^1=p^2$ and $f(t_1,...,s_{j-1}+s_j,
-s_j-s_{j+1}, t_{j+1},...,t_n) = (-1)^{\kappa }
f(t_1,...,s_{j-1}-s_j, s_j-s_{j+1},t_{j+1},...,t_n)$ is an even with
$\kappa =2$ or odd with $\kappa =1$ function by the
$s_j=(t_j+...+t_n)$ variable for each $t=(t_1,...,t_n)\in {\bf
R}^n$, where $t_j=s_j-s_{j+1}$ for $j=1,...,n$,
$s_{n+1}=s_{n+1}(n;t)=0$. From this and Formulas 2$(1,2,4)$ the
second and the third statements of this corollary follow.

\par {\bf 5. Remark.} For a subset $U$ in ${\cal A}_r$ we put
$\pi _{{\sf s},{\sf p},{\sf t}}(U):= \{ {\sf u}: z\in U,
z=\sum_{{\sf v}\in \bf b}w_{\sf v}{\sf v},$ ${\sf u}=w_{\sf s}{\sf
s}+w_{\sf p}{\sf p} \} $ for each ${\sf s}\ne {\sf p}\in \bf b$,
where ${\sf t}:=\sum_{{\sf v}\in {\bf b}\setminus \{ {\sf s}, {\sf
p} \} } w_{\sf v}{\sf v} \in {\cal A}_{r,{\sf s},{\sf p}}:= \{ z\in
{\cal A}_r:$ $z=\sum_{{\sf v}\in \bf b} w_{\sf v}{\sf v},$ $w_{\sf
s}=w_{\sf p}=0 ,$ $w_{\sf v}\in \bf R$ $\forall {\sf v}\in {\bf b}
\} $, where ${\bf b} := \{ i_0,i_1,...,i_{2^r-1} \} $ is the family
of standard generators of the Cayley-Dickson algebra ${\cal A}_r$.
That is, geometrically $\pi _{{\sf s},{\sf p},{\sf t}}(U)$ means the
projection on the complex plane ${\bf C}_{{\sf s},{\sf p}}$ of the
intersection $U$ with the plane ${\tilde \pi }_{{\sf s},{\sf p},{\sf
t}}\ni {\sf t}$, ${\bf C}_{{\sf s},{\sf p}} := \{ a{\sf s}+b{\sf
p}:$ $a, b \in {\bf R} \} $, since ${\sf s}{\sf p}^*\in {\hat
b}:={\bf b}\setminus \{ 1 \} $. Recall that in \S \S 2.5-7
\cite{ludfov} for each continuous function $f: U\to {\cal A}_r$ it
was defined the operator ${\hat f}$ by each variable $z\in {\cal
A}_r$. For the non-commutative integral transformations consider,
for example, the left algorithm of calculations of integrals. \par A
Hausdorff topological space $X$ is said to be $n$-connected for
$n\ge 0$ if each continuous map $f: S^k\to X$ from the
$k$-dimensional real unit sphere into $X$ has a continuous extension
over $\bf R^{k+1}$ for each $k\le n$  (see also \cite{span}). A
$1$-connected space is also said to be simply connected.
\par  It is supposed further, that a domain $U$ in ${\cal A}_r$
has the property that $U$ is $(2^r-1)$-connected; $\pi _{{\sf
s},{\sf p},{\sf t}}(U)$ is simply connected in $\bf C$ for each
$k=0,1,...,2^{r-1}$, ${\sf s}=i_{2k}$, ${\sf p}=i_{2k+1}$, ${\sf
t}\in {\cal A}_{r,{\sf s},{\sf p}}$ and ${\sf u}\in {\bf C}_{{\sf
s},{\sf p}}$, for which there exists $z={\sf u}+{\sf t}\in U$.

\par {\bf 6. Theorem.} {\it If a function $f(t)$ is an original
(see Definition 1), such that $\mbox{ }_NF_u^n(p;\zeta )$ is its
image multiparameter non-commutative transform, where the functions
$f$ and $F_u^n$ are written in the forms given by 3$(3, 3.1)$,
$f({\bf R}^n)\subset {\cal A}_r$ over the Cayley-Dickson algebra
${\cal A}_r$, where $1\le r\in \bf N$, $2^{r-1}\le n\le 2^r-1$.
\par  Then at each point $t$, where $f(t)$ satisfies the H\"older
condition the equality is accomplished :
$$(1)\quad f(t) = \{ [(2\pi N_n)^{-1}
\int_{-N_n\infty }^{N_n\infty }](... ([(2\pi N_1)^{-1}
\int_{-N_1\infty }^{N_1\infty }] \mbox{ }_NF_u^n(a+p;\zeta )$$
$$\exp \{ u(a+p,t;\zeta ) \} )...)dp \} =: ({\cal F}^n )^{-1} (\mbox{ }_NF_u^n(a+p;\zeta ),u,t;\zeta ) ,$$
where either $u(p,t;\zeta )=<p,t) + \zeta $ or $u(p,t;\zeta )=p_0
s_1 + M_N(p,t;\zeta )+\zeta _0$ (see \S \S 1 and 2), the integrals
are taken along the straight lines $p(\tau _j)=N_j\tau _j\in {\cal
A}_r$, $\tau _j\in \bf R$ for each $j=1,...,n$; $a_1< Re (p) = a <
a_{-1}$ and this integral is understood in the sense of the
principal value, $t=(t_1,...,t_n)\in {\bf R}^n$,
$dp=(...((d[p_1N_1])d[p_2N_2])...)d[p_nN_n]$.}
\par {\bf Proof.} In Integral $(1)$ an integrand $\eta (p)dp$
certainly corresponds to the iterated integral as $(...(\eta
(p)d[p_1N_1])...)d[p_nN_n]$, where $p = p_1N_1+...+p_nN_n$,
$p_1,...,p_n\in \bf R$. Using Decomposition 3$(3.1)$ of a function
$f$ it is sufficient to consider the inverse transformation of the
real valued function $f_j$, which we denote for simplicity by $f$.
We put $$\mbox{ }_NF^n_{u,j}(p;\zeta ) := \int_{{\bf R}^n}f_j(t)\exp
(-u(p,t;\zeta ))dt .$$  If $\eta $ is a holomorphic function of the
Cayley-Dickson variable, then locally in a simply connected domain
$U$ in each ball $B({\cal A}_r,z_0,R)$ with the center at $z_0$ of
radius $R>0$ contained in the interior $Int (U)$ of the domain $U$
there is accomplished the equality
\par $(\partial [\int_{z_0}^z\eta (a+\zeta )
d\zeta ]/\partial z).1=\eta (a+z)$, \\
where the integral depends only on an initial $z_0$ and a final $z$
points of a rectifiable path in $B({\cal A}_r,z_0,R)$, $a\in {\bf
R}$ (see also Theorem 2.14 \cite{lutsltjms}). Therefore, along the
straight line $N_j{\bf R}$ the restriction of the antiderivative has
the form $\int_{\theta _0}^{\theta }\eta (a+N_j\tau _j)d\tau _j$,
since \par $(2)$ $\int_{z_0=N_j\theta _0}^{z=N_j\theta }\eta
(a+\zeta )d\zeta =\int_{\theta _0}^{\theta } {\hat
{\eta }}(a+N_j\tau _j).N_jd\tau _j$, \\
where $\partial \eta (a+z)/\partial \theta =(\partial \eta
(a+z)/\partial z).N_j$ for the (super)differentiable by $z\in U$
function $\eta (z)$, when $z=\theta N_j$, $\theta \in {\bf R}$. For
the chosen branch of the line integral specified by the left
algorithm this antiderivative is unique up to a constant from ${\cal
A}_r$ with the given $z$-representation $\nu $ of the function $\eta
$ \cite{ludfov,ludoyst,lutsltjms}. On the other hand, for analytic
functions with real expansion coefficients in their power series
non-commutative integrals specified by left or right algorithms
along straight lines coincide with usual Riemann integrals by the
corresponding variables. The functions $\sin (z)$, $\cos (z)$ and
$e^z$ participating in the multiparameter non-commutative transform
are analytic with real expansion coefficients in their series by
powers of $z\in {\cal A}_r$.

\par Using Formula 4$(1)$ we reduce
the consideration to $\chi _{U_v}(t)f(t)$ instead of $f(t)$. By
symmetry properties of such domains and integrals and utilizing
change of variables it is sufficient to consider $U_v$ with
$v=(1,...,1)$. In this case $\int_{{\bf R}^n}$ for the direct
multiparameter non-commutative transform 1$(7)$ and 2$(4)$ reduces
to $\int_0^{\infty }...\int_0^{\infty }$. Therefore, we consider in
this proof below the domain $U_{1,...,1}$ only. Using Formulas 3$(3,
3.1)$ and 2$(1,2,2.1)$ we mention that any real algebra with
generators $N_0=1$, $N_k$ and $N_j$ with $1\le k\ne j$ is isomorphic
with the quaternion skew field $\bf H$, since $Re (N_jN_k)=0$ and
$|N_j|=1$, $|N_k|=1$ and $|N_jN_k|= 1$. Then $\exp (\alpha + M\beta
) \exp (\gamma + M \omega ) = \exp ((\alpha + \gamma ) + M (\beta
+\omega ))$ for each real numbers $\alpha , \beta , \gamma , \delta
$ and a purely imaginary Cayley-Dickson number $M$.
\par The octonion algebra $\bf O$ is alternative, while the real field $\bf
R$ is the center of the Cayley-Dickson algebra ${\cal A}_r$. We
consider the integral
\par $(3)$ $g_b(t) := [(2\pi N_n)^{-1}
\int_{-N_nb}^{N_nb}](... ([(2\pi N_1)^{-1} \int_{-N_1b}^{N_1b}]
\mbox{ }_NF_{u,j}^n(a+p;\zeta )\exp \{ u(a+p,t;\zeta ) \} )...)dp$
\\ for each positive value of the parameter $0<b<\infty $. With the
help of generators of the Cayley-Dickson algebra ${\cal A}_r$ and
the Fubini Theorem for real valued components of the function the
integral can be written in the form:
$$(4)\quad g_b(t) = [(2\pi N_n)^{-1} \int_0^{\infty }d\tau _n
\int_{-N_nb}^{N_nb}](...([(2\pi N_1)^{-1} \int_0^{\infty }d\tau _1
\int_{-N_1b}^{N_1b}]$$  $$ f(\tau )\exp \{ - u_N(a+p,t;\zeta ) \}
\exp \{ u_N(a+p,\tau ;\zeta ) \} )...)dp ,$$  since the integral
$\int_{U_{1,...,1}}f(\tau )\exp \{ - u_N(a+p,\tau ;\zeta ) \} d\tau
$ for any marked $0<\delta < (a_{-1} - a_1)/3$ is uniformly
converging relative to $p$ in the domain $a_1+\delta \le Re (p) \le
a_{-1} - \delta $ in ${\cal A}_r$ (see also Proposition 2.18
\cite{lutsltjms}).
\par If take marked $t_k$ for each $k\ne j$ and $S=N_j$ for some
$j\ge 1$ in Lemma 2.17 \cite{lutsltjms} considering the variable
$t_j$, then with a suitable (${\bf R}$-linear) automorphism ${\bf
v}$ of the Cayley-Dickson algebra ${\cal A}_r$ an expression for
${\bf v} (M(p,t;\zeta ))$ simplifies like in the complex case with
${\bf C}_{K} := {\bf R}\oplus {\bf R}K$ for a purely imaginary
Cayley-Dickson number $K$, $|K|=1$, instead of ${\bf C} := {\bf R}
\oplus {\bf R}i_1$, where ${\bf v}(x)=x$ for each real number $x\in
{\bf R}$. But each equality $\alpha = \beta $ in ${\cal A}_r$ is
equivalent to ${\bf v}(\alpha ) = {\bf v}(\beta )$. Then \par $(5)$
$Re [(N_jN_q) (N_jN_l)^*]= Re (N_qN_l^*)=\delta _{q,l}$ for each $q,
l$.  \par If $S^j = \sum_{0\le l\le n; l\ne j} \alpha _lN_l$, $N^j=
\sum_{0\le l\le n; l\ne j} \beta _lN_l$ with $j\ge 1$ and real
numbers $\alpha _l, \beta _l\in {\bf R}$ for each $l$, then
\par $(6)$ $Re [(N_jS^j)(N_jN^j)^*]=Re [S^j(N^j)^*] =\sum_l\alpha
_l\beta _l$.
\par The latter identity can be applied to either
$S^k=M_{k+1}(p_{k+1}N_{k+1}+...+p_nN_n, (t_{k+1},...,t_n); \zeta
_{k+1}N_{k+1}+...+\zeta _nN_n)$ and $N^k =
M_{k+1}(p_{k+1}N_{k+1}+...+p_nN_n, (\tau _{k+1},...,\tau _n); \zeta
_{k+1}N_{k+1}+...+\zeta _nN_n)$, or $S^k= (p_{k+1}t_{k+1}+\zeta
_{k+1})N_{k+1}+...+ (p_nt_n+\zeta _n)N_n$ and $N^k = (p_{k+1}\tau
_{k+1}+\zeta _{k+1})N_{k+1}+...+ (p_n\tau _n+\zeta _n)N_n$, where
\par $(7)$ $M_{k+1}(p_{k+1}N_{k+1}+...+p_nN_n, (t_{k+1},...,t_n); \zeta
_{k+1}N_{k+1}+...+\zeta _nN_n) = (p_{k+1} s_{1,k+1} + \zeta _{k+1})
[N_{k+1} \cos (p_{k+2} s_{2,k+1} + \zeta _{k+2})+...$\par $+ N_n
\sin(p_{k+2} s_{2,k+1} + \zeta _{k+2})...\sin (p_n s_{n-k,k+1} + \zeta _n)]$, \\
\par $(8)$ $s_{j,k+1} = s_{j,k+1}(n;t) = t_{k+j}+...+t_n=s_{k+j}(n;t)$
for each $j=1,...,n-1$; $s_{n-k,k+1}=s_{n-k,k+1}(n;t)=t_n$.
\par  We take the limit of $g_b(t)$ when $b$ tends to the
infinity. Evidently, $s_k(n;\tau ) -s_j(n;\tau ) = s_k(j-1;\tau
)=\tau _k+...+\tau _{j-1}$ for each $1\le k<j\le n$. By our
convention $s_k(n;\tau ) =s_1(n;\tau )$ for $k<1$, while $s_k(n;\tau
)=0$ for $k>n$. Put
\par $(9)$ $u_{n,j}(p_0+p_jN_j+...+p_nN_n, (\tau _j,...,\tau _n);
\zeta _0 + \zeta _jN_j+...+\zeta _nN_n) = \zeta _0 + p_0s_{1,j} +
M_j(p_jN_j+...+p_nN_n, (\tau _j,...,\tau _n); \zeta _0 + \zeta
_jN_j+...+\zeta _nN_n)$ \\ for $u_N$ given by 2$(1,2,2.1)$, where
$M_j$ is prescribed by $(7)$, $s_{k,j}=s_{k,j}(n;\tau )$;
\par $(10)$ $u_{n,j}(p_0+p_jN_j+...+p_nN_n, (\tau _j,...,\tau _n);
\zeta _0 + \zeta _jN_j+...+\zeta _nN_n) = \zeta _0 + p_0s_{1,j} +
\sum_{k=j}^n (p_k\tau _k+\zeta _k)N_k$
\\ for $u=u_N$ given by 1$(8,8.1)$. For $j>1$ the parameter $\zeta
_0$ for $u=u_N$ given by 1$(8,8.1)$ or 2$(1,2,2.1)$ can be taken
equal to zero.
\par  When $t_1,...,t_{j-1}, t_{j+1},...,t_n$ and
$p_1,...,p_{j-1}, p_{j+1},...,p_n$ variables are marked, we take the
parameter \par $\zeta ^j := \zeta ^j (p_jN_j+...+p_nN_n, (\tau
_j,...,\tau _n); \zeta _0+\zeta _jN_j+...+\zeta _nN_n ) := (\zeta
_0+\zeta _jN_j+...+\zeta _nN_n ) + (a+p_0) s_{j+1} + p_{j+1}s_{j+1}
N_{j+1} +...+p_ns_nN_n$ for $u(p,\tau ;\zeta )$ given by Formulas
2$(1,2,2.1)$ or \par $\zeta ^j := \zeta ^j (p_jN_j+...+p_nN_n),(\tau
_j,...,\tau _n); \zeta _0+\zeta _jN_j+...+\zeta _nN_n ) := (\zeta
_0+\zeta _jN_j+...+\zeta _nN_n ) + (a+p_0) s_{j+1} + p_{j+1}\tau
_{j+1} N_{j+1}+...+p_n\tau _nN_n$ for $u(p,\tau ;\zeta )$ described
in 1$(8,8.1)$. Then the integral operator
\\ $\lim_{b\to \infty } [(2\pi N_j)^{-1} \int_0^{\infty }d\tau _j
\int_{-N_jb}^{N_jb}]...(dp_jN_j)$ (see also Formula $(4)$ above)
applied to the function $f(t_1,...,t_{j-1},\tau _j,...,\tau _n)\exp
\{ - u_{N,j}(a+p_0+p_jN_j+...+p_nN_n, (t_j,...,t_n);\zeta _0+\zeta
_jN_j+...+\zeta _nN_n ) \} \exp \{
u_{N,j}(a+p_0+p_jN_j+...+p_nN_n,(\tau _j,...,\tau _n); \zeta
_0+\zeta _jN_j+...+\zeta _nN_n ) \}$ with the parameter $\zeta ^j$
instead of $\zeta $ treated by Theorems 2.19 and 3.15
\cite{lutsltjms} gives the inversion formula corresponding to the
real variable $t_j$ for $f(t)$ and to the Cayley-Dickson variable
$p_0N_0+p_jN_j$ restricted on the complex plane ${\bf C}_{N_j} =
{\bf R}\oplus {\bf R}N_j$, since $d(\tau _j+c) = d\tau _j$ for each
(real) constant $c$. After integrations with $j=1,...,k$ with the
help of Formulas $(6-10)$ and 3$(1,2)$ we get the following:
$$(11)\quad \lim_{b\to \infty} g_b(t) = Re [(2\pi N_n)^{-1}
\int_0^{\infty }d\tau _n \int_{-N_n\infty }^{N_n\infty }](...([(2\pi
N_{k+1})^{-1} \int_0^{\infty }d\tau _{k+1} \int_{-N_{k+1}\infty
}^{N_{k+1}\infty }]$$   $$ f(t_1,...,t_k,\tau _{k+1},...,\tau
_n)\exp \{ - u_{N,k+1} ((a+p_0+p_{k+1}N_{k+1}+...+p_nN_n),
(t_{k+1},...,t_n);$$ $$(\zeta _0+\zeta _{k+1}N_{k+1}+...+\zeta _nN_n
)) \} \exp \{ u_{N,k+1}((a+p_0+p_{k+1}N_{k+1}+...+p_nN_n),$$ $$
(\tau _{k+1},...,\tau _n); (\zeta _0+\zeta _{k+1}N_{k+1}+...+\zeta
_nN_n )) \} )...)dp .$$ Moreover, $Re (f_q)=f_q$ for each $q$ and in
$(11)$ the function $f=f_q$ stands for some marked $q$ in accordance
with Decompositions 3$(3,3.1)$ and the beginning of this proof.
\par Mention, that the algebra $alg_{\bf R}(N_j,N_k,N_l)$ over the real field
with three generators $N_j$, $N_k$ and $N_l$ is alternative. The
product $N_kN_l$ of two generators is also the corresponding
generator $(-1)^{\xi (k,l)} N_m$ with the definite number $m=m(k,l)$
and the sign multiplier $(-1)^{\xi (k,l)}$, where $\xi (k,l)\in \{
0, 1 \} $. On the other hand, $N_{k_1}[{\tilde
N}_j(N_j(N_{k_2}N_l))] = N_{k_1}(N_{k_2}N_l)$. We use decompositions
$(7-10)$ and take $k_2=l$ due to Formula $(11)$, where $Re$ stands
on the right side of the equality, since $Re (N_kN_l)=0$  and $Re
[{\tilde N}_j (N_j(N_kN_l))] =0$ for each $k\ne l$.
 Thus the repeated
application of this procedure by $j=1, 2, ..., n$ leads to Formula
$(1)$ of this theorem.
\par {\bf 6.1. Corollary.} {\it If the conditions of Theorem 6 are
satisfied, then
$$(1)\quad f(t) = (2\pi )^{-n}
\int_{{\bf R}^n} F_u^n(a+p;\zeta ) \exp \{ u(a+p,t;\zeta ) \}
dp_1...dp_n$$  $$ = ({\cal F}^n )^{-1} (\mbox{ }_NF_u^n(a+p;\zeta
),u,t;\zeta ).$$}
\par {\bf Proof.} Each algebra $alg_{\bf R}(N_j,N_k,N_l)$ is alternative.
Therefore, in accordance with \S 6 and Formulas 1$(8,8.1)$ and
2$(1-4)$ for each non-commutative integral given by the left
algorithm we get
$$(2)\quad N_j^{-1}\int_{-N_jb}^{N_jb}[f(\tau )\exp \{ - u_N(a+p,t ;\zeta )
\} ] \exp \{ u_N(a+p,\tau ;\zeta ) \} d(p_jN_j)$$
$$\sum_{l=0}^{2^r-1} {\tilde N}_j [N_j (\int_{-N_jb}^{N_jb}[N_l f_l(\tau )\exp \{ - u_N(a+p,t ;\zeta )
\} ] \exp \{ u_N(a+p,\tau ;\zeta ) \} dp_j)]$$
$$=\int_{-b}^{b} [ f(\tau ) \exp \{ - u_N(a+p,t ;\zeta ) \} ]
\exp \{ u_N(a+p,\tau ;\zeta ) \} dp_j$$ for each $j=1,...,n$, since
the real field is the center of the Cayley-Dickson algebra ${\cal
A}_r$, while the functions $\sin $ and $\cos $ are analytic with
real expansion coefficients. Thus
\par $(3)$ $g_b(t) = (2\pi )^{-n}
[\int_0^{\infty }d\tau _n \int_{-b}^{b}](...([ \int_0^{\infty }d\tau
_1 \int_{-b}^{b}] f(\tau )\exp \{ - u_N(a+p,t ;\zeta ) \} $
\par $\exp \{ u_N(a+p,\tau ;\zeta ) \} )...)dp_1...dp_n$, \\
hence taking the limit with $b$ tending to the infinity implies,
that the non-commutative iterated (multiple) integral in Formula
6$(1)$ reduces to the principal value of the usual integral by real
variables $(\tau _1,...,\tau _n)$ and $(p_1,...,p_n)$ 6.1$(1)$.

\par {\bf 7. Theorem.} {\it An original $f(t)$ with $f({\bf R}^n)\subset
{\cal A}_r$ over the Cayley-Dickson algebra ${\cal A}_r$ with $1\le
r \in \bf N$ is completely defined by its image $\mbox{
}_NF_u^n(p;\zeta )$  up to values at points of discontinuity, where
the function $u(p,t; \zeta )$ is given by 1$(8,8.1)$ or 2$(1,2,
2.1)$.}
\par {\bf Proof.} Due to Corollary 6.1 the value $f(t)$
at each point $t$ of continuity of $f(t)$ has the expression
throughout $\mbox{ }_NF_u^n(p;\zeta )$ prescribed by Formula
$6.1(1)$. Moreover, values of the original at points of
discontinuity do not influence on the image $\mbox{ }_NF_u^n(p;\zeta
)$, since on each bounded interval in $\bf R$ by each variable $t_j$
a number of points of discontinuity is finite and by our supposition
above the original function $f(t)$ is $\lambda _n$-almost everywhere
on ${\bf R}^n$ continuous.

\par {\bf 8. Theorem.} {\it Suppose that a function $\mbox{ }_NF_u^n(p;\zeta )$ is analytic
by the variable $p\in {\cal A}_r$ in a domain $W := \{ p\in {\cal
A}_r: a_1< Re (p)< a_{-1} \} $, where $2\le r\in \bf N$, $2^{r-1}\le
n\le 2^r-1$, $f({\bf R}^n)\subset {\cal A}_r$, either $u(p,t;\zeta
)=<p,t) + \zeta $ or $u(p,t;\zeta ) := p_0 s_1 + M(p,t;\zeta )+\zeta
_0$ (see \S \S 1 and 2).
 Let $\mbox{ }_NF^n_u(p;\zeta )$ be written in the form $\mbox{
}_NF^n_u(p;\zeta )=\mbox{ }_NF^{n,0}_u(p;\zeta ) + \mbox{
}_NF^{n,1}_u(p;\zeta )$, where $\mbox{ }_NF^{n,0}_u(p;\zeta )$ is
holomorphic by $p$ in the domain $a_1<Re (p)$. Let also $\mbox{
}_NF^{n,1}_u(p;\zeta )$ be holomorphic by $p$ in the domain $Re
(p)<a_{-1}$. Moreover, for each $a>a_1$ and $b<a_{-1}$ there exist
constants $C_a>0$, $C_b>0$ and $\epsilon _a
>0$ and $\epsilon _b>0$ such that
\par $(1)$ $|\mbox{ }_NF^{n,0}_u(p;\zeta )|\le C_a\exp (-\epsilon _a |p|)$
for each $p\in {\cal A}_r$ with $Re (p)\ge a$,
\par $(2)$ $|\mbox{ }_NF^{n,1}_u(p;\zeta )|\le C_b\exp (-\epsilon _b |p|)$
for each $p\in {\cal A}_r$ with $Re (p)\le b$, the integral \par
$(3)$ $\int_{-N_n\infty }^{N_n\infty }...\int_{-N_1\infty
}^{N_1\infty }\mbox{ }_NF_u^{n,k}(w+p;\zeta )dp$ converges
absolutely for $k=0$ and $k=1$ and each $a_1<w<a_{-1}$. Then $\mbox{
}_NF_u^n(w+p;\zeta )$ is the image of the function
$$(4)\quad f(t)=[(2\pi )^{-1}{\tilde N_n}\int_{-N_n\infty
}^{N_n\infty }] (...([(2\pi )^{-1}{\tilde N_1}\int_{-N_1\infty
}^{N_1\infty }] \mbox{ }_NF_u^n(w+p;\zeta )\exp \{ u(w+p,t;\zeta )
\} )...)dp$$  $$= ({\cal F}^n )^{-1} (\mbox{ }_NF_u^n(w+p;\zeta
),u,t;\zeta ).$$}

\par {\bf Proof.} For the function
$\mbox{ }_NF^{n,1}_u(p;\zeta )$ we consider the substitution of the
variable $p=-g$, $-a_{-1}<Re (g)$. Thus the proof reduces to the
consideration of $\mbox{ }_NF^{n,0}_u(w+p;\zeta )$.
\par An integration by $dp$ in the iterated integral
$(4)$ is treated as in \S 6. Take marked values of variables
$p_1,...,p_{j-1},p_{j+1},...,p_n$ and $t_1,...,t_{j-1},
t_{j+1},...,t_n$, where $s_k=s_k(n;\tau )$ for each $k=1,...,n$ (see
\S 6 also). For a given parameter $\zeta ^j := (\zeta _0+\zeta
_jN_j+...+\zeta _nN_n) + (w+p_0) s_{j+1} +
p_{j+1}s_{j+1}N_{j+1}+...+p_ns_nN_n $ for $u(p,\tau ;\zeta )$
prescribed by Formulas 2$(1,2,2.1)$ or $\zeta ^j := (\zeta _0+\zeta
_jN_j +...+\zeta _nN_n) + (w+p_0) s_{j+1} + p_{j+1}\tau
_{j+1}N_{j+1}+...+p_n\tau _n N_n$ for $u(p,t;\zeta )$ given by
1$(8,8.1)$ instead of $\zeta $ and any non-zero Cayley-Dickson
number $\beta \in {\cal A}_r$ we have $\lim_{\tau _j\to \infty }
[\beta \tau _j +\zeta ^j]/[\beta \tau _j + \zeta ]=1$. \par For any
locally $z$-analytic function $g(z)$ in a domain $U$ satisfying
conditions of \S 5 the homotopy theorem for a non-commutative line
integral over ${\cal A}_r$, $2\le r$, is satisfied (see
\cite{ludoyst,ludfov}). In particular if $U$ contains the straight
line $w+{\bf R} N_j$ and the path $ \gamma _j (t_j) :=  \zeta ^j +
t_jN_j$, then $\int_{-N_j\infty }^{N_j\infty } g(z)dz = \int_{\gamma
_j} g(w+z)dz$, when ${\hat g}(z)\to 0$ while $|z|$ tends to the
infinity, since $|\zeta ^j|$ is a finite number (see Lemma 2.23 in
\cite{lutsltjms}). We apply this to the integrand in Formula $(4)$,
since $\mbox{ }_NF_u^n(w+p;\zeta )$ is locally analytic by $p$ in
accordance with Theorem 4 and Conditions $(1,2)$ are satisfied.
\par Then the integral operator $[(2\pi N_j)^{-1}
\int_{-N_j\infty }^{N_j\infty }]$ on the $j$-th step with the help
of Theorems 2.22 and 3.16 \cite{lutsltjms} gives the inversion
formula corresponding to the real parameter $t_j$ for $f(t)$ and to
the Cayley-Dickson variable $p_0N_0+p_jN_j$ which is restricted on
the complex plane ${\bf C}_{N_j} = {\bf R}\oplus {\bf R}N_j$ (see
also Formulas 6$(4,11)$ above). Therefore, an application of this
procedure by $j=1, 2, ..., n$ as in \S 6 implies Formula $(4)$ of
this theorem.
\par Thus there exist originals $f^0$ and $f^1$ for
functions $\mbox{ }_NF^{n,0}_u(p;\zeta )$ and $\mbox{
}_NF^{n,1}_u(p;\zeta )$ with a choice of $w\in \bf R$ in the common
domain $a_1<Re (p)<a_{-1}$. Then $f=f^0+f^1$ is the original for
$\mbox{ }_NF^n_u(p;\zeta )$ due to the distributivity of the
multiplication in the Cayley-Dickson algebra ${\cal A}_r$ leading to
the additivity of the considered integral operator in Formula $(4)$.
\par {\bf 8.1. Corollary.} {\it Let the conditions of Theorem 8 be
satisfied, then $$(1)\quad f(t)=(2\pi )^{-n} \int_{{\bf R}^n} \mbox{
}_NF_u^n(w+p;\zeta )\exp \{ u(w+p,t;\zeta ) \} dp_1...dp_n$$  $$=
({\cal F}^n )^{-1} (\mbox{ }_NF_u^n(w+p;\zeta ),u,t;\zeta ).$$}
\par {\bf Proof.} In accordance with \S \S 6 and 6.1 each non-commutative
integral given by the left algorithm reduces to the principal value
of the usual integral by the corresponding real variable:
$$(2)\quad (2\pi )^{-1}{\tilde N_j}\int_{-N_j\infty }^{N_j\infty }
\mbox{ }_NF_u^n(w+p;\zeta )\exp \{ u(w+p,t;\zeta ) \} d(p_jN_j)$$
$$ = (2\pi )^{-1}\int_{-\infty }^{\infty } \mbox{ }_NF_u^n(w+p;\zeta )\exp
\{ u(w+p,t;\zeta ) \} dp_j$$ for each $j=1,...,n$. Thus Formula
8$(4)$ with the non-commutative iterated (multiple) integral reduces
to Formula 8.1$(1)$ with the principal value of the usual integral
by real variables $(p_1,...,p_n)$.

\par {\bf 9. Note.} In Theorem 8 Conditions $(1,2)$ can be replaced on
\par $(1)$ $\lim_{n\to \infty }\sup_{p\in C_{R(n)}} \| {\hat F}(p) \|
=0,$ \\
where $C_{R(n)} := \{ z\in {\cal A}_r: |z| =R(n), a_1<Re (z)<a_{-1}
\} $ is a sequence of intersections of spheres with a domain $W$,
where $R(n)<R(n+1)$ for each $n$, $\lim_{n\to \infty } R(n)=\infty
$. Indeed, this condition leads to the accomplishment of the ${\cal
A}_r$ analog of the Jordan Lemma for each $r\ge 2$ (see also Lemma
2.23 and Remark 2.24 \cite{lutsltjms}).
\par Subsequent properties of quaternion, octonion and
general ${\cal A}_r$ multiparameter non-commutative analogs of the
Laplace transform are considered below. We denote by \par $(2)$ $W_f
= \{ p\in {\cal A}_r: ~ a_1(f)< Re (p) <a_{-1} (f) \} $ a domain of
$\mbox{ }_NF_u^n(p;\zeta )$ by the $p$ variable, where $a_1=a_1(f)$
and $a_{-1} = a_{-1}(f)$ are as in \S 1. For an original
\par $(3)$ $f(t)\chi _{U_{1,...,1}}(t)$ we put $W_f = \{ p\in {\cal A}_r: ~
a_1(f) <Re (p) \} ,$ \\ that is $a_{-1} = \infty $. Cases may be,
when either the left hyperplane $Re (p)=a_1$ or the right hyperplane
$Re (p)=a_{-1}$ is (or both are) included in $W_f$. It may also
happen that a domain reduces to the hyperplane $W_f = \{ p: ~ Re
(p)=a_1=a_{-1} \} $.

\par {\bf 10. Proposition.} {\it If images $\mbox{
}_NF_u^n(p;\zeta )$ and $\mbox{ }_NG_u^n(p;\zeta )$  \\
of functions-originals $f(t)$ and $g(t)$ exist in domains $W_f$ and
$W_g$ with values in ${\cal A}_r$, where the function $u(p,t; \zeta
)$ is given by 1$(8,8.1)$ or 2$(1,2, 2.1)$, then for each $\alpha ,
\beta \in {\cal A}_r$ in the case ${\cal A}_2=\bf H$; as well as $f$
and $g$ with values in $\bf R$ and each $\alpha , \beta \in {\cal
A}_r$ or $f$ and $g$ with values in ${\cal A}_r$ and each $\alpha ,
\beta \in \bf R$ in the case of ${\cal A}_r$ with $r\ge 3$; the
function $\alpha \mbox{ }_NF_u(p;\zeta ) + \beta \mbox{
}_NG_u(p;\zeta ) $ is the image of the function $\alpha f(t) +\beta
g(t)$ in a domain $W_f\cap W_g$.}
\par {\bf Proof.} Since the transforms
$\mbox{ }_NF_u^n(p;\zeta )$ and $\mbox{ }_NG_u^n(p;\zeta )$ exist,
then the integral  $$\int_{{\bf R}^n}(\alpha f(t)+\beta g(t)) \exp
(-u(p,t;\zeta ))dt= \int_{{\bf R}^n}\alpha f(t) \exp (-u(p,t;\zeta
))dt$$   $$ + \int_{{\bf R}^n}\beta g(t) \exp (-u(p,t;\zeta ))dt $$
converges in the domain \par $W_f\cap W_g = \{ p\in {\cal A}_r: ~
\max (a_1(f), a_1(g)) < Re (p) <\min (a_{-1}(f),a_{-1}(g)) \} $. \\
We have $t\in {\bf R}^n$, $2^{r-1}\le n\le 2^r-1$, while $\bf R$ is
the center of the Cayley-Dickson algebra ${\cal A}_r$. The
quaternion skew field $\bf H$ is associative. Thus, under the
imposed conditions the constants $\alpha , \beta $ can be carried
out outside integrals.

\par {\bf 11. Theorem.} {\it Let $\alpha =const >0$, let also $F^n(p;\zeta )$
be an image of an original function $f(t)$ with either $u=<p,t) +
\zeta $ or $u$ given by Formulas 2$(1,2)$ over the Cayley-Dickson
algebra ${\cal A}_r$ with $2\le r<\infty $, $2^{r-1}\le n \le
2^r-1$. Then an image $F^n(p/ \alpha ;\zeta )/ \alpha ^n$ of the
function $f(\alpha t)$ exists.}
\par {\bf Proof.} Since $p_js_j+\zeta _j=p_j({s'}_j/\alpha )+
\zeta _j = (p_j/\alpha ) {s'}_j +\zeta _j$ for each $j=1,...,n$,
where $s_j\alpha ={s'}_j$, $s_j=s_j(n;t)$, ${s'}_j = s_j(n;\tau )$,
$\tau _j= \alpha t_j$ for each $j=1,...,n$. Then changing of these
variables implies:
\par $\int_{{\bf R}^n} f(\alpha t) e^{-u(p,t;\zeta )}dt= \int_{{\bf
R}^n}f(\tau )e^{-u(p,\tau /\alpha ;\zeta )} d\tau /\alpha ^n=
F^n(p/\alpha ;\zeta )/\alpha ^n$ \\ due to the fact that the real
filed $\bf R$ is the center $Z({\cal A}_r)$ of the Cayley-Dickson
algebra ${\cal A}_r$.

\par {\bf 12. Theorem.} {\it Let $f(t)$ be a function-original on
the domain $U_{1,...,1}$ such that $\partial f(t)/\partial t_k$ also
for $k=j-1$ and $k=j$ satisfies Conditions 1$(1-4)$. Suppose that
$u(p,t;\zeta )$ is given by 2$(1,2, 2.1)$ or 1$(8,8.1)$ over the
Cayley-Dickson algebra ${\cal A}_r$ with $2\le r<\infty $,
$2^{r-1}\le n \le 2^r-1$. Then
$$(1)\quad {\cal F}^n((\partial f(t)/\partial t_j) \chi
_{U_{1,...,1}}(t),u; p;\zeta ) = - {\cal F}^{n-1; t^j} (f(t)\chi
_{U_{1,...,1}}(t^j), u(p,t^j;\zeta );p;\zeta )$$  $$ +  [p_0 +
\sum_{k=1}^j  p_k {\sf S}_{e_k}] {\cal F}^n (f(t)\chi
_{U_{1,...,1}}(t),u; p;\zeta )$$ in the ${\cal A}_r$ spherical
coordinates or
$$(1.1)\quad {\cal F}^n((\partial f(t)/\partial t_j) \chi
_{U_{1,...,1}}(t),u; p;\zeta ) = - {\cal F}^{n-1; t^j} (f(t)\chi
_{U_{1,...,1}}(t^j), u(p,t^j;\zeta );p;\zeta )$$  $$ + [p_0 + p_j
{\sf S}_{e_j}] {\cal F}^n (f(t)\chi _{U_{1,...,1}}(t),u; p;\zeta )$$
in the ${\cal A}_r$ Cartesian coordinates in a domain $W= \{ p\in
{\cal A}_r: ~ \max ( a_1(f), a_1(\partial f/\partial t_j)) <Re (p)
\} $, where $t^j := (t_1,...,t_j,...,t_n: ~ t_j=0)$, ${\sf S}_{e_k}=
-
\partial /\partial \zeta _k$ for each $k\ge 1$.}
\par {\bf Proof.}  Certainly,
\par $(2)$ $\partial f(t(s))/\partial s_1=\partial f(t)/\partial t_1$
and
\par $(2.1)$ $\partial f(t)/\partial t_j =
\sum_{k=1}^n (\partial f(t(s))/\partial s_k)(\partial s_k/\partial
t_j)=  \sum_{k=1}^j \partial f(t(s))/\partial s_k$
\\ for each $j=2,...,n$, since $t_j=s_j-s_{j+1}$, $t_1 = s_1 - s_2$,
where $s_j=s_j(n;t)$, $s_{n+l}=0$ for each $l\ge 1$. From Formulas
30$(6,7)$ \cite{lutsltjms} we have the equality in the ${\cal A}_r$
spherical coordinates:
\par $(3)$  $\partial \exp (-u(p,t;\zeta ))/\partial s_j = - p_0
\delta _{1,j} \exp (-u(p,t;\zeta )) - p_j{\sf S}_{e_j}\exp
(-u(p,t;\zeta )) $, since
\par $\exp (-u(p,t;\zeta ))= \exp \{ - p_0 s_1 - \zeta _0 \} \exp (-
M(p,t;\zeta ))$,
\par $\partial \exp (-p_0s_1-\zeta _0)/\partial s_j = -p_0
\delta _{1,j} \exp (-p_0s_1-\zeta _0)$, \par  $\partial [\cos
(p_js_j+\zeta _j)-\sin (p_js_j+\zeta _j)i_j]/\partial s_j =\partial
\exp (-(p_js_j+\zeta _j)i_j)/
\partial s_j$  \\  $ = -p_ji_j\exp (-(p_js_j+\zeta _j)i_j)=
-p_j\exp (-(p_js_j+\zeta _j-\pi /2)i_j)$ \\  $=-p_j[\cos
(p_js_j+\zeta _j-\pi /2) - \sin (p_js_j+\zeta _j-\pi /2)i_j]=$ \\
$ - p_j {\sf S}_{e_j} [\cos (p_js_j+\zeta _j) - \sin (p_js_j+\zeta
_j)i_j],$ \\ since $s_j$ and $s_k$ are real independent variables
for each $k\ne j$, where $\delta _{j,k}=0$ for $j\ne k$, while
$\delta _{j,j}=1$, \par $(3.1)$ ${\sf S}_{e_j} [\cos (p_js_j+\zeta
_j) - \sin (p_js_j+\zeta _j)i_j]=$ \par  $ - \partial [\cos
(p_js_j+\zeta _j) - \sin (p_js_j+\zeta _j)i_j]/\partial \zeta _j$
\par $ = [\cos (p_js_j+\zeta _j-\pi /2) - \sin (p_js_j+\zeta _j-\pi
/2)i_j]$. \par In the ${\cal A}_r$ Cartesian coordinates we take
$t_j$ instead of $s_j$ in $(3.1)$. If $\phi (z)$ is a differentiable
function by $z_j$ for each $j$, $\phi : {\cal A}_r\to {\cal A}_r$,
$z_j=p_jt_j+\zeta _j$, then
\par $(3.2)$ $\partial \exp ( - \phi (z))/\partial (qt_j) =
- q [d\exp (\xi )/d\xi ]|_{\xi = -\phi }.(\partial \phi (z)/\partial
z_j)p_j$
\par $= - q p_j [ \sum_{n=1}^{\infty } \sum_{k=1}^{n-1} ((\xi (z))^k
(\partial \phi (z)/\partial z_j)) (\xi (z))^{n-1-k}/n!]|_{\xi =
-\phi }$
\par $ = - q p_j ( - \partial \exp (- \phi (z))/\partial \zeta _j) =
- p_j {\sf S}_{qe_j}\exp (- \phi (z)),$ where either $q=1$ or
$q=-1$, since $\partial z_j/\partial \zeta _j=1$.
\\  That is
\par $(3.3)$ ${\sf S}_{e_j}^x \exp ( - i_k (\phi _k+\zeta _k)) =0$ for
each $j\ne k\ge 1$ and any positive number $x>0$,
\par $(3.4)$ ${\sf S}_{e_j}^x \exp ( - i_j (\phi _j+\zeta _j)) =
\exp ( - i_j (\phi _j+\zeta _j - x \pi /2))$ and \par ${\sf
S}_{-e_j}^x \exp ( - i_j (\phi _j+\zeta _j)) = \exp ( - i_j (\phi
_j+\zeta _j + x \pi /2))$ \\ for each non-negative real number $x\ge
0$, $\phi _k$ and $\zeta _k\in {\bf R}$, where ${\sf S}_{e_j} = {\sf
S}_{e_j}(\zeta _j)$, the zero power ${\sf S}_{e_j}^0=I$ is the unit
operator;
\par $(3.5)$ ${\sf S}_{qe_j} e^{-u(p,t;\zeta )} = e^{-p_0s_1-\zeta _0} $ \\  $ T_j^q [i_0 \delta
_{j,1} \cos (p_1s_1+\zeta _1) + (1-\delta _{j,1}) i_{j-1}\sin
(p_1s_1+\zeta _1)...\cos (p_js_j+\zeta _j) + \{ \sum_{k=j}^{2^r-2}
i_k \sin (p_1s_1+\zeta _1)... \cos (p_{k+1}s_{k+1}+\zeta _{k+1}) \}
+ i_{2^r-1}\sin (p_1s_1+\zeta _1)...\sin (p_{2^r-1}s_{2^r-1} + \zeta
_{2^r-1})]$ \\ in the ${\cal A}_r$ spherical coordinates, where
either $q=1$ or $q=-1$ and \par $(3.6)$ $T_j^x\xi (\zeta _j) := \xi
(\zeta _j-x\pi /2)$
\\ for any function $\xi (\zeta _j)$ and any real number $x\in {\bf
R}$, where $j\ge 1$. Then in accordance with Formula $(3.2)$ we
have:
\par $(3.7)$ ${\sf S}_{qe_j} \exp ( - u(p,t;\zeta )) =$
\par $= [ \sum_{n=1}^{\infty } \sum_{k=1}^{n-1} ((\xi (z))^k
qi_j) (\xi (z))^{n-1-k}/n!]|_{\xi = - u(p,t;\zeta )}$
\\ for $u(p,t;\zeta )$ given by Formulas 1$(8,8.1)$ in the ${\cal A}_r$ Cartesian coordinates, where
either $q=1$ or $q=-1$.
\par The integration by parts theorem (Theorem 2 in \S II.2.6 on p.
228 \cite{kamyn}) states: if $a<b$ and two functions $f$ and $g$ are
Riemann integrable on the segment $[a,b]$, $F(x)= A+ \int_a^xf(t)dt
$ and $G(x)=B+\int_a^xg(t)dt$, where $A$ and $B$ are two real
constants, then $\int_a^b F(x)g(x)dx = F(x)G(x)|^b_a - \int_a^b
f(x)G(x)dx$.
\par Therefore, the integration by parts gives
$$(4)\quad \int_0^{\infty }(\partial f(t)/\partial t_j) \exp
(-u(p,t;\zeta ))dt_j = f(t) \exp (-u(p,t;\zeta
))|_{t_j=0}^{t_j=\infty }$$  $$ - \int_0^{\infty }[f(t)(\partial
\exp (-u(p,t;\zeta ))/\partial t_j)]dt_j .$$ Using the change of
variables $t\mapsto s$ with the unit Jacobian $\partial
(t_1,...,t_n)/\partial (s_1,...,s_n)$ and applying the Fubini's
theorem componentwise to $f_ji_j$ we infer:
$$(5)\quad \int_{U_{1,...,1}}(\partial f(t)/\partial t_j)
\exp (-u(p,t;\zeta ))dt= \int_{s_1\ge s_2\ge ...\ge s_n\ge 0}
(\partial f(t)/\partial t_j) \exp (-u(p,t;\zeta ))ds$$
$$= \int_0^{\infty }...\int_0^{\infty } [\int_{s_{j+1}}^{\infty }(\partial
f(t)/\partial t_j) \exp (-u(p,t;\zeta ))ds_j]dt^j$$  $$ = -
[\int_0^{\infty }...\int_0^{\infty } f(t^j) \exp (-u(p,t^j;\zeta ))
dt^j]$$  $$ + [p_0 + \sum_{k=1}^j p_k {\sf S}_{e_k}] \int_0^{\infty
}...\int_0^{\infty } f(t) \exp (-u(p,t;\zeta ))dt$$ in the ${\cal
A}_r$ spherical coordinates, or
$$(5.1)\quad \int_{U_{1,...,1}}(\partial f(t)/\partial t_j)
\exp (-u(p,t;\zeta ))dt$$
$$ = - [\int_0^{\infty }...\int_0^{\infty } f(t^j) \exp (-u(p,t^j;\zeta ))
dt^j]$$  $$ + [p_0 + p_j {\sf S}_{e_j}] \int_0^{\infty
}...\int_0^{\infty } f(t) \exp (-u(p,t;\zeta ))dt$$ in the ${\cal
A}_r$ Cartesian coordinates, since $\partial \exp ( - (p_0s_1+\zeta
_0))/\partial t_j = - p_0\exp ( - (p_0s_1+\zeta _0))$ for each $1\le
j\le n$. This gives Formula $(1)$, where
$$(6)\quad {\cal F}^{n-1; t^j}(f(t^j)\chi _{U_{1,...,1}},
u(p,t^j;\zeta );p;\zeta ) = \int_0^{\infty }...\int_0^{\infty }
f(t^j) \exp (-u(p,t^j;\zeta )) dt^j$$  $$ = \int_0^{\infty
}dt_1...\int_0^{\infty }dt_{j-1} \int_0^{\infty
}dt_{j+1}...\int_0^{\infty } (dt_n) f(t^j) \exp
(-u(p,t^j;\zeta ))$$  \\
is the non-commutative transform by
$t^j=(t_1,...,t_{j-1},0,t_{j+1},...,t_n)$.
\par {\bf 12.1. Remark.}
Shift operators of the form $\xi (x+\phi ) = \exp (\phi d/dx)\xi
(x)$ in real variables are also frequently used in the class of
infinite differentiable functions with converging Taylor series
expansion in the corresponding domain.
\par It is possible to use also the following
convention. One can put $\cos (\phi _1+\zeta _1) = \cos (\phi
_1+\zeta _1) \cos (\psi _2)...\cos (\psi _{2^r-1})$,...,$\sin (\phi
_1+\zeta _1)... \cos (\phi _k+\zeta _k)=\sin (\phi _1+\zeta _1)...
\cos (\phi _k+\zeta _k)\cos (\psi _{k+1})...\cos (\psi _{2^r-1})$,
where $\psi _j=0$ for each $j\ge 1$, $2\le k <2^r-1$, so that
$T_j^l\cos (\phi _1+\zeta _1) =0$ for each $j>1$ and $l\ge 1$,
$T_j^l\sin (\phi _1+\zeta _1)... \cos (\phi _k+\zeta _k)=0$ for each
$j>k$ and $l\ge 1$, where $T_j^l\xi = T_j^{l-1}(T_j\xi )$ is the
iterated composition for $l>1$, $l\in {\bf N}$. Then
$T_j^le^{-u(p,t;\zeta )}$ gives with such convention the same result
as ${\sf S}_{e_j}^le^{-u(p,t;\zeta )}$, so one can use the symbolic
notation $T_j^le^{-u(p,t;\zeta )}=e^{-u(p,t;\zeta -i_j\pi l/2)}$.
But to avoid misunderstanding we shall use ${\sf S}_{e_j}$ and $T_j$
in the sense of Formulas 12$(3.1-3.7)$.
\par It is worth to mention that instead of 12$(3.7)$ also the
formulas
\par $(1)$ $\exp (p_1i_1+...+p_ni_n) = \cos (\phi ) + M\sin (\phi )$
with $\phi := \phi (p) := [p_1^2+...+p_n^2]^{1/2}$ and
$M=(p_1i_1+...+p_ni_n)/\phi $ for $\phi \ne 0$, $e^0=1$;
\par $(2)$ $\partial \exp (p_1 i_1+..+p_ni_n)/\partial p_j =
[ - \sin (\phi ) + M \cos (\phi ) ]p_j/\phi + (\phi i_j - Mp_j) \phi
^{-2} \sin (\phi )$ and $\partial (p_jt_j+\zeta _j)/\partial \zeta
_j=1$ can be used.

\par {\bf 13. Theorem.} {\it Let $f(t)$ be a function-original.
Suppose that $u(p,t;\zeta )$ is given by 2$(1,2, 2.1)$ or 1$(8,8.1)$
over the Cayley-Dickson algebra ${\cal A}_r$ with $2\le r<\infty $.
Then a (super)derivative of an image is given by the following
formula:
\par $(1)$  $(\partial {\cal F}^n(f(t),u; p;\zeta )/
\partial p).h =
- {\cal F}^n(f(t)s_1,u;p;\zeta )h_0 $ \\
$- {\sf S}_{e_1} {\cal F}^n(f(t)s_1,u; p;\zeta )h_1 -...- {\sf
S}_{e_n} {\cal F}^n(f(t)s_n,u; p;
\zeta )h_n$\\
in the ${\cal A}_r$ spherical coordinates, or
\par $(1.1)$  $(\partial {\cal F}^n(f(t),u; p;\zeta )/
\partial p).h =
- {\cal F}^n(f(t)s_1,u;p;\zeta )h_0 $ \\
$- {\sf S}_{e_1} {\cal F}^n(f(t)t_1,u; p;\zeta )h_1 -...- {\sf
S}_{e_n} {\cal F}^n(f(t)t_n,u; p;
\zeta )h_n$\\
in the ${\cal A}_r$ Cartesian coordinates for each
$h=h_0i_0+...+h_ni_n\in {\cal A}_r$, where $h_0,...,h_n\in {\bf R}$,
$2^{r-1}\le n \le 2^r-1$, $p\in W_f$.}
\par {\bf Proof.} The inequalities
$a_1(f)<Re (p)<a_{-1}(f)$ are equivalent to the inequalities
$a_1(f(t)|t|)<Re (p)<a_{-1}(f(t)|t|)$, since $\lim_{|t|\to + \infty
}\exp (-b|t|)|t|=0$ for each $b>0$. An image ${\cal F}^n(f(t),u;
p;\zeta )$ is a holomorphic function by $p$ for $a_1(f)<Re
(p)<a_{-1}(f)$ by Theorem 4, also $|\int_0^{\infty
}e^{-ct}t^ndt|<\infty $ for each $c>0$ and $n=0,1,2,...$. Thus it is
possible to differentiate under the sign of the integral:
$$(2)\quad (\partial (\int_{{\bf R}^n} f(t)\exp (-u(p,t;\zeta ))dt)/\partial
p).h =$$
$$\sum_{v\in \{ -1, 1 \} ^n }
(\partial (\int_{U_v} f(t)\exp (-u(p,t;\zeta )) \chi _{U_v}
dt)/\partial p).h =$$  $$=
 \int_{{\bf R}^n}f(t)(\partial \exp (-u(p,t;\zeta
))/\partial p).hdt .$$
 Due to Formulas 12$(3, 3.2)$ we get:
\par $(3)$  $(\partial \exp (-u(p,t;\zeta ))/\partial p).h= - \exp
(-u(p,t;\zeta ))s_1h_0 - {\sf S}_{e_1}\exp (-u(p,t;\zeta ))s_1h_1
-... -{\sf S}_{e_n}\exp (-u(p,t;\zeta ))s_nh_n$ \\ in the ${\cal
A}_r$ spherical coordinates, or
\par $(4)$  $(\partial \exp (-u(p,t;\zeta ))/\partial p).h= - \exp
(-u(p,t;\zeta ))s_1h_0 - {\sf S}_{e_1}\exp (-u(p,t;\zeta ))t_1h_1
-... -{\sf S}_{e_n}\exp (-u(p,t;\zeta ))t_nh_n$ \\ in the ${\cal
A}_r$ Cartesian coordinates.\\
Thus from Formulas $(2,3)$ we deduce Formula $(1)$.
\par {\bf 14. Theorem.} {\it If $f(t)$ is a function-original, then
\par $(1)$  ${\cal F}^n(f(t-\tau ),u;p;\zeta )=
{\cal F}^n(f(t),u;p; \zeta + <p,\tau ])$ for either \par $(i)$
$u(p,t;\zeta )= p_0s_1
 + M(p,t;\zeta )+ \zeta _0$ or \par $(ii)$ $u(p,t;\zeta )= <p,t) + \zeta $ over
${\cal A}_r$ with $2\le r<\infty $ in a domain $p\in W_f$, where
$\tau \in {\bf R}^n$, $2^{r-1}\le n\le 2^r-1$, \par $(2)$ $<p,\tau ]
= p_0s_1+p_1s_1i_1+...+p_ns_ni_n$ with $s_j=s_j(n;\tau )$ for each
$j$ in the first $(i)$ and $<p,\tau ]=<p,\tau )$ in the second
$(ii)$ case (see also Formulas 1$(8)$, 2$(1,2,2.1)$).}

\par {\bf Proof.} For $p$  in the domain $Re (p)>a_1$
the identities are satisfied: $$(3)\quad {\cal F}^n((f\chi
_{U_{1,...,1}}) (t-\tau ),u;p;\zeta ) = \int_{\tau _1}^{\infty
}...\int_{\tau _n}^{\infty }f(t-\tau ) e^{-u(p,t;\zeta )}dt$$  $$ =
\int_{U_{1,...,1}} f(t) e^{-u(p,\xi ;\zeta + <p,\tau ])}d\xi ={\cal
F}^n ((f\chi _{U_{1,...,1}})(t),u;p;\zeta +<p,\tau ]), $$  due to
Formulas 1$(7,8)$ and 2$(1,2,2.1,4)$, since $p_0s_1(n;t) + \zeta _0=
p_0 s_1(n;\xi ) + \zeta _0 + p_0s_1(n;\tau )$ and $p_jt_j + \zeta _j
= p_j\xi _j + (\zeta _j+p_j \tau _j)$ and $p_js_j(n;t) + \zeta _j =
p_js_j(n;\xi ) + (\zeta _j+p_j s_j(n;\tau ))$ for each $j=
1,...,2^r-1$, where $t=\xi +\tau $. Symmetrically we get $(2)$ for
$U_v$ instead of $U_{1,...,1}$. Naturally, that the multiparameter
non-commutative Laplace integral for an original $f$ can be
considered as the sum of $2^n$ integrals by the sub-domains $U_v$:
$$(4)\quad \int_{{\bf R}^n} f(t) \exp (-u(p,t;\zeta ))dt= \sum_{v\in
\{ -1 , 1 \} ^n } \int_{{\bf R}^n}f(t)\exp (-u(p,t;\zeta ))\chi
_{U_v}(t) dt .$$  The summation by all possible $v\in \{ -1 , 1 \}
^n$ gives Formula $(1)$.

\par {\bf 15. Note.} In view of the definition of the non-commutative
transform ${\cal F}^n$ and $u(p,t;\zeta )$ and Theorem 14 the term
$\zeta _1 i_1+...+\zeta _{2^r-1}i_{2^r-1}$ has the natural
interpretation as the initial phase of a retardation.

\par {\bf 16. Theorem.} {\it If $f(t)$ is a function-original
with values in ${\cal A}_r$ for $2\le r<\infty $, $2^{r-1}\le n \le
2^r-1$, $b\in \bf R$, then
$$(1)\quad {\cal F}^n(e^{b(t_1+...+t_n)}f(t), u;p;\zeta )= {\cal
F}^n(f(t), u;p-b;\zeta )$$ for each $a_{-1}+b > Re (p)
>a_1+b$, where $u$ is given by 1$(8, 8.1)$ or 2$(1,2)$.}

\par {\bf Proof.} In accordance with Expressions 1$(8,8.1)$ and 2$(1,2,2.1)$
one has $u(p,t;\zeta ) - b (t_1+...+t_n) = u(p-b,t;\zeta )$. If
$a_{-1}+b>Re (p)>a_1+b$, then the integral
$$(2)\quad {\cal F}^n(e^{b(t_1+...+t_n)}f(t)\chi
_{U_v}(t), u;p;\zeta )= \int_{U_v} f(t)e^{b(t_1+...+t_n)}\exp (-
u(p,t;\zeta ))dt$$  $$= \int_{U_v} f(t)\exp (- u(p-b,t;\zeta ) )dt =
{\cal F}^n(f(t)\chi _{U_v}(t), u;p-b;\zeta )$$ converges. Applying
Decomposition 14$(4)$ we deduce Formula $(1)$.

\par {\bf 17. Theorem.} {\it  Let a function $f(t)$ be a real valued
original, \\ $F(p;\zeta ) = {\cal F}^n(f(t);u;p;\zeta )$, where the
function $u(p,t;\zeta )$ is given by 1$(8,8.1)$ or 2$(1,2,2.1)$. Let
also $G(p;\zeta )$ and $q(p)$ be locally analytic functions such
that
\par $(1)$ ${\cal F}^n(g(t,\tau );u;p;\zeta ) = G(p;\zeta )
\exp (-u(q(p),\tau ;\zeta ))$ \\  for $u=<p,t)+ \zeta $ or
$u=p_0(t_1+...+t_n)+ M(p,t;\zeta ) + \zeta _0$, then
\par $(2)$ ${\cal F}^n (\int_{{\bf R}^n}
g(t,\tau )f(\tau )d\tau ;u;p;\zeta ) = G(p;\zeta )F(q(p);\zeta )$ \\
for each $p\in W_g$ and $q(p) \in W_f$, where $2\le r<\infty $,
$2^{r-1}\le n \le 2^r-1$.}
\par {\bf Proof.} If $p \in W_g$ and $q(p)\in W_f$,
then in view of the Fubini's theorem and the theorem conditions a
change of an integration order gives the equalities:
$$\int_{{\bf R}^n} (\int_{{\bf R}^n} g(t,\tau ) f(\tau )d\tau
)\exp (-u(p,t;\zeta ))dt $$   $$= \int_{{\bf R}^n} (\int_{{\bf R}^n}
g(t,\tau ) \exp (-u(p,t;\zeta ))dt )f(\tau )d\tau $$
$$= \int_{{\bf R}^n} G(p;\zeta ) \exp (- u(q(p),\tau ;\zeta ))
f(\tau )d\tau $$
$$ =G(p;\zeta )\int_{{\bf R}^n}f(\tau )\exp (-u(q(p),\tau ;\zeta )) d\tau =
G(p;\zeta )F(q(p);\zeta ) ,$$  since $t, \tau \in {\bf R}^n$ and the
center of the algebra ${\cal A}_r$ is $\bf R$.
\par {\bf 18. Theorem.} {\it  If a function
$f(t)\chi _{U_{1,...,1}}$ is original together with its derivative
$\partial ^nf(t)\chi _{U_{1,...,1}}(t)/\partial s_1...\partial s_n$
or $\partial ^nf(t)\chi _{U_{1,...,1}}(t)/\partial t_1...\partial
t_n$, where $F^n_u(p;\zeta )$ is an image function of $f(t)\chi
_{U_{1,...,1}}$ over the Cayley-Dickson algebra ${\cal A}_r$ with
$2\le r\in \bf N$, $2^{r-1}\le n \le 2^r-1$, for $u=p_0s_1 +
M(p,t;\zeta )+ \zeta _0$ given by 2$(1,2,2.1)$, then
$$(1)\quad \lim_{p\to \infty } \{ [p_0 + p_1{\sf S}_{e_1}] p_2{\sf S}_{e_2}...p_n
{\sf S}_{e_n} F^n_u(p;\zeta ) + \sum_{m=0}^{n-1} (-1)^m $$
$$ \sum_{1 \le j_1 <...< j_{n-m} \le n; ~ 1\le l_1<...<l_m\le n; ~
l_{\alpha }\ne j_{\beta } ~ \forall \alpha , \beta } [p_0\delta
_{1,j_1} + p_{j_1}{\sf S}_{e_{j_1}}] p_{j_2}{\sf
S}_{e_{j_2}}...p_{j_{n-m}}{\sf S}_{e_{j_{n-m}}}$$
$$ F^{n-m}_u(p^{(l)}; \zeta ) \} = (-1)^{n+1} f(0)e^{-u(0,0;\zeta
)},$$ or
$$(1.1)\quad \lim_{p\to \infty } \{ [p_0 + p_1{\sf S}_{e_1}]
[p_0+p_2{\sf S}_{e_2}]...[p_0+p_n {\sf S}_{e_n}] F^n_u(p;\zeta ) +
\sum_{m=0}^{n-1} (-1)^m $$
$$ \sum_{1 \le j_1 <...< j_{n-m} \le n; ~ 1\le l_1<...<l_m\le n; ~
l_{\alpha }\ne j_{\beta } ~ \forall \alpha , \beta } [p_0 +
p_{j_1}{\sf S}_{e_{j_1}}] [p_0+p_{j_2}{\sf
S}_{e_{j_2}}]...[p_0+p_{j_{n-m}}{\sf S}_{e_{j_{n-m}}}]$$
$$ F^{n-m}_u(p^{(l)}; \zeta ) \} = (-1)^{n+1} f(0)e^{-u(0,0;\zeta
)}$$ for $u(p,t;\zeta )$ given by 1$(8,8.1)$, where $f(0) =
\lim_{t\in U_{1,...,1}; t\to 0} f(t)$, $p$ tends to the infinity
inside the angle $|Arg (p)|<\pi /2-\delta $ for some $0<\delta <\pi
/2$, $1\le j\le 2^r-1$, $p^{(l)} = \sum_{j=0, j\notin (l)}^n
p_ji_j$, $(l) = (l_1,...,l_m)$. If the restriction
\par $f(t)|_{t_{j_1}=0,...,t_{j_m}=0; t_k=\infty \forall k\notin \{
j_1,...,j_m \}} = \lim_{t\in U_{1,...,1}; t_{j_1}\to
0,...,t_{j_m}\to 0; t_k\to \infty ~ \forall k\notin \{ j_1,...,j_m
\} } f(t)$ exists for all $1\le j_1<...<j_m\le n$, then
$$(2)\quad \lim_{p\to 0} \{ [p_0 + p_1{\sf S}_{e_1}] p_2{\sf S}_{e_2}...p_n{\sf S}_{e_n}
F^n_u(p;\zeta ) + \sum_{m=0}^{n-1} (-1)^m
$$  $$ \sum_{1 \le
j_1 <...< j_{n-m} \le n; ~ 1\le l_1<...<l_m\le n; ~ l_{\alpha }\ne
j_{\beta } ~ \forall \alpha , \beta } [p_0\delta _{1,j_1} +
p_{j_1}{\sf S}_{e_{j_1}}] p_{j_2}{\sf S}_{e_{j_2}}...p_{j_{n-m}}{\sf
S}_{e_{j_{n-m}}}$$  $$ F^{n-m}_u(p^{(l)}; \zeta ) \} $$
$$ =\sum_{m=0}^{n-1} (-1)^m \sum_{1\le j_1<...<j_m\le n}
f(t)|_{t_{j_1}=0,...,t_{j_m}=0; t_k=\infty \forall k\notin \{
j_1,...,j_m \}} e^{-u(0,0,\zeta )}$$ in the ${\cal A}_r$ spherical
coordinates or
$$(2.1)\quad \lim_{p\to 0} \{ [p_0 + p_1{\sf S}_{e_1}] [p_0+p_2{\sf S}_{e_2}]...[p_0+p_n{\sf
S}_{e_n}] F^n_u(p;\zeta ) + \sum_{m=0}^{n-1} (-1)^m
$$  $$ \sum_{1 \le
j_1 <...< j_{n-m} \le n; ~ 1\le l_1<...<l_m\le n; ~ l_{\alpha }\ne
j_{\beta } ~ \forall \alpha , \beta } [p_0 + p_{j_1}{\sf
S}_{e_{j_1}}] [p_0+p_{j_2}{\sf S}_{e_{j_2}}]...[p_0+p_{j_{n-m}}{\sf
S}_{e_{j_{n-m}}}]$$  $$ F^{n-m}_u(p^{(l)}; \zeta ) \} $$
$$ =\sum_{m=0}^{n-1} (-1)^m \sum_{1\le j_1<...<j_m\le n}
f(t)|_{t_{j_1}=0,...,t_{j_m}=0; t_k=\infty \forall k\notin \{
j_1,...,j_m \}} e^{-u(0,0,\zeta )}$$ in the ${\cal A}_r$ Cartesian
coordinates, where $p\to 0$ inside the same angle.}
\par {\bf Proof.} In accordance with Theorem 12 the equality follows:
$$(3) \quad {\cal F}^n((\partial f(t)/\partial s_j)\chi
_{U_{1,...,1}}(t),u;p;\zeta ) = [p_0 \delta _{1,j} + p_j {\sf
S}_{e_j}] {\cal F}^n (f(t)\chi _{U_{1,...,1}}(t), u(p,t;\zeta ),
p;\zeta )$$ $$- {\cal F}^{n-1; t^j} (f(t^j)\chi _{U_{1,...,1}},
u(p,t^j;\zeta );p;\zeta )$$ for $u= u(p,t;\zeta ) =
p_0s_1+M(p,t;\zeta )+\zeta _0$ in the ${\cal A}_r$ spherical
coordinates, or
$$(3.1) \quad {\cal F}^n((\partial f(t)/\partial t_j)\chi
_{U_{1,...,1}}(t),u;p;\zeta ) = [p_0 + p_j {\sf S}_{e_j}] {\cal F}^n
(f(t)\chi _{U_{1,...,1}}(t), u(p,t;\zeta ), p;\zeta )$$  $$- {\cal
F}^{n-1; t^j} (f(t^j)\chi _{U_{1,...,1}}, u(p,t^j;\zeta );p;\zeta
)$$ in the ${\cal A}_r$ Cartesian coordinates, since \par $(3.2)$
$\partial f(t(s))/\partial s_j = -
\partial f(t)/\partial t_{j-1} + \partial f(t)/\partial t_j$ for
each $j\ge 2$, $\partial f(t(s))/\partial s_1
=  \partial f(t)/\partial t_1$, \\
where $p=p_0+p_1i_1+...+p_{2^r-1}i_{2^r-1}\in {\cal A}_r$,
$p_0,...,p_{2^r-1}\in {\bf R}$, $ \{ i_0,...,i_{2^r-1} \} $ are the
generators of the Cayley-Dickson algebra ${\cal A}_r$, $s_{n+l}=0$
for each $l\ge 1$, the zero power ${\sf S}_{e_j}^0=I$ is the unit
operator. For short we write $f$ instead of $f\chi _{U_{1,...,1}}$.
Thus the limit exists:
$$(4)\quad {\cal F}^{n-1;t^j} (f(t^j), u(p,t^j;\zeta );p;\zeta ) =$$
$$\lim_{t_j\to +0} \int_0^{\infty }dt_1...\int_0^{\infty }dt_{j-1}
\int_0^{\infty }dt_{j+1}...\int_0^{\infty } (dt_n) f(t) \exp
(-u(p,t;\zeta )).$$  Mention, that $(...((t^1)^2)...)^j =
(0,...,0,t_j,...,t_n: t_j=0)$ for every $1\le j\le n$, since
$t_k=s_k-s_{k+1}$ for each $1\le k\le n$. We apply these Formulas
$(3,4)$ by induction $j=1,...,n$, $2^{r-1} \le n \le 2^r-1$, to
$\partial ^nf(t)/\partial s_1...\partial s_n$,...,$\partial ^{n-j+1}
f(t)/\partial s_j...\partial s_n$,\\ ...,$\partial f(t)/\partial
s_n$ instead of $\partial f(t)/\partial s_j$.
\par From Note 8 \cite{lutsltjms} it follows, that in the ${\cal A}_r$ spherical coordinates
$$\lim_{p\to
\infty , |Arg (p)|<\pi /2-\delta } {\cal F}^n((\partial
^nf(t)/\partial s_1...\partial s_n)\chi _{U_{1,..,1}},u;p;\zeta
)=0,$$  also in the ${\cal A}_r$ Cartesian coordinates $$\lim_{p\to
\infty , |Arg (p)|<\pi /2-\delta } {\cal F}^n((\partial
^nf(t)/\partial t_1...\partial t_n)\chi _{U_{1,..,1}},u;p;\zeta
)=0,$$ which gives the first statement of this theorem, since
$u(p,0,\zeta ) = u(0,t;\zeta ) = u(0,0,\zeta ) $ and
$F^0_u(p^{(1,...,1)};\zeta ) = f(0) e^{-u(0,0,\zeta )}$, while
$F^n_u(p;\zeta )$ is defined for each $Re (p)>0$.
\par If the limit $f(t^{<j>})$ exists, where $t^{<j>} := (t_1,...,t_j,...,t_n: ~ t_j=\infty )$, then
$$(5)\quad \lim_{t_j\to \infty } \int_0^{\infty
}dt_1...\int_0^{\infty }dt_{j-1} \int_0^{\infty
}dt_{j+1}...\int_0^{\infty } (dt_n) f(t) \exp (-u(p,t;\zeta ))$$
$$ =: {\cal F}^{n-1; <t^j>} (f(t^{<j>}), u(p,t^{<j>};\zeta );p;\zeta ).$$
 Certainly,
$(...((t^{<1>})^{<2>})...)^{<j>} = (t_1, ...,t_n: t_1=\infty
,...,t_j=\infty )$ for each $1\le j\le n$. Therefore, the limit
exists:
$$\lim_{p\to 0, |Arg (p)|<\pi /2 -\delta }\int_{U_{1,...,1}}
(\partial ^nf(t)/\partial s_1...\partial s_n)\exp (-p_0s_1-\zeta _0
-M(p,t;\zeta ))$$   $$= \int_{U_{1,...,1} } (\partial
^nf(t)/\partial s_1...\partial s_n)e^{-u(0,0;\zeta )} dt$$
$$= \sum_{m=0}^n (-1)^m \sum_{1\le j_1<...<j_m\le n}
f(t)|_{t_{j_1}=0,..., t_{j_m}=0; t_k=\infty ~ \forall  k\notin \{
j_1,...,j_m \} } $$  $$=\lim_{p\to 0, |Arg (p)|<\pi /2 -\delta } \{
[p_0 + p_1{\sf S}_{e_1}] p_2{\sf S}_{e_2}...p_n{\sf S}_{e_n}
F^n_u(p;\zeta )$$  $$ + \sum_{m=0}^{n-1} (-1)^m
$$  $$ \sum_{1 \le
j_1 <...< j_{n-m} \le n; ~ 1\le l_1<...<l_m\le n; ~ l_{\alpha }\ne
j_{\beta } ~ \forall \alpha , \beta } [p_0\delta _{1,j_1} +
p_{j_1}{\sf S}_{e_{j_1}}] p_{j_2}{\sf S}_{e_{j_2}}...p_{j_{n-m}}{\sf
S}_{e_{j_{n-m}}}$$
$$ F^{n-m}_u(p^{(l)}; \zeta )
 + (-1)^n f(0) e^{ -u(0,0,\zeta )} \} ,$$ from which the second
statement of this theorem follows in the ${\cal A}_r$ spherical
coordinates and analogously in the ${\cal A}_r$ Cartesian
coordinates using Formula $(3.1)$.

\par {\bf 19. Definitions.}
\par Let $X$ and $Y$ be two $\bf R$ linear normed spaces which are
also left and right ${\cal A}_r$ modules, where $1\le r$. Let $Y$ be
complete relative to its norm.  We put $X^{\otimes k} := X\otimes
_{\bf R} ... \otimes _{\bf R} X$ is the $k$ times ordered tensor
product over $\bf R$ of $X$. By $L_{q,k}(X^{\otimes k},Y)$ we denote
a family of all continuous $k$ times $\bf R$ poly-linear and ${\cal
A}_r$ additive operators from $X^{\otimes k}$ into $Y$. Then
$L_{q,k}(X^{\otimes k},Y)$ is also a normed $\bf R$ linear and left
and right ${\cal A}_r$ module complete relative to its norm. In
particular, $L_{q,1}(X,Y)$ is denoted also by $L_q(X,Y)$.
\par We present $X$ as the direct sum $X=X_0i_0\oplus ... \oplus
X_{2^r-1} i_{2^r-1}$, where $X_0$,...,$X_{2^r-1}$ are pairwise
isomorphic real normed spaces. If $A\in L_q(X,Y)$ and $A(xb)=(Ax)b$
or $A(bx)=b(Ax)$ for each $x\in X_0$ and $b\in {\cal A}_r$, then an
operator $A$ we call right or left ${\cal A}_r$-linear respectively.
\par An $\bf R$ linear space of left (or right) $k$ times ${\cal A}_r$
poly-linear operators is denoted by $L_{l,k}(X^{\otimes k},Y)$ (or
$L_{r,k}(X^{\otimes k},Y)$ respectively).
\par We consider a space of test function ${\cal D} := {\cal D}({\bf
R}^n,Y)$ consisting of all infinite differentiable functions $f:
{\bf R}^n\to Y$ on ${\bf R}^n$ with compact supports. A sequence of
functions $f_n\in {\cal D}$ tends to zero, if all $f_n$ are zero
outside some compact subset $K$ in the Euclidean space ${\bf R}^n$,
while on it for each $k=0,1,2,...$ the sequence $ \{ f^{(k)}_n: ~
n\in {\bf N} \} $ converges to zero uniformly. Here as usually
$f^{(k)}(t)$ denotes the $k$-th derivative of $f$, which is a $k$
times $\bf R$ poly-linear symmetric operator from $({\bf
R}^n)^{\otimes k}$ to $Y$, that is $f^{(k)}(t).(h_1,...,h_k)=
f^{(k)}(t).(h_{\sigma (1)},...,h_{\sigma (k)})\in Y$ for each
$h_1,...,h_k\in {\bf R}^n$ and every transposition $\sigma : \{
1,...,k \} \to \{ 1,...,k \}$, $\sigma $ is an element of the
symmetric group $S_k$, $t\in {\bf R}^n$. For convenience one puts
$f^{(0)}=f$. In particular, $f^{(k)}(t).(e_{j_1},...,e_{j_k})=
\partial ^kf(t)/\partial t_{j_1}...\partial t_{j_k}$ for all
$1\le j_1,...,j_k\le n$, where $e_j = (0,...,0,1,0,...,0)\in {\bf
R}^n$ with $1$ on the $j$-th place.
\par Such convergence in $\cal D$ defines closed subsets in this space $\cal D$, their
complements by the definition are open, that gives the topology on
$\cal D$. The space ${\cal D}$ is $\bf R$ linear and right and left
${\cal A}_r$ module.

\par By a generalized function of class ${\cal D}' := [{\cal D}({\bf R}^n,Y)]'$ is called a continuous
$\bf R$-linear ${\cal A}_r$-additive function $g: {\cal D} \to {\cal
A}_r$. The set of all such functionals is denoted by ${\cal D}'$.
That is, $g$ is continuous, if for each sequence $f_n\in \cal D$,
converging to zero, a sequence of numbers $g(f_n)=: [g,f_n) \in
{\cal A}_r$ converges to zero for $n$ tending to the infinity. \par
A generalized function $g$ is zero on an open subset $V$ in ${\bf
R}^n$, if $[g,f)=0$ for each $f\in {\cal D}$ equal to zero outside
$V$. By a support of a generalized function $g$ is called the
family, denoted by $supp (g)$, of all points $t\in {\bf R}^n$ such
that in each neighborhood of each point $t\in supp (g)$ the
functional $g$ is different from zero. The addition of generalized
functions $g, h$ is given by the formula: \par $(1)$ $[g+h,f):=
[g,f)+ [h,f)$.
\par The multiplication $g\in {\cal D}'$ on an infinite differentiable
function $w$ is given by the equality:
\par $(2)$ $[gw,f)=[g, wf)$ either for $w: {\bf R}^n\to {\cal A}_r$ and
each test function $f\in \cal D$ with a real image $f({\bf
R}^n)\subset {\bf R}$, where $\bf R$ is embedded into $Y$; or $w:
{\bf R}^n\to {\bf R}$ and $f: {\bf R}^n\to Y$. \par  A generalized
function $g'$ prescribed by the equation: \par $(3)$ $[g',f):= -
[g,f')$ is called a derivative $g'$ of a generalized function $g$,
where $f' \in {\cal D}({\bf R}^n,L_q({\bf R}^n,Y))$, $g'\in [{\cal
D}({\bf R}^n,L_q({\bf R}^n,Y))]'$.
\par Another space ${\cal B} := {\cal B}({\bf R}^n,Y)$ of test functions consists of all
infinite differentiable functions  $f: {\bf R}^n\to Y$ such that the
limit $\lim_{|t|\to +\infty } |t|^m f^{(j)}(t)=0$ exists for each
$m=0,1,2,...$, $j=0,1,2,...$. A sequence $f_n\in \cal B$ is called
converging to zero, if the sequence $|t|^mf_n^{(j)}(t)$ converges to
zero uniformly on ${\bf R}^n\setminus B({\bf R}^n,0,R)$ for each $m,
j=0,1,2,...$ and each $0<R< + \infty $, where $B(Z,z,R) := \{ y\in
Z: ~ \rho (y,z)\le R \} $ denotes a ball with center at $z$ of
radius $R$ in a metric space $Z$ with a metric $\rho $. The family
of all $\bf R$-linear and ${\cal A}_r$-additive functionals on $\cal
B$ is denoted by ${\cal B}'$.
\par In particular we can take $X={\cal A}_r^{\alpha }$, $Y= {\cal
A}_r^{\beta }$ with $1\le \alpha , \beta \in \bf Z$. Analogously
spaces ${\cal D}(U,Y)$, $[{\cal D}(U,Y)]'$, ${\cal B}(U,Y)$ and
$[{\cal B}(U,Y)]'$ are defined for domains $U$ in ${\bf R}^n$, for
example, $U=U_v$ (see also \S 1).

\par A generalized function $f\in {\cal B}'$ we call a
generalized original, if there exist real numbers $a_1<a_{-1}$ such
that for each $a_1 < w_{-1}, w_1,...,w_{-n}, w_n < a_{-1}$ the
generalized function \par $(4)$ $f(t)\exp (-(q_v,t))\chi _{U_v}$ is
in $[{\cal B}(U_v,Y)]'$ for all $v= (v_1,...,v_n)$, $v_j\in \{ -1, 1
\}$ for every $j=1,...,n$ for each $t\in {\bf R}^n$ with $t_j v_j
\ge 0$ for each $j=1,...,n$, where $q_v =
(v_1w_{v_11},...,v_nw_{v_nn})$.
\par By an image of such original we call a function
\par $(5)$ ${\cal F}^n(f,u;p;\zeta ):= [f, \exp (-u(p,t;\zeta )))$
of the variable $p\in {\cal A}_r$ with the parameter $\zeta \in
{\cal A}_r$, defined in the domain $W_f = \{ p\in {\cal A }_r: ~
a_1< Re (p) <a_{-1} \} $ by the following rule. For a given $p\in
W_f$ choose $a_1 < w_1,...,w_n < Re (p) < w_{-1},...,w_{-n} <
a_{-1}$, then
\par $(6)$  $[f,\exp (-u(p,t;\zeta
)) := \sum_v [f\exp (- (q_v,t)) , \exp \{ - [u(p,t;\zeta )- (q_v,t)]
\} \chi _{U_v} )$, \\
since $\exp \{ - [ u(p,t;\zeta ) - (q_v,t) ] \} \in {\cal
B}(U_v,Y)$,  where in each term \\ $[f\exp (- (q_v,t)) , \exp \{ -
[u(p,t;\zeta )- (q_v,t)] \} \chi _{U_v} )$ the generalized function
belongs to $[{\cal B}(U_v,Y)]'$ by Condition $(4)$, while the sum in
$(6)$ is by all admissible vectors $v\in \{ -1, 1 \} ^n$.

\par {\bf 20. Note and Examples.} Evidently the transform
${\cal F}^n(f,u;p;\zeta )$ does not depend on a choice of $\{
w_{-1}, w_1,...,w_{-n}, w_n \} $, since \par $[f\exp (-(q_v,t),\exp
(-[u(p,t;\zeta )- (q_v,t)])\chi _{U_v})=$ \par $[f\exp
(-(q_v,t)-(b_v,t)),\exp (-[u(p,t;\zeta )-(q_v,t)-(b_v,t)])\chi
_{U_v})$ \\ for each $b\in {\bf R}^n$ such that $a_1<w_j + b_j< Re
(p) < w_{-j} + b_{-j} < a_{-1}$ for each $j=1,...,n$, because $\exp
(-(b_v,t))\in \bf R$. At the same time the real field $\bf R$ is the
center of the Cayley-Dickson algebra ${\cal A}_r$, where $2\le r\in
\bf N$.
\par  Let $\delta $ be the Dirac delta function, defined by the
equation \par $(DF)$ $[\delta (t),\phi (t)) := \phi (0)$ for each
$\phi \in {\cal B}$. Then
\par $(1)$ ${\cal F}^n(\delta ^{(j)}(t-\tau ),u;p;\zeta )
=\sum_{v\in \{ -1,1 \} ^n} [\delta ^{(j)}(t-\tau )
\exp (- (q_v,t)), \exp (-[u(p,t;\zeta )- (q_v,t)])\chi _{U_v})$ \\
$= (-1)^j [\partial ^j_t \exp (-[u(p,t;\zeta )])|_{t=\tau }$, \\
since it is possible to take $- \infty <a_1<0<a_{-1}<\infty $ and
$w_k=0$ for each $k\in \{ -1, 1, -2, 2,..., -n, n \} $, where $\tau
\in {\bf R}^n$ is the parameter, $\partial ^j_t := \partial
^{|j|}/\partial t_1^{j_1}...\partial t_1^{j_1}$. In particular, for
$j=0$ we have
\par $(2)$ ${\cal F}^n(\delta (t-\tau
),u;p;\zeta ) = \exp (-u(p,\tau ;\zeta ))$. \\
In the general case:
\par $(3)$ ${\cal F}^n(\partial ^{|j|}\delta (t)/\partial s_1^{j_1}...
\partial s_n^{j_n},u;p;\zeta )=$ \par
$\sum_{0\le k_1\le j_1} {{j_1}\choose {k_1}} p_0^{j_1-k_1}(p_1{\sf
S}_{e_1})^{k_1}(p_2{\sf S}_{e_2})^{j_2}... (p_n{\sf
S}_{e_n})^{j_n}\exp (- \zeta _0 - M(p,0;\zeta ))$ \\ in the ${\cal
A}_r$ spherical coordinates, or
\par $(3.1)$ ${\cal F}^n(\partial ^{|j|}\delta (t)/\partial t_1^{j_1}...
\partial t_n^{j_n},u;p;\zeta )=$ \par
$ (p_0+p_1{\sf S}_{e_1})^{j_1} (p_0+p_2{\sf S}_{e_2})^{j_2}...
(p_0+p_n{\sf S}_{e_n})^{j_n}\exp (- u(p,0;\zeta ))$ \\ in the ${\cal
A}_r$ Cartesian coordinates,  where $j_1+...+j_n=|j|$, $k_1,
j_1,..., j_n$ are nonnegative integers, $2^{r-1}\le n\le 2^r-1$,
${l\choose m} := l!/[m!(l-m)!]$ denotes the binomial coefficient,
$0!=1$, $1!=1$, $2!=2$; $l!=1 \cdot 2\cdot ... \cdot l$ for each
$l\ge 3$, $s_j=s_j(n;t)$.
\par The transform ${\cal F}^n(f)$ of any generalized function $f$ is
the holomorphic function by $p\in W_f$ and by $\zeta \in {\cal
A}_r$, since the right side of Equation 19$(5)$ is holomorphic by
$p$ in $W_f$ and by $\zeta $ in view of Theorem 4. Equation 19$(5)$
implies, that Theorems 11 - 13 are accomplished also for generalized
functions.
\par
For $a_1=a_{-1}$ the region of convergence reduces to the vertical
hyperplane in ${\cal A}_r$ over $\bf R$. For $a_{-1} < a_1$ there is
no any common domain of convergence and $f(t)$ can not be
transformed.

\par {\bf 21. Theorem.} {\it If $f(t)$ is an original function on ${\bf R}^n$,
$F^n(p;\zeta )$ is its image, $\partial ^{|j|} f(t)/\partial
s_1^{j_1}...\partial s_n^{j_n}$ or $\partial ^{|j|} f(t)/\partial
t_1^{j_1}...\partial t_n^{j_n}$ is an original, $|j|=j_1+...+j_n$,
$0\le j_1,...,j_n\in {\bf Z}$, $2^{r-1}\le n \le 2^r-1$; then
$$(1)\quad {\cal F}^n(\partial ^{|j|} f(t)/\partial s_1^{j_1}...\partial
s_n^{j_n}, u; p;\zeta ) = \sum_{0\le k_1\le j_1} $$
$$ {{j_1}\choose {k_1}}
p_0^{j_1-k_1}(p_1{\sf S}_{e_1})^{k_1}(p_2{\sf
S}_{e_2})^{j_2}...(p_n{\sf S}_{e_n})^{j_n} {\cal F}^n (f(t),u;
p;\zeta )$$ for $u(p,t;\zeta ) := p_0s_1 + M(p,t;\zeta )+\zeta _0$
given by 2$(1,2,2.1)$, or
$$(1.1)\quad {\cal F}^n(\partial ^{|j|} f(t)/\partial t_1^{j_1}...\partial
t_n^{j_n}, u; p;\zeta ) = $$
$$ (p_0 + p_1{\sf S}_{e_1})^{j_1}(p_0+p_2{\sf
S}_{e_2})^{j_2}...(p_0+p_n{\sf S}_{e_n})^{j_n} {\cal F}^n (f(t),u;
p;\zeta )$$ for $u(p,t;\zeta )$ given by 1$(8,8.1)$ over the
Cayley-Dickson algebra ${\cal A}_r$ with $2\le r<\infty $. Domains,
where Formulas $(1, 1.1)$ are true may be different from a domain of
the multiparameter noncommutative transform for $f$, but they are
satisfied in the domain $a_1< Re (p)<a_{-1}$, where \par $a_{-1} =
\min (a_{-1}(f),a_{-1}(\partial ^{|m|} f(t)/\partial \phi
_1^{m_1}...\partial \phi _n^{m_n}): |m|\le |j|, 0\le m_l\le
j_l\forall l)$; \par $a_1= \max (a_1(f),a_1(\partial ^{|m|}
f(t)/\partial \phi _1^{m_1}...\partial \phi _n^{m_n}): |m|\le |k|,
0\le m_l\le j_l\forall l)$, if $a_1<a_{-1}$, where $\phi _j=s_j$ or
$\phi _j=t_j$ for each $j$ correspondingly.}

\par {\bf Proof.} To each domain $U_v$ the domain
$U_{-v}$ symmetrically corresponds. The number of different vectors
$v\in \{ -1, 1 \} ^n$ is even $2^n$. Therefore, for $u= p_0t + \zeta
_0 + M(p,t; \zeta )$ due to Theorem 12 the equality
$$(2)\quad \int_{{\bf R}^n} (\partial f(t)/\partial s_j)
e^{-u(p,t;\zeta )}ds = \int_{{\bf R}^n} (\partial f(t)/\partial s_j)
e^{-u(p,t;\zeta )}dt =$$   $$\int_{{\bf R}^{n-1}}(dt^j)
[f(t)e^{-u(p,t;\zeta )}]|_{-\infty }^{\infty } - \int_{{\bf
R}^{n-1}}(dt^j) (\int_{-\infty }^{\infty } f(t) [\partial
e^{-u(p,t;\zeta )}/\partial s_j] ds_j) $$ is satisfied in the ${\cal
A}_r$ spherical coordinates, since the absolute value of the
Jacobian $\partial t/\partial (t^j,s_j)$ is unit. Since for $a_1<Re
(p)<a_{-1}$ the first additive is zero, while the second integral
converts with the help of Formulas 12$(2, 2.1)$, Formula $(1)$
follows for $k=1$:
\par $(3) \quad {\cal F}^n(\partial f(t)/\partial s_j, u; p; \zeta )
= p_0 \delta _{1,j}{\cal F}^n(f(t), u; p; \zeta ) + p_j {\sf
S}_{e_j}{\cal F}^n(f(t), u; p; \zeta )$.
\par To accomplish the derivation we use Theorem 14 so that
$$\lim_{\tau \to 0}[{\cal F}^n(f(t),u;p;\zeta ) - {\cal
F}^n(f(t-\tau e_j),u;p;\zeta )]/\tau $$  $$= \lim_{\tau \to 0}[{\cal
F}^n(f(t),u;p;\zeta ) - {\cal F}^n(f(t),u;p;\zeta + \tau (p_0
+p_1i_1+...+ p_j i_j))]/\tau $$  $$=\lim_{\tau \to 0} \int_{{\bf
R}^n} f(t) [e^{-u(p,t;\zeta )} - e^{-u(p,t;\zeta + \tau
(p_0+p_1i_1+...+p_ji_j))}]\tau ^{-1}dt ,$$ where $e_j
=(0,...,0,1,0,..,0)\in {\bf R}^n$ with $1$ on the $j$-th place. If
the original $\partial ^{|j|} f(t)/\partial s_1^{j_1}...\partial
s_n^{j_n}$ exists, then $\partial ^{|m|} f(t)/\partial
s_1^{m_1}...\partial s_n^{m_n}$ is continuous for $0\le |m|\le
|j|-1$ with $0\le m_l\le j_l$ for each $l=1,...,n$, where $f^0:=f$.
The interchanging of $\lim_{\tau \to 0} $ and $\int_{{\bf R}^n}$ may
change a domain of convergence, but in the indicated in the theorem
domain $a_1< Re (p)<a_{-1}$, when it is non void, Formula $(3)$ is
valid. Applying Formula $(3)$ in the ${\cal A}_r$ spherical
coordinates by induction to $(\partial ^{|m|} f(t)/\partial
s_1^{m_1}...\partial s_n^{m_n}): |m|\le |j|, 0\le m_l\le j_l\forall
l)$ with the corresponding order subordinated to $\partial ^{|j|}
f(t)/\partial s_1^{j_1}...\partial s_n^{j_n}$, or in the ${\cal
A}_r$ Cartesian coordinates using Formula 12$(1.1 )$ for the partial
derivatives $(\partial ^{|m|} f(t)/\partial t_1^{m_1}...\partial
t_n^{m_n}): |m|\le |j|, 0\le m_l\le j_l\forall l)$ with the
corresponding order subordinated to $\partial ^{|j|} f(t)/\partial
t_1^{j_1}...\partial t_n^{j_n}$ we deduce Expressions $(1)$ and
$(1.1)$ with the help of Statement 6 from \S XVII.2.3 \cite{zorich}
about the differentiation of an improper integral by a parameter and
\S 2.

\par {\bf 22. Remarks.} For the entire Euclidean space ${\bf R}^n$
Theorem 21 for $\partial f(t)/\partial s_j$ gives only one or two
additives on the right side of 21$(1)$ in accordance with 21$(3)$.
\par Evidently Theorems 4, 11 and Proposition 10 are accomplished for
${\cal F}^{k; t_{j(1)},...,t_{j(k)}} (f,u;p;\zeta )$ also.  \par
Theorem 12 is satisfied for ${\cal F}^{k; t_{j(1)},...,t_{j(k)}}$
and any $j\in \{ j(1),..., j(k) \} $, so that $s_l =
s_l(k;t)=t_{j(l)}+...+t_{j(k)}$ for each $1\le l\le k$, $p_m=0$ and
$\zeta _m=0$ for each $1\le m \notin \{ j(1),...,j(k) \} $ (the same
convention is in 13, 14, 17, 21, see also below).  For ${\cal F}^{k;
t_{j(1)},...,t_{j(k)}}$ in Theorem 13 in Formula 13$(1)$ it is
natural to put $t_m=0$ and $h_m =0$ for each $1\le m \notin \{
j(1),...,j(k) \} $, so that only $(k+1)$ additives with $h_0$,
$h_{j(1)}$,...,$h_{j(k)}$ on the right side generally may remain.
Theorems 14 and 17 and 21 modify for ${\cal F}^{k;
t_{j(1)},...,t_{j(k)}}$ putting in 14$(1)$ and 17$(1, 2)$ and
21$(1)$ $t_j=0$ and $\tau _j =0$ respectively for each $j\notin \{
j(1),...,j(k) \} $.
\par To take into account boundary conditions for domains different from
$U_v$, for example, for bounded domains $V$ in ${\bf R}^n$ we
consider a bounded noncommutative multiparameter transform
\par $(1)$ ${\cal F}^n(f(t)\chi _V, u; p; \zeta )=: {\cal F}^n_V(f(t), u;
p; \zeta )$. \\  For it evidently Theorems 4, 6-8, 11, 13, 14, 16,
17, Proposition 10 and Corollary 4.1 are satisfied as well taking
specific originals $f$ with supports in $V$.

\par At first take domains $W$ which are quadrants, that is canonical
closed subsets affine diffeomorphic with $Q^n = \prod_{j=1}^n [{\sf
a}_j,b_j]$, where $-\infty \le {\sf a}_j <b_j \le \infty $, $[{\sf
a}_j,b_j] := \{ x\in {\bf R}: ~ {\sf a}_j\le x \le b_j \} $ denotes
the segment in $\bf R$. This means that there exists a vector $w\in
{\bf R}^n$ and a linear invertible mapping $C$ on ${\bf R}^n$ so
that $C(W)-w = Q$.  We put $t^{j,1} := (t_1,...,t_j,...,t_n: ~
t_j={\sf a}_j)$, $t^{j,2} := (t_1,...,t_j,...,t_n: ~ t_j=b_j)$.
Consider $t=(t_1,...,t_n)\in Q^n$.

\par {\bf 23. Theorem.} {\it Let $f(t)$ be a function-original with
a support by $t$ variables in $Q^n$ and zero outside $Q^n$ such that
$\partial f(t)/\partial t_j$ also satisfies Conditions 1$(1-4)$.
Suppose that $u(p,t;\zeta )$ is given by 2$(1,2,2.1)$ or 1$(8,8.1)$
over ${\cal A}_r$ with $2\le r<\infty $, $2^{r-1}\le n \le 2^r-1$.
Then
$$(1)\quad {\cal F}^n((\partial f(t)/\partial t_j) \chi
_{Q^n}(t),u; p;\zeta ) =$$
$${\cal F}^{n-1; t^{j,2}} (f(t^{j,2})\chi _{Q^n}(t^{j,2}), u;p;\zeta ) -
{\cal F}^{n-1; t^{j,1}} (f(t^{j,1})\chi _{Q^n}(t^{j,1}), u;p;\zeta
)$$
$$ +  [p_0 + \sum_{k=1}^j p_k
{\sf S}_{e_k}] {\cal F}^n (f(t)\chi _{Q^n}(t),u; p;\zeta )$$ in the
${\cal A}_r$ spherical coordinates, or
$$(1.1)\quad {\cal F}^n((\partial f(t)/\partial t_j) \chi
_{Q^n}(t),u; p;\zeta ) =$$
$${\cal F}^{n-1; t^{j,2}} (f(t^{j,2})\chi _{Q^n}(t^{j,2}), u;p;\zeta ) -
{\cal F}^{n-1; t^{j,1}} (f(t^{j,1})\chi _{Q^n}(t^{j,1}), u;p;\zeta
)$$
$$ + [p_0  + p_j
{\sf S}_{e_j}] {\cal F}^n (f(t)\chi _{Q^n}(t),u; p;\zeta )$$ in the
${\cal A}_r$ Cartesian coordinates in a domain $W\subset  {\cal
A}_r$; if ${\sf a}_j= - \infty $ or $b_j = + \infty $, then the
addendum with $t^{j,1}$ or $t^{j,2}$ correspondingly is zero.}
\par {\bf Proof.} Here the domain $Q^n$ is bounded and $f$ is almost
everywhere continuous and satisfies Conditions 1$(1-4)$, hence
$f(t)\exp (- u(p,t;\zeta ))\in L^1({\bf R}^n,\lambda _n,{\cal A}_r)$
for each $p\in {\cal A}_r$, since $\exp (- u(p,t;\zeta ))$ is
continuous and $supp (f(t))\subset Q^n$.
\par Analogously to \S 12 the integration by parts gives  $$(2)\quad
\int_{{\sf a}_j}^{b_j}(\partial f(t)/\partial t_j) \exp
(-u(p,t;\zeta ))dt_j = f(t)) \exp (-u(p,t;\zeta ))|_{t_j={\sf
a}_j}^{t_j=b_j}$$  $$ - \int_{{\sf a}_j}^{b_j}[f(t)(\partial \exp
(-u(p,t;\zeta ))/\partial t_j)]dt_j ,$$ where $t=(t_1,...,t_n)$.
Then the Fubini's theorem implies:
$$(3)\quad \int_{Q^n}(\partial f(t)/\partial t_j) \exp (-u(p,t;\zeta
))dt = $$  $$ \int_{{\sf a}_1}^{b_1}...\int_{{\sf
a}_{j-1}}^{b_{j-1}} \int_{{\sf a}_{j+1}}^{b_{j+1}}\int_{{\sf
a}_n}^{b_n} [\int_{{\sf a}_j}^{b_j}(\partial f(t)/\partial t_j) \exp
(-u(p,t;\zeta ))dt_j]dt^j$$  $$ = [\int_{t\in Q^n, ~ t_j=b_j}
f(t^{j,2}) \exp (-u(p,t^{j,2};\zeta )) dt^j] - [\int_{t\in Q^n, ~
t_j={\sf a}_j} f(t^{j,1}) \exp (-u(p,t^{j,1};\zeta )) dt^j]$$
$$ +  [p_0 + \sum_{k=1}^j p_k {\sf S}_{e_k}] \int_{{\sf
a}_1}^{b_1}...\int_{{\sf a}_n}^{b_n} f(t) \exp (-u(p,t;\zeta ))dt
$$ in the ${\cal A}_r$ spherical coordinates or
$$(3.1)\quad \int_{Q^n}(\partial f(t)/\partial t_j) \exp (-u(p,t;\zeta
))dt $$    $$ = [\int_{t\in Q^n, ~ t_j=b_j} f(t^{j,2}) \exp
(-u(p,t^{j,2};\zeta )) dt^j] - [\int_{t\in Q^n, ~ t_j={\sf a}_j}
f(t^{j,1}) \exp (-u(p,t^{j,1};\zeta )) dt^j]$$
$$ +  [p_0  + p_j {\sf S}_{e_j}] \int_{{\sf
a}_1}^{b_1}...\int_{{\sf a}_n}^{b_n} f(t) \exp (-u(p,t;\zeta ))dt
$$ in the ${\cal A}_r$ Cartesian coordinates, where as usually
$t^j=(t_1,...,t_{j-1},0,t_{j+1},...,t_n)$,
$dt^j=dt_1...dt_{j-1}dt_{j+1}...dt_n$. This gives Formulas $(1,
1.1)$, where
$$(4)\quad {\cal F}^{n-1; t^{j,k}}(f(t^{j,k})\chi _{Q^n}(t^{j,k}),
u(p,t^{j,k};\zeta );p;\zeta ) =$$
$$\int_{{\sf a}_1}^{b_1}...\int_{{\sf a}_{j-1}}^{b_{j-1}}
\int_{{\sf a}_{j+1}}^{b_{j+1}} \int_{{\sf a}_n}^{b_n} f(t^{j,k})
\exp (-u(p,t^{j,k};\zeta )) dt^{j,k}$$  is the non-commutative
transform by $t^{j,k}$, $2^{r-1}\le n \le 2^r-1$, $dt^{j,k}$ is the
Lebesgue volume element on ${\bf R}^{n-1}$.

\par {\bf 24. Theorem.} {\it  If a function
$f(t)\chi _{Q^n}(t)$ is original together with its derivative
$\partial ^nf(t)\chi _{Q^n}(t)/\partial s_1...\partial s_n$ or
$\partial ^nf(t)\chi _{Q^n}(t)/\partial t_1...\partial t_n$, where
$F^n_u(p;\zeta )$ is an image function of $f(t)\chi _{Q^n}(t)$ over
the Cayley-Dickson algebra ${\cal A}_r$ with $2\le r\in \bf N$,
$2^{r-1}\le n \le 2^r-1$, for the function $u(p,t;\zeta )$ given by
2$(1,2,2.1)$ or 1$(8,8.1)$, $Q^n = \prod_{j=1}^n [0,b_j]$, $b_j>0$
for each $j$, then
$$(1)\quad \lim_{p\to \infty }
\{ [p_0 + p_1{\sf S}_{e_1}] p_2{\sf S}_{e_2}...p_n{\sf S}_{e_n}
F^n_u(p;\zeta ) + \sum_{m=0}^{n-1} (-1)^m $$
$$ \sum_{1 \le j_1 <...< j_{n-m} \le n; ~ 1\le l_1<...<l_m\le n; ~
l_{\alpha }\ne j_{\beta } ~ \forall \alpha , \beta } [p_0\delta
_{1,j_1} + p_{j_1}{\sf S}_{e_{j_1}}] p_{j_2}{\sf
S}_{e_{j_2}}...p_{j_{n-m}}{\sf S}_{e_{j_{n-m}}}$$
$$F^{n-m}_u(p^{(l)}; \zeta ) \} = (-1)^{n+1} f(0) e^{ -
u(0,0;\zeta )}$$ in the ${\cal A}_r$ spherical coordinates, or
$$(1.1)\quad \lim_{p\to \infty }
\{ [p_0 + p_1{\sf S}_{e_1}] [p_0+p_2{\sf S}_{e_2}]...[p_0+p_n{\sf
S}_{e_n}] F^n_u(p;\zeta ) + \sum_{m=0}^{n-1} (-1)^m $$
$$ \sum_{1 \le j_1 <...< j_{n-m} \le n; ~ 1\le l_1<...<l_m\le n; ~
l_{\alpha }\ne j_{\beta } ~ \forall \alpha , \beta } [p_0 +
p_{j_1}{\sf S}_{e_{j_1}}] [p_0+p_{j_2}{\sf
S}_{e_{j_2}}]...[p_0+p_{j_{n-m}}{\sf S}_{e_{j_{n-m}}}]$$
$$F^{n-m}_u(p^{(l)}; \zeta ) \} = (-1)^{n+1} f(0) e^{ -
u(0,0;\zeta )}$$ in the ${\cal A}_r$ Cartesian coordinates,  where
$f(0) = \lim_{t\in Q^n, ~ t\to 0 } f(t)$, $p$ tends to the infinity
inside the angle $|Arg (p)|<\pi /2-\delta $ for some $0<\delta <\pi
/2$. }
\par {\bf Proof.} In accordance with Theorem 23 we have Equalities 23$(1,1.1)$.
Therefore we infer that
$$(2)\quad {\cal F}^{n-1;t^{j,k}}
(f(t^{j,k})\chi _{Q^n}(t^{j,k}), u(p,t^{j,k};\zeta );p;\zeta ) =$$
$$\lim_{t_j\to \beta _{j,k}+0} \int_{{\sf a}_1}^{b_1}dt_1...
\int_{{\sf a}_{j-1}}^{b_{j-1}}dt_{j-1} \int_{{\sf
a}_{j+1}}^{b_{j+1}} dt_{j+1}...\int_{{\sf a}_n}^{b_n} (dt_n) f(t)
\exp (-u(p,t;\zeta )),$$ where $\beta _{j,1} ={\sf a}_j=0$, $\beta
_{j,2}=b_j>0$, $k=1, 2$. Mention, that
$(...((t^{1,l_1})^{2,l_2})...)^{j,l_j} = (t: ~ t_1=\beta
_{1,l_1},...,t_j=\beta _{j,l_j})$ for every $1\le j\le n$.
Analogously to \S 12 we apply Formula $(2)$ by induction
$j=1,...,n$, $2^{r-1}\le n \le 2^r-1$, to \par $\partial
^nf(t(s))/\partial s_1...\partial s_n$,...,$\partial ^{n-j+1}
f(t(s))/\partial s_j...\partial s_n$,...,$\partial f(t(s))/\partial
s_n$ \\ instead of $\partial f(t(s))/\partial s_j$, $s_j=s_j(n;t)$
as in \S 2, or applying to the partial derivatives \par $\partial
^nf(t)/\partial t_1...\partial t_n$,...,$\partial ^{n-j+1}
f(t)/\partial t_j...\partial t_n$,...,$\partial f(t)/\partial t_n$
\\  instead of $\partial f(t)/\partial t_j$ correspondingly. If $s_j>0$
for some $j\ge 1$, then $s_1>0$ for $Q^n$ and $\lim_{p\to \infty }
e^{-u(p,t^{(l)};\zeta )}=0$ for such $t^{(l)}$, where $t =
(t_1,...,t_n)$, $(l) = (l_1,...,l_n)$, $|l| = l_1+...+l_n$, $t^{(l)}
= (t^{(l)}_1,...,t^{(l)}_n)$, $t^{(l)}_j={\sf a}_j$ for $l_j=1$ and
$t^{(l)}_j=b_j$ for $l_j=2$, $1\le j\le 2^r-1$. Therefore,
$$\lim_{p\to \infty } \sum_{l_j\in \{ 1, 2 \}; ~ j=1,...,n }
(-1)^{|l|} f(t^{(l)}) e^{-u(p,t^{(l)}; \zeta ) } = (-1)^n f(0) e^{
-u(0,0;\zeta )},$$  since $u(p,0;\zeta ) = u(0,0;\zeta )$, where
$f(^{(l)}) = \lim_{t\in Q^n; t\to t^{(l)}} f(t)$.
\par In accordance with Note 8 \cite{lutsltjms}  $$\lim_{p\to
\infty , |Arg (p)|<\pi /2-\delta } {\cal F}^n((\partial
^nf(t)/\partial s_1...\partial s_n)\chi _{Q^n}(t),u(p,t;\zeta
);p;\zeta )=0$$ in the ${\cal A}_r$ spherical coordinates and
$$\lim_{p\to \infty , |Arg (p)|<\pi /2-\delta } {\cal F}^n((\partial
^nf(t)/\partial t_1...\partial t_n)\chi _{Q^n}(t),u(p,t;\zeta
);p;\zeta )=0$$ in the ${\cal A}_r$ Cartesian coordinates, which
gives the statement of this theorem.

\par {\bf 25.} Suppose that
$f(t)\chi _{Q^n}(t)$ is an original function, $F^n(p;\zeta )$ is its
image, $\partial ^{|j|} f(t)\chi _{Q^n}(t)/\partial
t_1^{j_1}...\partial t_n^{j_n}$ is an original, $|j|=j_1+...+j_n$,
$0\le j_1,...,j_n\in {\bf Z}$, $2^{r-1}\le n \le 2^r-1$, $ - \infty
\le {\sf a}_k<b_k\le \infty $ for each $k=1,...,n$,
$(l)=(l_1,...,l_n)$, $l_k\in \{ 0, 1, 2 \}$, $W={\cal A}_r$ for
bounded $Q^n$. Let $W= \{ p\in {\cal A}_r: ~ a_1< Re (p) \} $ for
$b_k=\infty $ for some $k$ and finite ${\sf a}_k$ for each $k$;
$W=\{ p \in {\cal A}_r: ~ Re (p)<a_{-1} \} $ for ${\sf a}_k= -\infty
$ for some $k$ and finite $b_k$ for each $k$; $W= \{ p \in {\cal
A}_r: ~ a_1< Re (p) < a_{-1} \} $ when ${\sf a}_k = -\infty $ and
$b_l=+\infty $ for some $k$ and $l$;
$t^{(l)}=(t^{(l)}_1,...,t^{(l)}_n)$. \par We put $t^{(l)}_k=t_k$ and
$q_k=0$ for $l_k=0$, $t^{(l)}_k={\sf a}_k$ for $l_k=1$,
$t^{(l)}_k=b_k$ for $l_k=2$, $ ~ ~ (q)=(q_1,...,q_n)$, $ ~ ~
|q|=q_1+...+q_n$,
\par $a_1= \max (a_1(f),a_1(\partial ^{|m|} f(t)/\partial
t_1^{m_1}...\partial t_n^{m_n}): ~ |m|\le |j|, 0\le m_k\le
j_k\forall k)$, \par $a_{-1} = \min (a_{-1}(f),a_{-1}(\partial
^{|m|} f(t)/\partial t_1^{m_1}...\partial t_n^{m_n}): ~ |m|\le |j|,
0\le m_k\le j_k ~ \forall k)$ if $a_1<a_{-1}$. \par If ${\sf a}_k= -
\infty $ and $b_k = + \infty $ for $Q^n$ with a given $k$, then
$l_k=0$. If either ${\sf a}_k>-\infty $ or $b_k<+\infty $ for a
marked $k$, then $l_k \in \{ 0, 1, 2 \} $. We also put
$h_k=h_k(l)=sign (l_k)$ for each $k$, where $sign (x)=-1$ for $x<0$,
$sign (0)=0$, $sign (x)=1$ for $x>0$, $h=h(l)$,
$|h|=|h_1|+...+|h_n|$, \par $(lj) := (l_1 sign (j_1),...,l_n sign
(j_n))$. \par Let the vector $(l)$ enumerate faces $\partial
Q^n_{(l)}$ in $\partial Q^n_{k-1}$ for $|h(l)|= k\ge 1$, so that
$\partial Q^n_{k-1} = \bigcup_{|h(l)|=k} Q^n_{(l)}$, $\partial
Q^n_{(l)}\cap \partial Q^n_{(m)}=\emptyset $ for each $(l)\ne (m)$
(see also more detailed notations in \S 28).
\par Let the shift operator be defined:
\par $T_{(m)}F(p;\zeta ) := F(p;\zeta - (i_1 m_1 +...+i_nm_n)\pi /2)$, also
the operator
\par $(SO)$ ${\sf S}_{(m)}F(p;\zeta ) := {\sf
S}_{e_1}^{m_1}...{\sf S}_{e_n}^{m_n}F(p;\zeta )$, \\
where $(m)=(m_1,...,m_n)\in [0,\infty )^n\subset {\bf R}^n$, ${\sf
S}_{(m)}^k={\sf S}_{k(m)}$ for each positive number $0<k\in {\bf
R}$, ${\sf S}_0=I$ is the unit operator for $(m)=0$ (see also
Formulas 12$(3.1-3.7)$). As usually let
$e_1=(1,0,...,0)$,...,$e_n=(0,...,0,1)$ be the standard orthonormal
basis in ${\bf R}^n$ so that $(m)=m_1e_1+...+m_ne_n$.
\par {\bf Theorem.} {\it Then
$$(1)\quad {\cal F}^n(\partial ^{|j|}
f(t)\chi _{Q^n}(t)/\partial t_1^{j_1}...\partial t_n^{j_n},
u(p,t;\zeta ); p;\zeta ) = $$
$$ {\sf R}_{e_1}^{j_1} {\sf R}_{e_2}^{j_2}
...{\sf R}_{e_n}^{j_n} {\cal F}^n (f(t)\chi _{Q^n}(t),u; p;\zeta )$$
$$+ \sum_{1\le |(lj)|; ~ m_k+q_k+h_k=j_k; 0\le m_k; ~ 0\le q_k; ~ h_k= sign (l_kj_k);
~ q_k=0 \mbox{ for } l_kj_k=0 \mbox{, for each } k=1,...,n; ~ (l)\in
\{ 0, 1, 2 \} ^n} $$
$$ (-1)^{|(lj)|} {\sf R}_{e_1}^{m_1} {\sf R}_{e_2}^{m_2}
...{\sf R}_{e_n}^{m_n} {\cal F}^{n -|h(lj)|} (\partial ^{|q|}
f(t^{(lj)})\chi _{\partial Q^n_{(lj)}}(t^{(lj)})/\partial
t_1^{q_1}...\partial t_n^{q_n},u; p;\zeta )$$ for $u(p,t;\zeta )$ in
the ${\cal A}_r$ spherical coordinates or the ${\cal A}_r$ Cartesian
coordinates over the Cayley-Dickson algebra ${\cal A}_r$ with $2\le
r<\infty $, where \par $(1.1)$ ${\sf R}_{e_1} := p_0 + p_1{\sf
S}_{e_1}$, ${\sf R}_{e_2} := p_0 + p_1{\sf S}_{e_1} + p_2 {\sf
S}_{e_2}$,..., ${\sf R}_{e_n} := p_0 + p_1{\sf S}_{e_1} + p_2 {\sf
S}_{e_2}+ ...+p_n {\sf S}_{e_n}$ in the ${\cal A}_r$ spherical
coordinates, while
\par $(1.2)$ ${\sf R}_{e_1} := p_0 + p_1{\sf S}_{e_1}$, ${\sf R}_{e_2} := p_0
+ p_2 {\sf S}_{e_2}$,..., ${\sf R}_{e_n} := p_0 +p_n {\sf S}_{e_n}$
in the ${\cal A}_r$ Cartesian coordinates, i.e. ${\sf R}_{e_j}= {\sf
R}_{e_j}(p)$ are operators depending on the parameter $p$. If
$t^{(l)}_j=\infty $ for some $1\le j \le n$, then the corresponding
addendum on the right of $(1)$ is zero.}
\par {\bf Proof.} In view of Theorem 23 we get the equality
$$(2)\quad \int_{Q^n} [(\partial ^{|m|+1}f(t)/\partial
t_1^{m_1}...\partial t_{k-1}^{m_{k-1}}\partial t_k^{m_k+1}\partial
t_{k+1}^{m_{k+1}}... \partial t_n^{m_n}) e^{-u(p,t;\zeta )} ] dt =$$
$$\int_{{\bf R}^{n-1}\cap Q^n}(dt^k)
[(\partial ^{|m|} f(t)/\partial t_1^{m_1}...
\partial t_n^{m_n}) e^{-u(p,t;\zeta )}]|_{{\sf a}_k}^{b_k}$$
$$ - \int_{{\bf
R}^{n-1}\cap Q^n}(dt^k) (\int_{{\sf a}_k }^{b_k} (\partial ^{|m|}
f(t)/\partial t_1^{m_1}...
\partial t_n^{m_n}) [\partial
e^{-u(p,t;\zeta )}/\partial t_k] dt_k)$$ is satisfied for $0\le
m_k\le j_k$ for each $k=1,...,n$ with $|m|<|j|$. On the other hand,
for $p\in W$ additives on the right of $(2)$ convert with the help
of Formula 23$(1)$.  Each term of the form
$$\int_{{\bf R}^{n-|h(l)|}\cap Q^n}(dt^{(l)}) [(\partial
^{|q|} f(t^{(l)})\chi _{\partial Q^n_{(l)}}(t^{(l)})/\partial
t_1^{q_1}...\partial t_n^{q_n}e^{-u(p,t;\zeta )}]$$ can be further
transformed with the help of $(2)$ by the considered variable $t_k$
only in the case $l_k=0$. Applying Formula $(2)$ by induction to
partial derivatives $\partial ^{|j|}f/\partial t_1^{j_1}...\partial
t_n^{j_n}$, $\partial ^{|j|-j_1}f/\partial t_2^{j_2}...\partial
t_n^{j_n}$,...,$\partial ^{j_n}f/
\partial t_n^{j_n}$,...,$\partial f/\partial t_n$
as in \S 21 and using Theorem 14 and Remarks 22 we deduce $(1)$.

\par {\bf 26. Theorem.} {\it Let $f(t)\chi _{U_{1,...,1}}(t)$ be a
function-original with values in ${\cal A}_r$ with $2\le r<\infty $,
$2^{r-1}\le n \le 2^r-1$, $u$ is given by 2$(1,2,2.1)$ or
1$(8,8.1)$,
$$(1)\quad g(t) := \int_0^{t_1}...\int_0^{t_n} f(x)dx, \mbox{
then}$$
$$(2)\quad {\cal F}^n(f\chi _{U_{1,...,1}}(t),
u; p;\zeta ) = {\sf R}_{e_1} {\sf R}_{e_2}... {\sf R}_{e_n} {\cal
F}^n (g(t)\chi _{U_{1,...,1}}(t),u; p;\zeta )$$ in the domain $Re
(p)>\max (a_1,0)$, where the operators ${\sf R}_{e_j}$ are given by
Formulas 25$(1.1,1.2)$.}
\par {\bf Proof.} In view of Theorem 25 the equation
$$(3)\quad {\cal F}^n(f\chi _{U_{1,...,1}}(t),
u; p;\zeta ) = $$
$$ {\sf R}_{e_1} {\sf R}_{e_2}...
{\sf R}_{e_n} {\cal F}^n (g(t),u; p;\zeta )$$
$$ + \sum_{1\le |l|; ~ 0\le m_k\le 1; ~ m_k+h_k=1; ~ h_k=sign (l_k); ~
\mbox{ for each } k=1,...,n; ~ q_1=0,..., q_n=0}
$$
$$ (-1)^{|(l)|} {\sf R}_{e_1}^{m_1} {\sf R}_{e_2}^{m_2}...
{\sf R}_{e_n}^{m_n} {\cal F}^{n -|h(l)|} ( g(t^{(l)}),u; p;\zeta
),$$ is satisfied, since $\partial ^n g(t)/\partial t_1...\partial
t_n= (f\chi_{U_{1,...,1}})(t)$, where $j_1=1$,...,$j_n=1$, $l_j=1$
for each $j=1,...,n$. Equation $(3)$ is accomplished in the same
domain $Re (p)>\max (a_1,0)$, since $g(0)=0$ and $g(t)$ also
fulfills conditions of Definition 1, while $a_1(g)<\max
(a_1(f),0)+b$ for each $b>0$, where $a_1\in \bf R$. On the other
hand, $g(t)$ is equal to zero on $\partial U_{1,...,1}$ and outside
$U_{1,...,1}$ in accordance with formula $(1)$, hence all terms on
the right side of Equation $(3)$ with $|l|>0$ disappear and $supp
(g(t)) \subset U_{1,...,1}$. Thus we get Equation $(2)$.

\par {\bf 27. Theorem.} {\it Suppose that $F^k(p;\zeta )$ is an image
${\cal F}^{k;t_1,...,t_k}(f(t)\chi _{U_{1,...,1}}(t),u;p;\zeta )$ of
an original function $f(t)$ for $u$ given by 2$(1,2,2.1)$ in the
half space $W := \{ p\in {\cal A}_r: Re (p)>a_1 \} $ with $2\le
r<\infty $, $p_1=0$,...,$p_{j-1}=0$; $\zeta _1=\pi /2$,...,$\zeta
_{j-1}=\pi /2$ for each $j\ge 2$ in the ${\cal A}_r$ spherical
coordinates or $\zeta _1=0$,...,$\zeta _{j-1}=0$ for each $j\ge 2$
in the ${\cal A}_r$ Cartesian coordinates;
\par $(1)$ the integral $\int_{p_j i_j}^{\infty i_j} F^k(p_0 +
z;\zeta ) dz$ converges, where $p=p_0+p_1i_1+...+p_ki_k\in {\cal
A}_r$, $p_j\in {\bf R}$ for each $j=0,...,2^r-1$, $2^{r-1}\le k\le
2^r-1$, $U_{1,...,1} := \{ (t_1,...,t_k)\in {\bf R}^k: ~ t_1\ge 0,
...,t_k\ge 0 \} $. Let also
\par $(2)$ the function $F^k(p;\zeta )$ be continuous by the variable
$p\in {\cal A}_r$ on the open domain $W$, moreover, for each $w>a_1$
there exist constants ${C_w}'
>0$ and $\epsilon _w >0$ such that
\par $(3)$ $|F^k(p;\zeta ) |\le {C_w}'
\exp (-\epsilon _w |p|)$ for each $p\in S_{R(n)}$, $S_R := \{ z\in
{\cal A}_r: ~ Re (z)\ge w \} $, $0<R(n)<R(n+1)$ for each $n\in \bf
N$, $\lim_{n\to \infty }R(n)=\infty $, where $a_1$ is fixed, $\zeta
=\zeta _0i_0+...+\zeta _ki_k\in {\cal A}_r$ is marked, $\zeta _j\in
{\bf R}$ for each $j=0,...,k$. Then $$(4)\quad \int_{p_j
i_j}^{\infty i_j} F^k(p_0 + z;\zeta ) dz = {\cal S}_{ - e_j} {\cal
F}^{k;t_1,...,t_k}(f(t)\chi _{U_{1,...,1}}(t)/\xi _j, u;p;\zeta ),$$
where $p_1=0$,...,$p_{j-1}=0$ for each $j\ge 2$; $\zeta _1=\pi
/2$,...,$\zeta _{j-1}=\pi /2$ and $\xi _j=s_j(k;t)$ in the ${\cal
A}_r$ spherical coordinates, while $\zeta _1= 0$,...,$\zeta
_{j-1}=0$ and $\xi _j=t_j$ in the ${\cal A}_r$ Cartesian coordinates
correspondingly for each $j\ge 1$.}
\par {\bf Proof.} Take a path of an integration belonging
to the half space $Re (p)\ge w$ for some constant $w>a_1$. Then
$$|\int_{U_{1,...,1}}f(t)\exp (- u(p,t;\zeta ))dt|\le
C\int_{U_{1,...,1}} \exp (-(p_0-a_1)(t_1+...+t_k))dt < \infty $$
converges, where $C=const>0$, $p_0\ge w$. For $t_j>0$ for each
$j=1,...,k$ conditions of Lemma 2.23 \cite{lutsltjms} (that is of
the noncommutative analog over ${\cal A}_r$ of Jordan's lemma) are
satisfied. If $t_j\to \infty $, then $s_j\to \infty $, since all
$t_1$,...,$t_k$ are non-negative. Up to a set $\partial U_{1,...,1}$
of $\lambda _k$ Lebesgue measure zero we can consider that
$t_1>0$,...,$t_k>0$. If $s_j\to \infty $, then also $s_1\to \infty
$. The converging integral can be written as the following limit:
$$(5)\quad \int_{p_j i_j}^{\infty i_j}
F^k(p_0 + z;\zeta ) dz = \lim_{0<\kappa \to 0} \int_{p_j
i_j}^{\infty i_j} F^k(p_0 + z;\zeta ) \exp (-\kappa |z| ) dz$$ for
$1\le j\le k$, since the integral $\int_{-S\infty }^{S\infty }[F^k(w
+ z;\zeta ) ] dz$ is absolutely converging and the limit
$\lim_{\kappa \to 0}\exp (-\kappa |z|)=1$ uniformly by $z$ on each
compact subset in ${\cal A}_r$, where $S$ is a purely imaginary
marked Cayley-Dickson number with $|S|=1$. Therefore, in the
integral  $$(6)\quad \int_{p_j i_j}^{\infty i_j} F^k(p_0 + z;\zeta )
 dz = \int_{p_j i_j}^{\infty i_j} (\int_{U_{1,...,1}} f(t)[\exp
(- u(p_0+z,t;\zeta )) ]dt) dz$$ the order of the integration can be
changed in accordance with the Fubini's theorem applied
componentwise to an integrand $g=g_0i_0+...+g_ni_n$ with $g_l\in
{\bf R}$ for each $l=0,...,n$:
$$(7)\quad \int_{p_j i_j}^{\infty i_j}
F^k(p_0 + z;\zeta ) dz = \int_{U_{1,...,1}} dt (\int_{p_j
i_j}^{\infty i_j} (f(t) \exp (- u(p_0 + z,t;\zeta )) dz)$$  $$ =
\int_{U_{1,...,1}} f(t) \{ \int_{p_j i_j}^{\infty i_j} [e^{-
u(p_0+z,t;\zeta )} ] dz \} dt .$$  Generally, the condition
$p_1=0$,...,$p_{j-1}=0$ and $\zeta _1=\pi /2$,...,$\zeta _{j-1}=\pi
/2$ in the ${\cal A}_r$ spherical coordinates or $\zeta
_1=0$,...,$\zeta _{j-1}=0$ in the ${\cal A}_r$ Cartesian coordinates
for each $j\ge 2$ is essential for the convergence of such integral.
We certainly have
$$(8)\quad \int_{p_j i_j}^{b_j i_j} \cos (i_j^*z\xi _j+\zeta _j) dz =
[\sin (\theta _j \xi _j + \zeta _j)/\xi _j ]|_{\theta
_j=p_j}^{\theta _j=b_j} = [ - \cos (\theta _j\xi _j + \zeta _j + \pi
/2)/\xi _j ]|_{\theta _j=p_j}^{\theta _j=b_j}$$ and
$$(9)\quad \int_{p_j i_j}^{b_j i_j} \sin (i_j^*z \xi _j+\zeta _j) dz_j =
[- \cos (\theta _j \xi _j + \zeta _j)/\xi _j]|_{\theta
_j=p_j}^{\theta _j=b_j} = [ - \sin (\theta _j\xi _j + \zeta _j + \pi
/2)/\xi _j]|_{\theta _j=p_j}^{\theta _j=b_j}$$ for each $\xi _j>0$
and $-\infty < p_j<b_j<\infty $ and $j=1,...,k$.
 Applying Formulas $(5-9)$ and 2$(1,2,2.1)$ or 1$(8,8.1)$ and 12$(3.1-3.7)$ we deduce that:
$$\int_{p_j i_j}^{\infty i_j}
[F^k(p_0 + z;\zeta )] dz = {\sf S}_{-e_j}\int_{U_{1,...,1}}
[f(t)/\xi _j] \exp \{ - u(p,t;\zeta ) \} dt $$   $$  = {\sf
S}_{-e_j} {\cal F}^{k;t_1,...,t_k}(f(t)\chi _{U_{1,...,1}}(t)/\xi
_j, u;p;\zeta ) ,$$ where $t=(t_1,...,t_k)$, $s_j=t_j+...+t_k$ for
each $1\le j<k$, $s_k=t_k$, $\xi _j=s_j$ in the ${\cal A}_r$
spherical coordinates or $\xi _j=t_j$ in the ${\cal A}_r$ Cartesian
coordinates.

\par {\bf 28. Application of the noncommutative multiparameter transform
to partial differential equations.}
\par Consider a partial differential equation of the form:
\par $(1)$  $A[f](t) = g(t)$,  where
\par $(2)$  $A[f](t) := \sum_{|j|\le \alpha } {\bf a}_j(t) (\partial ^{|j|}
f(t)/\partial t_1^{j_1}...\partial t_n^{j_n}),$ \\ ${\bf a}_j(t)\in
{\cal A}_{\kappa }$ are continuous functions, where $0\le \kappa \in
{\bf Z}$, $j=(j_1,...,j_n)$, $|j| := j_1+...+j_n$, $0\le j_k\in {\bf
Z}$, $\alpha $ is a natural order of a differential operator $A$,
$2\le r$, $2^{r-1}\le n \le 2^r-1$. Since $s_k=s_k(n;t)
=t_k+...+t_n$ for each $k=1,...,n$, the operator $A$ can be
rewritten in $s$ coordinates as
\par $(2.1)$ $A[f](t(s)) :=
\sum_{|j|\le \alpha } {\bf b}_j(t) (\partial ^{|j|} f(t(s))/\partial
s_1^{j_1}...\partial s_n^{j_n}).$ \\ That is, there exists ${\bf
b}_j\ne 0$ for some $j$ with $|j|=\alpha $ and ${\bf b}_j=0$ for
$|j|>\alpha $, while a function $\sum_{j, |j|=\alpha } {\bf
b}_j(t(s)) s_1^{j_1}...s_n^{j_n}$ is not zero identically on the
corresponding domain $V$. We consider that \par $(D1)$ $U$ is a
canonical closed subset in the Euclidean space ${\bf R}^n$, that is
$U = cl ( Int (U))$, where $Int (U)$ denotes the interior of $U$ and
$cl (U)$ denotes the closure of $U$. \par Particularly, the entire
space ${\bf R}^n$ may also be taken. Under the linear mapping
$(t_1,...,t_n)\mapsto (s_1,...,s_n)$ the domain $U$ transforms onto
$V$.
\par We consider a manifold $W$ satisfying the following conditions $(i-v)$.
\par $(i)$. The manifold $W$ is continuous and
piecewise $C^{\alpha }$, where $C^l$ denotes the family of $l$ times
continuously differentiable functions. This means by the definition
that $W$ as the manifold is of class $C^0\cap C^{\alpha }_{loc}$.
That is $W$ is of class $C^{\alpha }$ on open subsets $W_{0,j}$ in
$W$ and $W\setminus (\bigcup_j W_{0,j})$ has a codimension not less
than one in $W$.
\par $(ii)$. $W=\bigcup_{j=0}^m W_{j}$, where $W_{0} = \bigcup_k
W_{0,k}$, $W_{j}\cap W_{k} = \emptyset $ for each $k\ne j$, $m =
dim_{\bf R} W$, $dim_{\bf R} W_{j} = m-j$, $W_{j+1}\subset
\partial W_{j}$. \par $(iii)$. Each $W_{j}$ with $j=0,...,m-1$ is an
oriented $C^{\alpha }$-manifold, $W_{j}$ is open in $\bigcup_{k=j}^m
W_{k}$. An orientation of $W_{j+1}$ is consistent with that of
$\partial W_{j}$ for each $j=0,1,...,m-2$. For $j>0$ the set $W_{j}$
is allowed to be void or non-void.
\par $(iv)$. A sequence $W^k$ of $C^{\alpha }$ orientable manifolds
embedded into ${\bf R}^n$, $\alpha \ge 1$, exists such that $W^k$
uniformly converges to $W$ on each compact subset in ${\bf R}^n$
relative to the metric $dist$. \par For two subsets $B$ and $E$ in a
metric space $X$ with a metric $\rho $ we put \par $(3)\quad dist
(B,E) := \max \{ \sup_{b\in B} dist (\{ b \} ,E), \sup_{e\in E} dist
(B,\{ e \} ) \} ,$ where \par $dist (\{ b \} ,E) := \inf_{e\in E}
\rho (b,e)$, $dist (B, \{ e \} ) := \inf_{b\in B} \rho (b,e)$, $b\in
B$, $e\in E$.
\par Generally, $dim_{\bf R} W=m\le n$. Let
$(e_1^k(x),...,e_m^k(x))$ be a basis in the tangent space $T_xW^k$
at $x\in W^k$ consistent with the orientation of $W^k$, $k\in {\bf
N}$.
\par We suppose that the sequence of orientation frames $(e^k_1(x_k),...,e_m^k(x_k))$
of $W^k$ at $x_k$ converges to $(e_1(x),...,e_m(x))$ for each $x\in
W_0$, where $\lim_kx_k = x\in W_0$, while $e_1(x)$,...,$e_m(x)$ are
linearly independent vectors in ${\bf R}^n$.
\par $(v)$. Let a sequence of Riemann volume elements $\lambda _k$
on $W^k$ (see \S XIII.2 \cite{zorich}) induce a limit volume element
$\lambda $ on $W$, that is, $\lambda (B\cap W) = \lim_{k\to \infty }
(B\cap W^k)$ for each compact canonical closed subset $B$ in ${\bf
R}^n$, consequently, $\lambda (W\setminus W_0)=0$. We shall consider
surface integrals of the second kind, i.e. by the oriented surface
$W$ (see $(iv)$), where each $W_j$, $j=0,...,m-1$ is oriented (see
also \S XIII.2.5 \cite{zorich}).
\par Recall, that a subset $V$ in ${\bf
R}^n$ is called convex, if from $a, b\in V$ it follows that
$\epsilon a + (1-\epsilon )b\in V$ for each $\epsilon \in [0,1]$.
\par $(vi)$. Let a vector $w\in Int (U)$ exist so that $U-w$ is convex in
${\bf R}^n$ and let $\partial U$ be connected. Suppose that a
boundary $\partial U$ of $U$ satisfies Conditions $(i-v)$ and \par
$(vii)$ let the orientations of $\partial U^k$ and $U^k$ be
consistent for each $k\in {\bf N}$ (see Proposition 2 and Definition
3 \cite{zorich}).
\par Particularly, the Riemann volume element $\lambda _k$ on
$\partial U^k$ is consistent with the Lebesgue measure on $U^k$
induced from ${\bf R}^n$ for each $k$. This induces the measure
$\lambda $ on $\partial U$ as in $(v)$.
\par Also the boundary conditions are imposed:
\par $(4)$ $f(t)|_{\partial U} = f_0(t'),$
$(\partial ^{|q|}f(t)/\partial s_1^{q_1}...\partial s_n^{q_n}
)|_{\partial U} = f_{(q)}(t')$ for $|q|\le \alpha -1$, where
$s=(s_1,...,s_n) \in {\bf R}^n$, $(q)=(q_1,...,q_n)$, $|q|
=q_1+...+q_n$, $0\le q_k\in {\bf Z}$ for each $k$, $t\in \partial U$
is denoted by $t'$, $f_0$, $f_{(q)}$ are given functions. Generally
these conditions may be excessive, so one uses some of them or their
linear combinations (see $(5.1)$ below). Frequently, the boundary
conditions
\par $(5)$ $f(t)|_{\partial U} = f_0(t'),$
$(\partial ^lf(t)/\partial \nu ^l)|_{\partial U} = f_l(t')$ for
$1\le l\le \alpha -1$ are also used, where $\nu $ denotes a real
variable along a unit external normal to the boundary $\partial U$
at a point $t'\in \partial U_0$. Using partial differentiation in
local coordinates on $\partial U$ and $(5)$ one can calculate in
principle all other boundary conditions in $(4)$ almost everywhere
on $\partial U$.
\par Suppose that a domain $U_1$ and its boundary $\partial U_1$ satisfy
Conditions $(D1,i-vii)$ and $g_1=g\chi _{U_1}$ is an original on
${\bf R}^n$ with its support in $U_1$. Then any original $g$ on
${\bf R}^n$ gives the original $g\chi _{U_2}=: g_2$ on ${\bf R}^n$,
where $U_2={\bf R}^n\setminus U_1$. Therefore, $g_1+g_2$ is the
original on ${\bf R}^n$, when $g_1$ and $g_2$ are two originals with
their supports contained in $U_1$ and $U_2$ correspondingly. Take
now new domain $U$ satisfying Conditions $(D1,i-vii)$ and $(D2-D5)$:
\par $(D2)$ $U\supset U_1$ and $\partial U\subset \partial U_1$;
\par $(D3)$ if a straight line $\xi $ containing a point $w_1$ (see
$(vi)$) intersects $\partial U_1$ at two points $y_1$ and $y_2$,
then only one point either $y_1$ or $y_2$ belongs to $\partial U$,
where $w_1\in U_1$, $U-w_1$ and $U_1-w_1$ are convex; if $\xi $
intersects $\partial U_1$ only at one point, then it intersects
$\partial U$ at the same point. That is,
\par $(D4)$  any straight line $\xi $ through the point $w_1$ either does not intersect
$\partial U$ or intersects the boundary $\partial U$ only at one
point. \par Take now $g$ with $supp (g)\subset U$, then $supp (g
\chi _{U_1})\subset U_1$. Therefore, any problem $(1)$ on $U_1$ can
be considered as the restriction of the problem $(1)$ defined on
$U$, satisfying $(D1-D4,i-vii)$. Any solution $f$ of $(1)$ on $U$
with the boundary conditions on $\partial U$ gives the solution as
the restriction $f|_{U_1}$ on $U_1$ with the boundary conditions on
$\partial U$. \par Henceforward, we suppose that the domain $U$
satisfies Conditions $(D1,D4,i-vii)$, which are rather mild and
natural. In particular, for $Q^n$ this means that either $a_k = -
\infty $ or $b_k = + \infty $ for each $k$. Another example is:
$U_1$ is a ball in ${\bf R}^n$ with the center at zero, $U=U_1\cup
({\bf R}^n\setminus U_{1,...,1})$, $w_1=0$; or $U=U_1\cup \{ t\in
{\bf R}^n: ~ t_n\ge - \epsilon \} $ with a marked number $0<\epsilon
<1/2$. But subsets $\partial U_{(l)}$ in $\partial U$ can also be
specified, if the boundary conditions demand it.
\par The complex field has the natural realization by $2\times 2$ real matrices
so that ${\bf i} = {{~0 ~1} \choose {-1 ~ 0}}$, ${\bf i}^2= - {{~1
~0} \choose {~ 0 ~ 1}}$. The quaternion skew field, as it is
well-known, can be realized with the help of $2\times 2$ complex
matrices with the generators $I = {{~1 ~0} \choose {~ 0 ~ 1}}$, $J =
{{~0 ~1} \choose {-1 ~ 0}}$, $K = {{ {\bf i} ~ ~ 0} \choose {0 ~ -
{\bf i} }}$, $L = {{0 ~ - {\bf i}} \choose { - {\bf i} ~ 0}}$, or
equivalently by $4\times 4$ real matrices. Considering matrices with
entries in the Cayley-Dickson algebra ${\cal A}_v$ one gets the
complexified or quaternionified Cayley-Dickson algebras $({\cal
A}_v)_{\bf C}$ or $({\cal A}_v)_{\bf H}$ with elements $z=aI+b{\bf
i}$ or $z=aI+bJ+cK+eL$, where $a, b, c, e \in {\cal A}_v$, such that
each $a\in {\cal A}_v$ commutes with the generators ${\bf i}$, $I$,
$J$, $K$ and $L$.
\par When $r=2$, $f$ and $g$ have values in ${\cal A}_2={\bf H}$ and $2\le n\le 4$ and coefficients
of differential operators belong to ${\cal A}_2$, then the
multiparameter noncommutative transform operates with the
associative case so that \par ${\cal F}^n(af)=a{\cal F}^n(f)$
\\ for each $a\in {\bf H}$. The left linearity property ${\cal F}^n(af)=a{\cal F}^n(f)$
for any $a\in {\bf H}_{J,K,L}$ is also accomplished for either
operators with coefficients in ${\bf R}$ or ${\bf C}_{\bf i} = I
{\bf R}\oplus {\bf i}{\bf R}$ or ${\bf H}_{J,K,L} = I {\bf R} \oplus
J {\bf R} \oplus K {\bf R} \oplus L {\bf R}$ and $f$ with values in
${\cal A}_v$ with $1\le n\le 2^v-1$; or vice versa $f$ with values
in ${\bf C}_{\bf i}$ or ${\bf H}_{J,K,L}$ and coefficients $a$ in
${\cal A}_v$ but with $1\le n\le 4$. Thus all such variants of
operator coefficients ${\bf a}_j$ and values of functions $f$ can be
treated by the noncommutative transform. Henceforward, we suppose
that these variants take place.
\par We suppose that $g(t)$ is an original function, that is satisfying
Conditions 1$(1-4)$. Consider at first the case of constant
coefficients ${\bf a}_j$ on a quadrant domain $Q^n$. Let $Q^n$ be
oriented so that ${\sf a}_k = - \infty $ and $b_k = + \infty $ for
each $k\le n - \kappa $; either ${\sf a}_k = - \infty $ or $b_k = +
\infty $ for each $k>n - \kappa $, where $0\le \kappa \le n$ is a
marked integer number. If conditions of Theorem 25 are satisfied,
then
$$(6)\quad {\cal F}^n(A[f](t), u;p;\zeta ) =
\sum_{|j|\le \alpha } {\bf a}_j \{ [{\sf R}_{e_1}(p)]^{j_1} [{\sf
R}_{e_2}(p)]^{j_2}... [{\sf R}_{e_n}(p)]^{j_n} {\cal F}^n (f(t)\chi
_{Q^n}(t),u; p;\zeta )$$
$$+ \sum_{1\le |(lj)|; ~ m_k+q_k+h_k=j_k; ~ 0\le m_k; ~ 0 \le q_k; ~
h_k= sign (l_kj_k); ~ q_k=0 \mbox{ for } l_kj_k=0 \mbox{, for each }
k=1,...,n; ~ (l)\in \{ 0, 1, 2 \} ^n} $$
$$ (-1)^{|(lj)|} [{\sf R}_{e_1}(p)]^{m_1}
[{\sf R}_{e_2}(p)]^{m_2}... [{\sf R}_{e_n}(p)]^{m_n} {\cal F}^{n
-|h(lj)|} (\partial ^{|q|} f(t^{(lj)})\chi _{\partial
Q^n_{(lj)}}(t^{(lj)})/\partial t_1^{q_1}...\partial t_n^{q_n},u;
p;\zeta ) \} $$  $$ = {\cal F}^n (g(t)\chi _{Q^n}(t),u;p;\zeta )$$
for $u(p,t;\zeta )$ in the ${\cal A}_r$ spherical or ${\cal A}_r$
Cartesian coordinates, where the operators ${\sf R}_{e_j}(p)$ are
given by Formulas 25$(1.1)$ or 25$(1.2)$. Here $(l)$ enumerates
faces $\partial Q^n_{(l)}$ in $\partial Q^n_{k-1}$ for $|h(l)|=k\ge
1$, so that $\partial Q^n_{k-1} = \bigcup_{|h(l)|=k} Q^n_{(l)}$,
$\partial Q^n_{(l)}\cap
\partial Q^n_{(m)}=\emptyset $ for each $(l)\ne (m)$ in accordance
with \S 25 and the notation of this section.
\par Therefore, Equation $(6)$ shows that the boundary conditions
are necessary: \par $(\partial ^{|q|} f(t^{(l)})/\partial
t_1^{q_1}...\partial t_n^{q_n})|_{\partial Q^n_{(l)}}$ for $|j|\le
\alpha $, $|(lj)|\ge 1$, ${\bf a}_j\ne 0$, $q_k=0$ for $l_kj_k=0$,
$m_k+q_k+h_k=j_k$, $~ h_k = sign (l_kj_k)$, $ ~ k=1,...,n$,
$t^{(l)}\in
\partial Q^n_{(l)}$. But $dim_{\bf R}
\partial Q^n =n-1$ for $\partial Q^n\ne \emptyset $, consequently,
$(\partial ^{|q|} f(t^{(l)})/\partial t_1^{q_1}...\partial
t_n^{q_n})|_{\partial Q^n_{(l)}}$ can be calculated if know
$(\partial ^{|\beta |} f(t^{(l)})/\partial t_{\gamma (1)}^{\beta
_1}...\partial t_{\gamma (m)}^{\beta _m})|_{\partial Q^n_{(l)}}$ for
$|\beta | = |q|$, where $\beta = (\beta _1,...,\beta _m)$,
$m=|h(l)|$, a number $\gamma (k)$ corresponds to $l_{\gamma (k)}>0$,
since $q_k=0$ for $l_k=0$ and $q_k>0$ only for $l_k j_k >0$ and
$k>n-\kappa $. That is, $t_{\gamma (1)}$,...,$t_{\gamma (m)}$ are
coordinates in ${\bf R}^n$ along unit vectors orthogonal to
$\partial Q^n_{(l)}$.
\par Take a sequence $U^k$ of sub-domains $U^k\subset
U^{k+1}\subset U$ for each $k\in {\bf N}$ so that each $U^k =
\bigcup_{l=1}^{m(k)} Q^n_{k,l}$ is the finite union of quadrants
$Q^n_{k,l}$, $m(k)\in {\bf N}$. We choose them so that each two
different quadrants may intersect only by their borders, each $U^k$
satisfies the same conditions as $U$ and
\par $(7)\quad \lim_{k\to \infty } dist (U,U^k)=0$. \par
Therefore, Equation $(6)$ can be written for more general domain $U$
also.
\par For $U$ instead of $Q^n$ we get a face $\partial U_{(l)}$ instead of
$\partial Q^n_{(l)}$ and local coordinates $\tau _{\gamma
(1)}$,...,$\tau _{\gamma (m)}$ orthogonal to $\partial U_{(l)}$
instead of $t_{\gamma (1)}$,...,$t_{\gamma (m)}$ (see Conditions
$(i-iii)$ above).
\par Thus the sufficient boundary conditions are:
\par $(5.1)$ $(\partial ^{|\beta |} f(t^{(lj)})/\partial \tau _{\gamma (1)}^{\beta _1}...\partial
\tau _{\gamma (m)}^{\beta _m})|_{\partial U_{(lj)}} = \phi _{\beta
,(lj)}(t^{(lj)})$ \\ for $|\beta | = |q|$, where $m=|h(lj)|$,
$|j|\le \alpha $, $|(lj)|\ge 1$, ${\bf a}_j\ne 0$, $q_k=0$ for
$l_kj_k=0$, $m_k+q_k+h_k=j_k$, $~ h_k = sign (l_kj_k)$, $0\le q_k\le
j_k-1$ for $k>n-\kappa $; $\phi _{\beta ,(l)}(t^{(l)})$ are known
functions on $\partial U_{(l)}$, $t^{(l)}\in \partial U_{(l)}$. In
the half-space $t_n\ge 0$ only
\par $(5.2)$ $\partial ^{\beta }f(t)/\partial t_n^{\beta }|_{t_n=0}$ \\ are
necessary for $\beta =|q|<\alpha $ and $q$ as above.
\par Depending on coefficients of the operator $A$ and the domain
$U$ some boundary conditions may be dropped, when the corresponding
terms vanish in Formula $(6)$. For example, if $A=
\partial ^2/\partial t_1\partial t_2$, $ ~ U=U_{1,1}$, $ ~ n=2$,
then $\partial f/\partial \nu |_{\partial U_0}$ is not necessary,
only the boundary condition $f|_{\partial U}$ is sufficient.
\par If $U={\bf R}^n$, then no any boundary condition appears.
Mention that
\par $(5.3)$ ${\cal F}^0(f(a);u;p;\zeta ) =f(a)e^{-u(p,a;\zeta )}$,
\\ which happens in $(6)$, when $a=t^{(l)}$ and $|h(l)|=n$.
\par Conditions in $(5.1)$ are given on disjoint for different
$(l)$ submanifolds $\partial U_{(l)}$ in $\partial U$ and partial
derivatives are along orthogonal to them coordinates in ${\bf R}^n$,
so they are correctly posed.
\par In ${\cal A}_r$ spherical coordinates due to Corollary 4.1
Equation $(6)$ with different values of the parameter $\zeta $ gives
a system of linear equations relative to unknown functions ${\sf
S}_{(m)} {\cal F}^n(f(t),u; p;\zeta )$, from which ${\cal
F}^n(f(t),u; p;\zeta )$ can be expressed through a family $$\{ {\sf
S}_{(m)} {\cal F}^n(g(t),u;p;\zeta ); ~ {\sf S}_{(m)} {\cal F}^{n
-|h(l)|} (\partial ^{|q|} f(t^{(l)})\chi _{\partial
Q^n_{(l)}}(t^{(l)})/\partial t_1^{q_1}...\partial t_n^{q_n},u;
p;\zeta ) : (m) \in {\bf Z}^n  \} $$ and polynomials of $p$, where
$\bf Z$ denotes the ring of integer numbers, where the corresponding
term ${\cal F}^{n -|h(l)|}$ is zero when $t_j^{(l)}=\pm \infty $ for
some $j$. In the ${\cal A}_r$ Cartesian coordinates there are not so
well periodicity properties generally, so the family may be
infinite. This means that ${\cal F}^n(f(t),u; p;\zeta )$ can be
expressed in the form:
$$(8)\quad {\cal F}^n(f(t),u; p;\zeta ) = \sum_{(m)} {\sf P}_{(m)}(p) {\sf S}_{(m)} {\cal
F}^n(g(t),u;p;\zeta )$$  $$ + \sum_{j,(q),(l), |(l)|\ge 1, (m)} {\sf
P}_{j,(q),(l),(m)} (p) {\sf S}_{(m)} {\cal F}^{n -|h(lj)|} (\partial
^{|q|} f(t^{(lj)})\chi _{\partial U_{(lj)}}(t^{(lj)})/\partial
t_1^{q_1}...\partial t_n^{q_n},u; p;\zeta ),$$ where ${\sf
P}_{(m)}(p)$ and ${\sf P}_{j,(q),(l),(m)}(p)$ are quotients of
polynomials of real variables $p_0, p_1,...,p_n$. The sum in $(8)$
is finite in the ${\cal A}_r$ spherical coordinates and may be
infinite in the ${\cal A}_r$ Cartesian coordinates. To the obtained
Equation $(8)$ we apply the theorem about the inversion of the
noncommutative multiparameter transform. Thus this gives an
expression of $f$ through $g$ as a particular solution of the
problem given by $(1,2,5.1)$ and it is prescribed by Formulas
6.1$(1)$ and 8.1$(1)$.
\par For ${\cal F}^n(f;u;p;\zeta )$ Conditions 8$(1,2)$ are
satisfied, since ${\sf P}_{(m)}(p)$ and ${\sf P}_{j,(q),(l),(m)}(p)$
are quotients of polynomials with real, complex or quaternion
coefficients and real variables, also $G^n$ and ${\cal F}^{n-|h(l)|}
$ terms on the right of $(6)$ satisfy them. Thus we have
demonstrated the theorem.
\par {\bf 28.1. Theorem.} {\it Suppose that
${\cal F}^n(f;u;p;\zeta )$ given by the right side of $(8)$
satisfies Conditions 8$(3)$. Then Problem $(1,2,5.1)$ has a solution
in the class of original functions, when $g$ and $\phi _{\beta
,(l)}$ are originals, or in the class of generalized functions, when
$g$ and $\phi _{\beta ,(l)}$ are generalized functions.} \par
Mention, that a general solution of $(1,2)$ is the sum of its
particular solution and a general solution of the homogeneous
problem $Af=0$. If $\phi _{\beta ,(l)}= \phi ^1_{\beta ,(l)}+\phi
^2_{\beta ,(l)}$, $g=g_1+g_2$, $f=f_1+f_2$, $Af_j=g_j$ and $f_j$ on
$\partial U_j$ satisfies $(5.1)$ with $\phi ^j_{\beta ,(l)}$, $j=1,$
$2$, then $Af=g$ and $f$ on $\partial U$ satisfies Conditions
$(5.1)$ with $\phi _{\beta ,(l)}$.
\par {\bf 28.2. Example.} We take the
partial differential operator of the second order
$$A= \sum_{h,m=1}^n {\bf a}_{h,m}\partial ^2/\partial \tau _h\partial \tau _m +
\sum_{h=1}^n \alpha _h\partial /\partial \tau _h + \omega ,$$ where
the quadratic form $a(\tau ) := \sum_{h,m} {\bf a}_{h,m} \tau _h\tau
_m$ is non-degenerate and is not always negative, because otherwise
we can consider $-A$. Suppose that ${\bf a}_{h,m}={\bf a}_{m,h}\in
{\bf R}$, $\alpha _h, \tau _h\in {\bf R}$ for each $h, m =1,...,n$,
$\omega \in {\cal A}_3$. Then we reduce this form $a(\tau )$ by an
invertible $\bf R$ linear operator $C$ to the sum of squares. Thus
$$(9)\quad A = \sum_{h=1}^n {\bf a}_h
\partial ^2/ \partial t_h^2 +\sum_{h=1}^n \beta _h \partial /\partial t_h +
\omega ,$$ where $(t_1,...,t_n) = (\tau _1,...,\tau _n)C$ with real
${\bf a}_h$ and $\beta _h$ for each $h$. If coefficients of $A$ are
constant, using a multiplier of the type $\exp (\sum_h \epsilon _h
s_h)$ it is possible to reduce this equation to the case so that if
${\bf a}_h\ne 0$, then $\beta _h=0$ (see \S 3, Chapter 4 in
\cite{rubinstb}).  Then we can simplify the operator with the help
of a linear transformation of coordinates and consider that only
$\beta _n$ may be non-zero if ${\bf a}_n=0$. For $A$ with constant
coefficients as it is well-known from algebra one can choose a
constant invertible real matrix $(c_{h,m})_{h,m =1,...,k}$
corresponding to $C$ so that ${\bf a}_h= 1$ for $h\le k_+$ and ${\bf
a}_h= - 1$ for $h>k_+$, where $0< k_+ \le n$. For $k_+ =n$ and
$\beta =0$ the operator is elliptic, for $k_+=n-1$ with ${\bf
a}_n=0$ and $\beta _n\ne 0$ the operator is parabolic, for $0<k_+<n$
and $\beta =0$ the operator is hyperbolic. Then Equation $(6)$
simplifies:
$$(10) \quad {\cal F}^n (A[f](t), u;p;\zeta ) =
\sum_{h=1}^n {\bf a}_h \{  [{\sf R}_{e_h}(p)]^2 {\cal F}^n (f(t)\chi
_{Q^n}(t),u; p;\zeta )$$
$$ + \sum_{l_h \in \{ 1, 2 \};
(l)=l_he_h } (-1)^{l_h} [ {\cal F}^{n - 1} (\partial f(t^{(l)})\chi
_{\partial Q^n_{(l)}}(t^{(l)})/\partial t_h,u; p;\zeta )$$  $$ +
[{\sf R}_{e_h}(p)] {\cal F}^{n - 1} (f(t^{(l)})\chi _{\partial
Q^n_{(l)}}(t^{(l)}),u; p;\zeta ) ] \} $$
$$ + \beta _n \{ {\cal F}^{n-1; t^{n,2}}
(f(t^{n,2})\chi _{\partial Q^n_{2e_n}}(t^{n,2}), u;p;\zeta ) - {\cal
F}^{n-1; t^{n,1}} (f(t^{n,1})\chi _{\partial Q^n_{e_n}}(t^{n,1}),
u;p;\zeta )$$
$$ + [{\sf R}_{e_n}(p)] {\cal F}^n (f(t)\chi _{Q^n}(t),u; p;\zeta )
\} + \omega {\cal F}^n (f(t)\chi _{Q^n}(t),u; p;\zeta ) = {\cal F}^n
(g(t),u;p;\zeta ) $$ in the ${\cal A}_r$ spherical or ${\cal A}_r$
Cartesian coordinates, where $e_h=(0,...,0,1,0,...,0) \in {\bf R}^n$
with $1$ on the $h$-th place, ${\sf S}_0=I$ is the unit operator,
the operators ${\sf R}_{e_h}(p)$ are given by Formulas 25$(1.1)$ or
25$(1.2)$ respectively.
\par We denote by $\delta _S(x)$
the delta function of a continuous piecewise differentiable manifold
$S$ in ${\bf R}^n$ satisfying conditions $(i-vi)$ so that
$$(\Delta)\quad \int_{{\bf R}^n} \eta (x) \delta _S(x)dx = \int_S
\eta (y)\lambda _m(dy)$$ for a continuous integrable function $\eta
(x)$ on ${\bf R}^n$, where $dim (S)=m<n$, $\lambda _m(dy)$ denotes a
volume element on the $m$ dimensional surface $S$ (see Condition
$(v)$ above). Thus we can consider a non-commutative multiparameter
transform on $\partial U$ for an original $f$ on $U$ given by the
formula: \par $(11)$ ${\cal F}^{n-1; t'}_{\partial U} (f(t')\chi
_{\partial U}(t'), u;p;\zeta ) = {\cal F}^{n;
t}(f(t)\delta _{\partial U}(t), u;p;\zeta )$. \\
Therefore, terms like ${\cal F}^{n-1}$ in $(10)$ correspond to the
boundary $\partial Q^n$. They can be simplified:
$$ (12) \quad  \beta _n \{
{\cal F}^{n-1; t^{n,2}} (f(t^{n,2})\chi _{\partial Q^n_{2e_n}}(t),
u;p;\zeta ) - {\cal F}^{n-1; t^{n,1}} (f(t^{n,1})\chi _{\partial
Q^n_{e_n}}(t), u;p;\zeta ) \}
$$
$$={\cal F}^{n-1; t'}_{\partial Q^n} (\beta (t') f(t')\chi
_{\partial Q^n}(t'), u;p;\zeta ),$$ where $\beta (t')$ is a
piecewise constant function on $\partial Q^n$ equal to $\beta _n$ on
the corresponding faces of $Q^n$ orthogonal to $e_n$ given by
condition: either $t_n={\sf a}_n$ or $t_n=b_n$; $\beta (t')=0$ is
zero otherwise.
\par If $ {\sf a}_k = - \infty $ or $b_k= +\infty $, then the
corresponding term disappears. If ${\bf R}^n$ embed into ${\cal
A}_r$ with $2^{r-1}\le n\le 2^r-1$ as ${\bf R}i_1\oplus ... \oplus
{\bf R}i_n$, then this induces the corresponding embedding $\Theta $
of $Q^n$ or $U$ into ${\cal A}_r$. This permits to make further
simplification:
$$(12.1)\quad  \sum_{h=1}^n {\bf a}_h \{ \sum_{l_h \in \{ 1, 2 \};
(l)=l_he_h } (-1)^{l_h} [ [{\sf R}_{e_h}(p)] {\cal F}^{n - 1}
(f(t^{(l)})\chi _{\partial Q^n_{(l)}}(t^{(l)}),u; p;\zeta )$$
$$ + {\cal F}^{n - 1} (\partial f(t^{(l)})\chi
_{\partial Q^n_{(l)}}(t^{(l)})/\partial t_h,u; p;\zeta ) ] \} $$
$$ = {\cal F}^{n-1}_{\partial Q^n} (a(t') (\partial
f(t')\chi _{\partial Q^n_0}(t')/\partial \nu ),u(p,t';\zeta ); p;
\zeta )$$
$$+ {\cal F}^{n-1}_{\partial Q^n} ({\sf P}(t') f(t')\chi
_{\partial Q^n_0}(t'),u; p; \zeta ),$$ where $\nu = \nu (t')$
denotes a real coordinate along an external unit normal $M(t')$ to
$\Theta (\partial U)$ at $\Theta (t')$, so that $M(t')$ is a purely
imaginary Cayley-Dickson number, $a(t')$ is a piecewise constant
function equal to ${\bf a}_h$ for the corresponding $t'$ in the face
$\partial Q^n_{l_he_h}$ with $l_h>0$; ${\sf P}(t',p) := {\sf P}(t')
:= {\sf R}_{e_h}(p)$ for $t' \in
\partial Q^n_{l_he_h}$, $h=1,...,n$, since $\sin (\psi +\pi ) = - \sin (\psi )$ and $\cos (\psi +\pi )
= - \cos (\psi )$ for each $\psi \in {\bf R}$. Certainly the
operator-valued function ${\sf P}(t')$ has a piecewise continuous
extension ${\sf P}(t)$ on $Q^n$. That is
$$(13)\quad {\cal F}^{n-1}_{\partial U} ( \xi (t') f(t')
\chi _{\partial U}(t'),u(p,t';\zeta ); p; \zeta )$$
$$ := \int_{{\bf R}^n} \xi (t) f(t) \delta _{\partial
U}(t) \exp \{ - u(p,t;\zeta ) \} dt
$$ for an integrable operator-valued function $\xi (t)$ so that
$[\xi (t) f(t)]$ is an original on $U$ whenever this integral
exists. For example, when $\xi $ is a linear combination of shift
operators ${\sf S}_{(m)}$ with coefficients $\epsilon _{(m)}(t,p)$
such that  each $\epsilon _{(m)}(t,p)$ as a function by $t\in U$ for
each $p\in W$ and $f(t)$ are originals or $f$ and $g$ are
generalized functions. For two quadrants $Q_{m,l}$ and $Q_{m,k}$
intersecting by a common face $\Upsilon $ external normals to it for
these quadrants have opposite directions. Thus the corresponding
integrals in ${\cal F}^{n-1}_{\partial Q_{m,l}}$ and ${\cal
F}^{n-1}_{\partial Q_{m,k}}$ restricted on $\Upsilon $ summands
cancel in ${\cal F}^{n-1}_{\partial (Q_{m,l}\cup Q_{m,k})}$.
\par  Using Conditions $(iv-vii)$ and the sequence $U^m$ and quadrants $Q_{m,l}^n$
outlined above we get for a boundary problem on $U$ instead of $Q^n$
the following equation:
$$(14) \quad {\cal F}^n(A[f](t), u;p;\zeta ) =
\{ \sum_{h=1}^n {\bf a}_h [{\sf R}_{e_h}(p)]^2 {\cal F}^n (f(t)\chi
_U(t),u; p;\zeta ) \} + $$  $$ \{ {\cal F}^{n-1}_{\partial U}
([\beta (t') + {\sf P}(t',p)] f(t')\chi _{\partial U_0}(t'), u; p;
\zeta ) + {\cal F}^{n-1}_{\partial U} ({\sf a}(t') (\partial
f(t')\chi _{\partial U_0}(t')/\partial \nu ),u; p; \zeta ) \} $$
$$ {\cal F}^n (\beta _n [{\sf R}_n(p)]
f(t)\chi _U(t),u; p;\zeta ) + \omega {\cal F}^n (f(t)\chi _U(t),u;
p;\zeta ) = {\cal F}^n (g(t),u;p;\zeta ) ,$$  where ${\sf
P}(t',p):={\sf P}(t') := \sum_{h=1}^n {\bf a}_h [{\sf
R}_h(p)](\partial \nu /\partial t_h)$ for each $t'\in
\partial U_0$ (see also Stokes' formula in \S XIII.3.4 \cite{zorich}
and Formulas $(14.2,14.3)$ below). Particularly, for the quadrant
domain $Q^n$ we have $a(t)={\bf a}_h$ for $t\in \partial
Q^n_{l_he_h}$ with $l_h>0$, $\beta (t)=\beta _n$ for $t\in \partial
Q^n_{l_ne_n}$ with $l_n>0$ and zero otherwise.
\par The boundary conditions are:
\par $(14.1)$ $f(t)|_{\partial U_0}=\phi (t)|_{\partial U_0}$,
$\quad (\partial f(t)/\partial \nu )|_{\partial U_0}  = \phi
_1(t)|_{\partial U_0}$. \\  The functions ${\sf a}(t)$ and ${\sf
\beta }(t)$ can be calculated from $ \{ {\bf a}_h: ~ h \} $ and
$\beta _n$ almost everywhere on $\partial U$ with the help of change
of variables from $(t_1,...,t_n)$ to $(y_1,...,y_{n-1},y_n)$, where
$(y_1,...,y_n)$ are local coordinates in $\partial U_0$ in a
neighborhood of a point $t' \in \partial U_0$, $y_n=\nu $, since
$\partial U_0$ is of class $C^1$. Consider the differential form
$\sum_{h=1}^n (-1)^{n-h} {\bf a}_h dt_1\wedge ... \wedge
\widehat{dt_h} \wedge ... \wedge dt_n = a dy_1\wedge ... \wedge
dy_{n-1}$ and its external product with $d\nu = \sum_{h=1}^n
(\partial \nu /\partial t_h)dt_h$, then
\par $(14.2)$ ${\sf a}(t)|_{\partial U_0} =
\sum_{h=1}^n {\bf a}_h (\partial \nu /\partial t_h)|_{\partial U_0}
~ ~$ and \par $(14.3)$ $\beta (t)|_{\partial U_0} = \beta _n\chi
_{U_{e_n}\cup U_{2e_n}}(\partial \nu /\partial t_n) |_{\partial U_0}$.\\
This is sufficient for the calculation of ${\cal F}^{n-1}_{\partial
U}$.
\par {\bf 28.3. Inversion procedure in the ${\cal A}_r$ spherical coordinates.}
\par When boundary conditions 28$(5.1)$ are specified, this equation
28$(6)$ can be resolved relative to ${\cal F}^n (f(t)\chi
_U(t),u(p,t;\zeta ); p;\zeta )$, particularly, for Equations
28.2$(14,14.1)$ also. The operators ${\sf S}_{e_j}$ and $T_j$ of \S
12 have the periodicity properties: ${\sf S}_{e_j}^{4+k}F(p;\zeta )=
{\sf S}_{e_j}^kF(p;\zeta ) ~$ and $~ T_j^{4+k}F(p;\zeta ) =
T_j^kF(p;\zeta ) ~$, $~ {\sf S}_{e_1}^{2+k}F(p;\zeta ) = - {\sf
S}_{e_1}^kF(p;\zeta ) ~$ and $~ T_1^{2+k}F(p;\zeta ) = -
T_1^kF(p;\zeta )$ for each positive integer number $k$ and $1\le
j\le 2^r-1$. We put \par $(6.1)$ ${\bf F}_j(p;\zeta ) := ({\sf
S}_{e_j}^4 - {\sf S}_{e_{j+1}}^4) F(p;\zeta )$ for any $1\le j\le
2^r-2$, \par $(6.2)$ ${\bf F}_{2^r-1}(p;\zeta ) := {\sf
S}_{e_{2^r-1}}^4 F(p;\zeta )$. Then from $(6)$ we get the following
equations:
$$(6.3) \sum_{|j|\le \alpha } {\bf a}_j \{ [p_0+p_1T_1]^{j_1}
[p_0+p_1T_1+p_2T_2]^{j_2}$$
$$... [p_0+p_1T_1+...+p_nT_n]^{j_n} \} |_{p_b=0 ~ \forall b>w} ~ {\bf F}_w(p;\zeta )  =
\{ - \sum_{|j|\le \alpha } {\bf a}_j $$  $$ \sum_{1\le |(lj)|; ~
m_k+q_k+h_k=j_k; ~ 0\le m_k; ~ 0 \le q_k; ~ h_k= sign (l_kj_k); ~
q_k=0 \mbox{ for } l_kj_k=0 \mbox{, for each } k=1,...,n; ~ (l)\in
\{ 0, 1, 2 \} ^n} $$
$$ (-1)^{|(lj)|} \{ [p_0+p_1T_1]^{m_1}
[p_0+p_1T_1+p_2T_2]^{m_2}... [p_0+p_1T_1+...+p_nT_n]^{m_n} \}
|_{p_b=0 ~ \forall b>w}$$
$${\cal F}^{n -|h(lj)|}_w (\partial ^{|q|} f(t^{(lj)})\chi _{\partial
Q^n_{(lj)}}(t^{(lj)})/\partial t_1^{q_1}...\partial t_n^{q_n},u;
p;\zeta ) \} + {\bf G}_w(p;\zeta ) $$ for each $w=1,...,n$, where
$F(p;\zeta ) = {\cal F}^n (f(t)\chi _{Q^n}(t),u; p;\zeta )$ and
$G(p;\zeta ) ={\cal F}^n (g(t)\chi _{Q^n}(t),u; p;\zeta )$. These
equations are resolved for each $w=1,...,n$ as it is indicated
below. Taking the sum  one gets the result
\par $(6.4)$ $F(p;\zeta )={\bf F}_1(p;\zeta )+...+{\bf F}_n(p;\zeta )$, \\ since
$\{ [\sum_{j=1}^{2^r-2} ({\sf S}_{e_j}^4 - {\sf S}_{e_{j+1}}^4)]+
{\sf S}_{e_{2^r-1}}^4 \} e^{-u(p,t;\zeta )} = {\sf
S}_{e_1}^4e^{-u(p,t;\zeta )}=e^{-u(p,t;\zeta )}$.
\par The analogous procedure is for Equation $(14)$ with the domain $U$
instead of $Q^n$.
\par From Equation $(6.3)$ or $(14)$ we get the linear equation:
$$(15)\quad \sum_{(l)} \psi _{(l)} x_{(l)} = \phi ,$$
where $\phi $ is the known function and depends on the parameter
$\zeta $, $\psi _{(l)}$ are known coefficients depending on $p$,
$x_{(l)}$ are indeterminates and may depend on $\zeta $, $l_1=0, 1$
for $h=1$, so that $x_{(l)+2e_1} = - x_{(l)}$; $l_h=0, 1, 2, 3$ for
$h>1$, where $x_{(l)+4e_h} = x_{(l)}$ for each $h>1$ in accordance
with Corollary 4.1, $(l)=(l_1,...,l_n)$. \par Acting on both sides
of $(6.3)$ or $(14)$ with the shift operators $T_{(m)}$ (see Formula
25$(SO)$), where $m_1=0,1$, $m_h=0,1,2,3$ for each $h>1$, we get
from $(15)$ a system of $2^{1+2(k-1)}$ linear equations with the
known functions $\phi _{(m)} := T_{(m)} \phi $ instead of $\phi $,
$\phi _{(0)}=\phi $:
\par $(15.1)$ $\sum_{(l)} \psi _{(l)} T_{(m)} x_{(l)} =  \phi
_{(m)}$ for each $(m)$.
\par Each such shift of $\zeta $ left coefficients $\psi _{(l)}$ intact
and $x_{(l)+(m)} = (-1)^{\eta } x_{(l')}$ with ${l'}_1 = l_1+m_1 ~
(mod ~ 2)$, ${l'}_h = l_h+m_h ~ (mod ~ 4)$ for each $h>1$, where
$\eta =1$ for $l_1+m_1-{l'}_1=2$, $\eta =2$ otherwise. This system
can be reduced, when a minimal additive group ${\cal G} := \{ (l): ~
l_1 ~ (mod ~ 2), ~ l_j ~ (mod ~ 4) ~ \forall 2\le j\le k;$ $\mbox{
generated by all } $ $(l)$ $\mbox{ with non-zero coefficients}$
$\mbox{ in Equation } $ $(15) \} $ is a proper subgroup of ${\sf
g}_2\times {\sf g}_4^{k-1}$, where ${\sf g}_h := {\bf Z}/ (h{\bf
Z})$ denotes the finite additive group for $0<h\in {\bf Z}$.
Generally the obtained system is non-degenerate for $\lambda _{n+1}$
almost all $p=(p_0,...,p_n)\in {\bf R}^{n+1}$ or in $W$, where
$\lambda _{n+1}$ denotes the Lebesgue measure on the real space
${\bf R}^{n+1}$.
\par  We consider the non-degenerate operator $A$ with real, complex
${\bf C}_{\bf i}$ or quaternion ${\bf H}_{J,K,L}$ coefficients.
Certainly in the real and complex cases at each point $p$, where its
determinate $\Delta = \Delta (p)$ is non-zero, a solution can be
found by the Cramer's rule.
\par  Generally, the system can be solved by the following algorithm.
We can group variables by $l_1$, $l_2$,...,$l_k$. For a given
$l_2,...,l_h$ and $l_1=0,1$ subtracting all other terms from both
sides of $(15)$ after an action of $T_{(m)}$ with $m_1=0,1$ and
marked $m_h$ for each $h>1$ we get the system of the form
\par $(16)$ $\alpha x_1 + \beta x_2 = {\bf b}_1$,
\par $ -\beta x_1 + \alpha x_2 ={\bf b}_2$, \\
which generally has a  unique solution for $\lambda _{n+1}$ almost
all $p$: \par $(17)$  $x_1 = (\alpha (\alpha ^2+\beta ^2)^{-1}) {\sf
b}_1 - (\beta (\alpha ^2+\beta ^2)^{-1}) {\sf b}_2)$; $x_2 = (\alpha
(\alpha ^2+\beta ^2)^{-1}){\sf b}_2 + (\beta (\alpha ^2+\beta
^2)^{-1}) {\sf b}_1$, \\  where ${\sf b}_1, {\sf b}_2\in {\cal
A}_r$, for a given set $(m_2,...m_n)$.
\par When $l_h$ are specified for each $1\le h\le k$ with $h\ne
h_0$, where $1<h_0\le k$, then the system is of the type:
\par $(18)$ $ax_1+bx_2+cx_3+dx_4={\sf b}_1$,
\par $dx_1+ax_2+bx_3+cx_4={\sf b}_2$,
\par $cx_1+dx_2+ax_3+bx_4={\sf b}_3$,
\par $bx_1+cx_2+dx_3+ax_4={\sf b}_4$, \\
where $a, b, c, d\in {\bf R}$ or ${\bf C}_{\bf i}$ or ${\bf
H}_{J,K,L}$, while ${\sf b}_1, {\sf b}_2, {\sf b}_3, {\sf b}_4\in
{\cal A}_r$. In the latter case of ${\bf H}_{J,K,L}$ it can be
solved by the Gauss' exclusion algorithm. In the first two cases of
${\bf R}$ or ${\bf C}_{\bf i}$ the solution is:
\par $(19)$ $x_j= \Delta _j/\Delta $, where
\par $\Delta = a\xi _1 - d \xi _2 + c\xi _3 - b\xi _4$,
\par $\Delta _1 = {\sf b}_1\xi _1 - {\sf b}_2\xi _2 + {\sf b}_3\xi _3 - {\sf b}_4\xi _4$,
\par $\Delta _2 = - {\sf b}_1 \xi _4 + {\sf b}_2 \xi _1 - {\sf b}_3\xi _2 + {\sf b}_4\xi
_3$,
\par $\Delta _3 = {\sf b}_1\xi _3 - {\sf b}_2\xi _4 + {\sf b}_3\xi _1 - {\sf b}_4\xi _2$,
\par $\Delta _4 = - {\sf b}_1\xi _2 + {\sf b}_2\xi _3 - {\sf b}_3\xi _4 + {\sf b}_4\xi _1$,
\par $\xi _1 = a^3+b^2c+cd^2-ac^2-2abd$,
\par $\xi _2=a^2b+bc^2+d^3-b^2d-2acd$,
\par $\xi _3 = ab^2+c^3+ad^2-a^2c-2bcd$,
\par $\xi _4=a^2d+b^3+c^2d-bd^2-2abc$.
\par Thus on each step either two or four indeterminates
are calculated and substituted into the initial linear algebraic
system that gives new linear algebraic system with a number of
indeterminates less on two or four respectively. May be pairwise
resolution on each step is simpler, because the denominator of the
type $(\alpha ^2+\beta ^2)$ should be $\lambda _{2^r}$ almost
everywhere by $p\in {\cal A}_r$ positive (see $(6)$, $(14)$ above).
This algorithm acts analogously to the Gauss' algorithm. Finally the
last two or four indeterminates remain and they are found with the
help of Formulas either $(17)$ or $(19)$ respectively. When for a
marked $h$ in $(6)$ or $(14)$ all $l_h=0 ~ (mod ~ 2)$ (remains only
$x_1$ for $h=1$, or remain $x_1$ and $x_3$ for $h>1$) or for some
$h>1$ all $l_h=0 ~ (mod ~ 4)$ (remains only $x_1$) a system of
linear equations as in $(15, 15.1)$ simplifies.
\par Thus a solution of the type prescribed by $(8)$ generally $\lambda
_{n+1}$ almost everywhere by $p\in W$ exists, where $W$ is a domain
$W=\{ p \in {\cal A}_r: a_1< Re (p) <a_{-1}, ~ p_j=0 ~\forall j>n \}
$ of convergence of the noncommutative multiparameter transform,
when it is non-void, $2^{r-1}\le n \le 2^r-1$, $Re (p)=p_0$,
$p=p_0i_0+...+p_ni_n$.
\par This domain $W$ is caused by properties of $g$ and initial
conditions on $\partial U$ and by the domain $U$ also. Generally $U$
is worthwhile to choose with its interior $Int (U)$ non-intersecting
with a characteristic surface $\phi (x_1,...,x_n)=0$, i.e. at each
point $x$ of it the condition is satisfied
\par $(CS)$ $\sum_{|j|=\alpha } {\bf a}_j(t(x)) (\partial \phi /\partial x_1)^{j_1}...
(\partial \phi /\partial x_n)^{j_n} =0  $ \\
and at least one of the partial derivatives $(\partial \phi
/\partial x_k)\ne 0$ is non-zero.
\par In particular, the boundary problem may be with the right side
$g=\varsigma f $ in $(2,2.1,14)$, where $\varsigma $ is a real or
complex ${\bf C}_{\bf i}$ or quaternion ${\bf H}_{J,K,L}$
multiplier, when boundary conditions are non-trivial. In the space
either ${\cal D}({\bf R}^n,{\cal A}_r)$ or ${\cal B}({\bf R}^n,{\cal
A}_r)$ (see \S 19) a partial differential problem simplifies,
because all boundary terms disappear. If $f\in {\cal B}({\bf
R}^n,{\cal A}_r),$ then $\{ p\in {\cal A}_r: ~ Re (p)\ge 0 \}
\subset W_f$. For $f\in {\cal D}({\bf R}^n,{\cal A}_r)$ certainly
$W_f = {\cal A}_r$ (see also \S 9).
\par {\bf 28.4. Examples.} Take  partial differential equations
of the fourth order. In this subsection the noncommutative
multiparameter transforms in ${\cal A}_r$ spherical coordinates are
considered. For
\par $(20)$ $A =
\partial ^3/\partial s_1^3 + \sum_{j=2}^n \gamma _j \partial
^4/\partial s_j^4$ \\ with constants $\gamma _j\in {\bf
H}_{J,K,L}\setminus \{ 0 \} $ on the space either ${\cal D}({\bf
R}^n,{\cal A}_r)$ or ${\cal B}({\bf R}^n,{\cal A}_r)$, where $n\ge
2$, Equation $(6)$ takes the form:
$$(21)\quad {\cal F}^n(A[f](t), u;p;\zeta ) =$$
$$ \{ p_0 (p_0^2+ 3(p_1{\sf S}_{e_1})^2) + \sum_{j=2}^n \gamma _j (p_j{\sf S}_{e_j})^4 \} {\cal
F}^n (f(t),u; p;\zeta ) + p_1 (3p_0^2 + (p_1{\sf S}_{e_1})^2) {\sf
S}_{e_1} {\cal F}^n (f(t),u; p;\zeta )$$   $$ = {\cal F}^n
(g(t),u;p;\zeta )$$  due to Corollary 4.1.  In accordance with
$(16,17)$ we get:
\par $(22)$ ${\bf F}_w(p;\zeta ) = (\alpha (\alpha ^2+\beta ^2)^{-1}){\bf G}_w
(p;\zeta ) - (\beta (\alpha ^2+\beta
^2)^{-1}) T_1 {\bf G}_w(p;\zeta ) )$ for each $w=1,...,n$, \\
where $\alpha _w = \alpha = [p_0 (p_0^2- 3p_1^2) + \sum_{j=2}^n
\gamma _j p_j^4]|_{p_b=0 ~ \forall b>w}$, $\beta _w = \beta = p_1
(3p_0^2 - p_1^2)|_{p_b=0 ~ \forall b>w}$. From Theorem 6, Corollary
6.1 and Remarks 24 we infer that:
$$(23)\quad f(t) = (2\pi )^{-n}
\int_{{\bf R}^n} F (a+p;\zeta )\exp \{ u(p,t;\zeta ) \} dp_1...dp_n
$$ supposing that the conditions of Theorem 6 and Corollary 6.1
are satisfied, where $F(p;\zeta ) = {\cal F}^n (f(t),u; p;\zeta )$.
\par If on the space either ${\cal D}({\bf R}^k,{\cal A}_r)$ or
${\cal B}({\bf R}^k,{\cal A}_r)$ an operator is as follows:
\par $(24)$ $A=
\partial ^4/\partial s_1^2\partial s_2^2  + \sum_{j=3}^n \gamma _j \partial
^4/\partial s_j^4,$ where $\gamma _j \in {\bf H}_{J,K,L}\setminus \{
0 \} $, where $n\ge 3$, then $(6)$ reads as:
\par $(25)$ ${\cal F}^n (Af(t),u; p;\zeta ) = p_2^2(p_0^2 + (p_1{\sf S}_{e_1})^2) {\sf S}^2_{e_2} {\cal
F}^n (f(t),u; p;\zeta )$ \\  $ + 2p_0p_1p_2^2 {\sf S}_{e_1} {\sf
S}^2_{e_2} {\cal F}^n (f(t)),u; p;\zeta ) + \sum_{j=3}^n \gamma _j
(p_j{\sf S}_{e_j})^4 {\cal F}^n (f(t)),u; p;\zeta )$ \\  $ = {\cal
F}^n (g(t),u; p;\zeta )$.
\par If on the same spaces an operator is:
\par $(26)$ $A=
\partial ^3/\partial s_1\partial s_2^2  + \sum_{j=3}^n \gamma _j \partial
^4/\partial s_j^4,$  where $n\ge 3$, then $(6)$ takes the form:
\par $(27)$ ${\cal F}^n (Af(t),u; p;\zeta ) = p_0p_2^2 {\sf S}_{e_2}^2 {\cal F}^n
(f(t),u; p;\zeta ) + p_1p_2^2 {\sf S}_{e_1} {\sf S}_{e_2}^2 {\cal
F}^n (f(t),u; p;\zeta ) + \sum_{j=3}^n \gamma _j (p_j{\sf
S}_{e_j})^4 {\cal F}^n (f(t),u; p;\zeta ) = {\cal F}^n (g(t),u;
p;\zeta )$.
\par To find ${\cal F}^n (f(t),u; p;\zeta )$ in $(25)$ or $(27)$
after an action of suitable shift operators $T_{(0,2,0,...,0)}$,
$T_{(1,0,...,0)}$ and $T_{(1,2,0,...,0)}$ we get the system of
linear algebraic equations:
\par $(28)$ $ax_1+bx_3+cx_4={\sf b}_1$,
\par $bx_1+cx_2+ax_3={\sf b}_2$,
\par $ax_2 - cx_3+bx_4={\sf b}_3$,
\par $-cx_1+bx_2+ax_4={\sf b}_4$ \\
with coefficients $a$, $b$ and $c$, and Cayley-Dickson numbers on
the right side ${\sf b}_1,...,{\sf b}_4\in {\cal A}_r$, where $x_1=
{\bf F}_w(p;\zeta )$, $x_2=T_1 {\bf F}_w(p;\zeta )$, $x_3= T_2^2
{\bf F}_w (p;\zeta )$, $x_4= T_1 T_2^2 {\bf F}_w(p;\zeta )$, ${\sf
b}_1= {\bf G}_w (p;\zeta )=({\cal F}^n (g(t),u; p;\zeta ))_w$, ${\sf
b}_2= T_2^2 {\bf G}_w (p;\zeta )$, ${\sf b}_3= T_1 {\bf G}_w
(p;\zeta )$, ${\sf b}_4= T_1 T_2^2 {\bf G}_w(p;\zeta )$.
Coefficients are: $a_w=a= [\sum_{j=3}^n \gamma _jp_j^4]|_{p_b=0 ~
\forall b>w}\in {\bf H}_{J,K,L}$, $b_w=b=p_2^2(p_0^2 - p_1^2)\in
{\bf R}$, $c_w=c=2p_0p_1p_2^2|_{p_b=0 ~ \forall b>w} \in {\bf R}$
for $A$ given by $(24)$; $a_w=a= [\sum_{j=3}^n \gamma _j
p_j^4]|_{p_b=0 ~ \forall b>w}\in {\bf H}_{J,K,L}$,
$b_w=b=p_0p_2^2|_{p_b=0 ~ \forall b>w}\in {\bf R}$,
$c_w=c=p_1p_2^2|_{p_b=0 ~ \forall b>w}\in {\bf R}$ for $A$ given by
$(26)$, $w=1,...,n$. If $a=0$ the system reduces to two systems with
two indeterminates $(x_1,x_2)$ and $(x_3,x_4)$ of the type described
by $(16)$ with solutions given by Formulas $(17)$. It is seen that
these coefficients are non-zero $\lambda _{n+1}$ almost everywhere
on ${\bf R}^{n+1}$. Solving this system for $a\ne 0$ we get:
\par $(29)$ ${\bf F}_w(p;\zeta ) = a^{-1} {\sf b}_1 -
[a^2-b^2+c^2)^2+4b^2c^2]^{-1} a^{-1} [ (a^2-b^2+c^2) ((c^2-b^2) {\sf
b}_1 + ab{\sf b}_2 - 2bc{\sf b}_3+ac{\sf b}_4) - 2bc
(2bc{\sf b}_1 - ac{\sf b}_2 +(c^2-b^2) {\sf b}_3  + ab{\sf b}_4)] $. \\
Finally Formula $(23)$ provides the expression for $f$ on the
corresponding domain $W$ for suitable known function $g$ for which
integrals converge. If $\gamma _j>0$ for each $j$, then $a>0$ for
each $p_3^2+...+p_w^2>0$.
\par For $(21,24)$ on a bounded domain with given boundary conditions
equations will be of an analogous type with a term on the right
${\cal F}^n (g(t),u;p;\zeta )$ minus boundary terms appearing in
$(6)$ in these particular cases.
\par For a partial differential equation
$$ (30)\quad {\bf a} (t_{n+1}) Af(t_1,...,t_{n+1}) +
\partial f(t_1,...,t_{n+1})/\partial
t_{n+1} = g(t_1,...,t_{n+1})$$ with octonion valued functions $f,
g$, where $A$ is a partial differential operator by variables
$t_1,...,t_n$ of the type given by $(2,2.1)$ with coefficients
independent of $t_1,...,t_n$, it may be simpler the following
procedure. If a domain $V$ is not the entire Euclidean space ${\bf
R}^{n+1}$ we impose boundary conditions as above in $(5.1)$. Make
the noncommutative transform ${\cal F }^{n;t_1,...,t_n}$ of both
sides of Equation $(30)$, so it takes the form:
$$ (31)\quad {\bf a} (t_{n+1}) {\cal F }^{n;t_1,...,t_n}(Af(t_1,...,t_{n+1}),u;p;\zeta ) +
\partial {\cal F }^{n;t_1,...,t_n}(f(t_1,...,t_{n+1}), u;p;\zeta )/\partial
t_{n+1}$$ $$ = {\cal F
}^{n;t_1,...,t_n}(g(t_1,...,t_{n+1}),u;p;\zeta ).$$  In the
particular case, when
\par ${\bf a}(t_{n+1}) \sum_{|j|\le \alpha } {\bf a}_j(t_{n+1})
\sum_{0\le k_1\le j_1} {{j_1}\choose {k_1}} S_{(k_1,j_2,...,j_k)}
e^{-u(p,t;\zeta )} = e^{-u(p,t;\zeta )}$ \\ for each $t_{n+1}$, $p$,
$t$ and $\zeta $, with the help of $(6,8)$ one can deduce an
expression of $F^n(p;\zeta ;t_{n+1}) :=$  ${\cal F
}^{n;t_1,...,t_n}(f(t_1,...,t_{n+1}), u;p;\zeta )$ through
$G^n(p;\zeta ;t_{n+1}) := {\cal F
}^{n;t_1,...,t_n}(g(t_1,...,t_{n+1}), u;p;\zeta )$ and boundary
terms in the following form:
$$ (32) \quad {\bf b} (p_0,...,p_n; t_{n+1}) F^n(p;\zeta ;t_{n+1}) +
\partial F^n(p;\zeta ;t_{n+1})/\partial t_{n+1} = Q(p_0,...,p_n; t_{n+1}),$$
where ${\bf b} (p_0,...,p_n; t_{n+1})$ is a real mapping and
$Q(p_0,...,p_n; t_{n+1})$ is an octonion valued function. The latter
differential equation by $t_{n+1}$ has a solution analogously to the
real case, since $t_{n+1}$ is the real variable, while ${\bf R}$ is
the center of the Cayley-Dickson algebra ${\cal A}_r$. Thus we
infer:
$$(33)\quad F^n(p;\zeta ;t_{n+1}) =
\exp \{  - \int_{\tau _0}^{t_{n+1}} {\bf b} (p_0,...,p_n; \xi )d\xi
\} $$ $$ \{ C_0 +[\int_{\tau _0}^{t_{n+1}} Q(p_0,...,p_n; \tau )
\exp \{ \int_{\tau _0}^{\tau } {\bf b} (p_0,...,p_n; \xi )d\xi \}
d\tau ] \} ,$$ since the octonion algebra is alternative and each
equation ${\bf b}x={\bf c}$ with non-zero ${\bf b}$ has the unique
solution $x= {\bf b}^{-1} {\bf c}$, where $C_0$ is an octonion
constant which can be specified by an initial condition. More
general partial differential equations as $(30)$, but with $\partial
^lf/\partial t_{n+1}^l$, $l\ge 2$, instead of $\partial f/\partial
t_{n+1}$ can be considered. Making the inverse transform $({\cal
F}^{n;t_1,...,t_n})^{-1}$ of the right side of $(33)$ one gets the
particular solution $f$.
\par {\bf 28.5. Integral kernel.}  We rewrite Equation 28$(6)$ in the
form: $$(34)\quad {\bf A}_{\cal S} {\cal F}^n (f\chi
_{Q^n},u;p;\zeta ) = {\cal F}^n(g\chi _{Q^n},u;p;\zeta ) - $$ $$
\sum_{|j|\le \alpha } {\bf a}_j \sum_{1\le |(lj)|, ~ 0\le m_k, ~
0\le q_k, ~ h_k = sign (l_kj_k), ~ m_k+q_k+h_k=j_k; ~ q_k =0 \mbox{
for } l_kj_k=0; ~ \forall k=1,...,n; ~ (l)\in \{ 0, 1, 2 \} ^n }$$
$$(-1)^{|(lj)|} {\cal S}^m {\cal F}^{n-|h(lj)|} ( \partial
^{|q|} f(t^{(lj)})/\partial t_1^{q_1}...\partial t_n^{q_n})\chi
_{\partial Q^n_{(lj)}}(t^{(lj)}), u ; p; \zeta ) ,\mbox{  where}$$
$(34.1)$ ${\cal S}_k(p) := {\cal S}_k := {\sf R}_{e_k}(p)$ \\ in the
${\cal A}_r$ spherical or ${\cal A}_r$ Cartesian coordinates
respectively (see also Formulas 25$(1.1,1.2)$), for each
$k=1,...,n$,
\par $(34.2)$ ${\cal S}^m(p) := {\cal S}^m := {\cal S}_1^{m_1}...{\cal S}_n^{m_n}$,
\par $(35)$ ${\bf A}_{\cal S} :=
\sum_{|j|\le \alpha } {\bf a}_j {\cal S}^j(p)$. \\
Then we have the integral formula:
\par $(36)$ ${\bf A}_{\cal S} {\cal F}^n (f\chi _{Q^n}, u; p;\zeta
) = \int_{Q^n} f(t) [ {\bf A}_{\cal S} \exp ( - u(p,t;\zeta ))]dt$
\\ in accordance with 1$(7)$ and 2$(4)$. Due to \S 28.3 the operator ${\bf
A}_{\cal S}$ has the inverse operator for $\lambda _{n+1}$ almost
all $(p_0,...,p_n)$ in ${\bf R}^{n+1}$. Practically, its calculation
may be cumbersome, but finding for an integral inversion formula its
kernel is sufficient. In view of the inversion Theorem 6 or
Corollary 6.1 and \S \S 19 and 20 we have
\par $(37)$ $(2\pi )^{-n} \int_{{\bf R}^n} \exp ( - u(a+p,t;\zeta ))
\exp ( u(a+p,\tau ; \zeta )) dp_1...dp_n = \delta (t;\tau ),$ where
\par $(38)$ $[\delta , f)(\tau ) = \int_{{\bf R}^n} f(t) \delta
(t;\tau ) dt_1...dt_n = f(\tau )$ \\
at each point $\tau \in {\bf R}^n$, where the original $f(\tau )$
satisfies H\"older's condition. That is, the functional $\delta
(t;\tau )$ is ${\cal A}_r$ linear. Thus the inversion of Equation
$(36)$ is:
$$(39) \quad \int_{{\bf R}^n}  ( \int_{{\bf R}^n} f(t) \chi _{Q^n}
(t) \{ [ {\bf A}_{\cal S} \exp ( - u(p+a,t;\zeta )) ] \xi
(p+a,t,\tau ;\zeta ) \} dt)dp_1...dp_n = f(\tau ),$$ so that
\par $(40)$ $[{\bf A}_{\cal S} \exp ( - u(p+a,t;\zeta )) ] \xi
(p+a,t,\tau ;\zeta ) = (2\pi )^{-n} \exp ( - u(p+a,t;\zeta )) \exp (
- u(p+a,\tau ;\zeta )) $, \\ where the coefficients of ${\bf
A}_{\cal S}$ commute with generators $i_j$ of the Cayley-Dickson
algebra ${\cal A}_r$ for each $j$. Consider at first the alternative
case, i.e. over the Cayley-Dickson algebra ${\cal A}_r$ with $r\le
3$. \par Let by our definition the adjoint operator ${\bf A}^*_{\cal
S}$ be defined by the formula
\par $(41)$ ${\bf A}^*_{\cal S} \eta (p,t;\zeta ) = \sum_{|j|\le
\alpha } {\bf a}^*_j {\cal S}^j \eta ^*(p,t;\zeta )$ for any
function $\eta : {\cal A}_r\times {\bf R}^n \times {\cal A}_r \to
{\cal A}_r$, where $t\in {\bf R}^n$, $p$ and $\zeta \in {\cal A}_r$,
${\cal S}^j\eta ^*(p,t;\zeta ) := [{\cal S}^j\eta (p,t;\zeta )]^*$.
Any Cayley-Dickson number $z\in {\cal A}_v$ can be written with the
help of the iterated exponent (see \S 3) in ${\cal A}_v$ spherical
coordinates
as \par $(42)$ $z = |z| \exp ( - u(0,0;\psi ))$, \\
where $v\ge r$, $\psi \in {\cal A}_v$, $u\in {\cal A}_v$, $Re (\psi
)=0$. Certainly the phase shift operator is isometrical:
\par $(43)$ $|T_1^{k_1}...T_n^{k_n} z| = |z|$ \\
for any $k_1,...,k_n\in {\bf R}$, since $ ~ |\exp ( - u(0,0; Im
(\psi ))| =1$ for each $\psi \in {\cal A}_v$, while
$T_1^{k_1}...T_n^{k_n} e^{-u(0,0;Im (\psi ))} = \exp \{  - u(0,0; Im
(\psi )- (k_1i_1+...+k_ni_n)\pi /2) \} $ (see \S 12).
\par In the ${\cal A}_r$ Cartesian coordinates each Cayley-Dickson
number can be presented as:
\par $(42.1)$ $z=|z|\exp (\phi M)$, where $\phi \in \bf R$ is a real parameter,
$M$ is a purely imaginary Cayley-Dickson number (see also \S 3 in
\cite{ludoyst,ludfov}).  Therefore, we deduce that
\par $(44)$ $|{\bf A}_{\cal S} \exp ( - u(p+a,t;\zeta ))| = \exp (
- (p_0+a) s_1- \zeta _0) |{\bf A}_{\cal S}\exp ( - u(Im (p),t;Im
(\zeta ) ))|$,
\\ since $\bf R$ is the center of the Cayley-Dickson algebra ${\cal
A}_v$ and $p_0, ~ a, ~ \zeta _0, ~ s_1 \in {\bf R}$, $ ~
s_1=s_1(t)$, where particularly ${\bf A}_{\cal S}1 := {\bf A}_{\cal
S}e^{-u(0,0;\zeta )}|_{\zeta =0}$ (see also Formulas 12$(3.1-3.7)$).
\par Then expressing $\xi $ from $(40)$ and using Formulas
$(41,42,42.1,44)$ we infer, that
\par $(45)$ $\xi (p,t,\tau ;\zeta ) = (2\pi )^{-n} [{\bf A}^*_{\cal
S} \exp ( - u(Im (p), t; Im (\zeta )) ]$ \par $ [\exp ( - u(Im (p),
t; Im (\zeta )) \exp ( u(p, \tau ; \zeta )) ] |{\bf A}_{\cal S}\exp
( - u(Im (p), t; Im (\zeta ))|^{-2}$,
\\ since $z^{-1} = z^*/|z|^2$ for each non-zero Cayley-Dickson
number $z\in {\cal A}_v$, $v\ge 1$, where $Im (p) =
p_1i_1+...+p_ni_n$, $p=p_0i_0+...+p_ni_n$, $p_0=Re (p)$.
\par Generally, for $r\ge 4$, Formula $(45)$ gives the integral
kernel $\xi (p,t,\tau ;\zeta )$ for any restriction of $\xi $ on the
octonion subalgebra $alg_{\bf R} (N_1,N_2,N_4)$ embedded into ${\cal
A}_r$. In view of \S 28.3 $\xi $ is unique and is defined by $(45)$
on each subalgebra $alg_{\bf R} (N_1,N_2,N_4)$, consequently,
Formula $(45)$ expresses $\xi $ by all variables $p, \xi \in {\cal
A}_r$ and $t, ~ \tau \in {\bf R}^n$. Applying Formulas $(39,45)$ and
28.2$(\Delta )$ to Equation $(34)$, when Condition 8$(3)$ is
satisfied, we deduce, that $$(46)\quad (f\chi _{Q^n})(\tau ) =
\int_{{\bf R}^n} (\int_{{\bf R}^n} g(t)\chi _{Q^n} (t) [ \exp ( -
u(p+a,t;\zeta )) \xi (p+a,t,\tau ;\zeta )]dt)dp_1...dp_n - $$
$$ \sum_{|j|\le \alpha } {\bf a}_j \sum_{1\le |(lj)|, ~ 0\le m_k, ~ 0\le q_k,
~ h_k = sign (l_kj_k); ~ m_k+q_k+h_k=j_k; ~ q_k=0 \mbox{ for }
l_kj_k=0, ~ \forall ~ k=1,...,n; ~ (l) \in \{ 0, 1, 2 \} ^n }
(-1)^{|(lj)|}$$  $$ \int_{{\bf R}^n} (\int_{\partial Q^n_{(lj)}}
[\partial ^{|q|} f(t^{(lj)}/\partial t_1^{q_1}...\partial t_n^{q_n}]
[ \{ {\cal S}^m(p) \exp ( - u(p+a,t^{(lj)};\zeta )) \}  \xi
(p+a,t^{(lj)},\tau ;\zeta )]dt^{(lj)})dp_1...dp_n , $$ where
$dim_{\bf R} \partial Q^n_{(lj)}=n-|h(lj)|$, $t^{(lj)}\in \partial
Q^n_{(lj)}$ in accordance with \S 28.1, ${\cal S}^m(p)$ is given by
Formulas $(34.1,34.2)$ above. \par For simplicity the zero phase
parameter $\zeta =0$ in $(46)$ can be taken. In the particular case
$Q^n = {\bf R}^n$ all terms with $\int_{\partial Q^n_{(lj)}}$
vanish.
\par Terms of the form $\int_{{\bf R}^n} [ \{ {\cal
S}^m(p) \exp ( - u(p+a,t;\zeta )) \}  \xi (p+a,t,\tau ;\zeta
)]dp_1...dp_n$ in Formula $(46)$ can be interpreted as left ${\cal
A}_r$ linear functionals due to Fubini's theorem and \S \S 19 and
20, where ${\cal S}^0=I$.
\par For the second order operator from $(14)$ one gets:
\par $(47)$ ${\bf A}_{\cal S} = (\sum_{h=1}^n {\bf a}_h [{\cal
S}_h(p)]^2) + \beta _n {\cal S}_n(p) +\omega $ and
$$(48)\quad  (f\chi _U)(t) = \int_{{\bf R}^n} ( \int_{{\bf R}^n}
g(t) \chi _U(t) [ \exp ( - u(p+a,t;\zeta )) \xi (p,t,\tau ;\zeta ) ]
dt) dp_1...dp_n - $$
$$\int_{{\bf R}^n} ( \int_{\partial U_0}
f(t') [ \{ (\beta (t')+ {\sf P}(t',p)) \exp ( - u(p+a,t;\zeta )) \}
\xi (p,t',\tau ;\zeta ) ] dt') dp_1...dp_n - $$
$$\int_{{\bf R}^n} ( \int_{\partial U_0}
a(t') (\partial f(t')/\partial \nu ) [ \exp ( - u(p+a,t;\zeta )) \xi
(p,t',\tau ;\zeta ) ] dt') dp_1...dp_n .$$ For a calculation of the
appearing integrals the generalized Jordan lemma (see \S \S 23 and
24 in \cite{lutsltjms}) and residues of functions at poles
corresponding to zeros $|{\bf A}_{\cal S}\exp ( - u(Im (p),t;Im
(\zeta ) ))|=0$ by variables $p_1,...,p_n$ can be used.
\par Take $g(t)=\delta (y;t)$, where $y\in {\bf R}^n$ is a
parameter, then $$(49)\quad \int_{{\bf R}^n} (\int_{{\bf R}^n}
\delta (y;t) [ \exp ( - u(p+a,t;\zeta )) \xi (p+a,t,\tau ;\zeta
)]dt)dp_1...dp_n $$  $$= \int_{{\bf R}^n} [ \exp ( - u(p+a,y;\zeta
)) \xi (p+a,y,\tau ;\zeta )]dp_1...dp_n =: {\cal E}(y;\tau )$$ is
the fundamental solution in the class of generalized functions,
where
\par $(50)$ $A_t {\cal E}(y;t) =\delta (y;t)$,
\par $(51)$ $\int_{{\bf R}^n} \delta (y;t) f(t) dt = f(y)$ \\
for each continuous function $f(t)$ from the space $L^2({\bf
R}^n,{\cal A}_r)$; $~ A_t$ is the partial differential operator as
above acting by the variables $t=(t_1,...,t_n)$ (see also \S \S 19,
20 and 33-35).
\par {\bf 29. The decomposition theorem of partial differential operators
over the Cayley-Dickson algebras.}
\par We consider a partial
differential operator of order $u$: $$(1)\quad Af(x)= \sum_{|\alpha
|\le u} {\bf a}_{\alpha }(x)\partial ^{\alpha } f(x),$$ where
$\partial ^{\alpha } f=\partial ^{|\alpha |}f(x)/\partial
x_0^{\alpha _0}...\partial x_n^{\alpha _n}$, $x=x_0i_0+...x_ni_n$,
$x_j\in {\bf R}$ for each $j$, $1\le n=2^r-1$, $\alpha = (\alpha
_0,...,\alpha _n)$, $|\alpha |=\alpha _0+...+\alpha _n$, $0\le
\alpha _j\in {\bf Z}$. By the definition this means that the
principal symbol
$$(2)\quad A_0 := \sum_{|\alpha |= u} {\bf a}_{\alpha }(x)\partial
^{\alpha }$$ has $\alpha $ so that $|\alpha |=u$ and ${\bf
a}_{\alpha }(x)\in {\cal A}_r$ is not identically zero on a domain
$U$ in ${\cal A}_r$. As usually $C^k(U,{\cal A}_r)$ denotes the
space of $k$ times continuously differentiable functions by all real
variables $x_0,...,x_n$ on $U$ with values in ${\cal A}_r$, while
the $x$-differentiability corresponds to the super-differentiability
by the Cayley-Dickson variable $x$.
\par Speaking about locally constant or locally differentiable
coefficients we shall undermine that a domain $U$ is the union of
subdomains $U^j$ satisfying conditions 28$(D1,i-vii)$ and $U^j\cap
U^k = \partial U^j\cap \partial U^k$ for each $j\ne k$. All
coefficients ${\bf a}_{\alpha }$ are either constant or
differentiable of the same class on each $Int (U^j)$ with the
continuous extensions on $U^j$. More generally it is up to a $C^u$
or $x$-differentiable diffeomorphism of $U$ respectively.
\par  If an operator $A$ is of the odd order $u=2s-1$, then an operator $E$ of
the even order $u+1=2s$ by variables $(t,x)$ exists so that \par
$(3)$ $Eg(t,x)|_{t=0}=Ag(0,x)$ for any $g\in C^{u+1}([c,d]\times
U,{\cal A}_r)$, where $t\in [c,d]\subset {\bf R}$, $c\le 0<d$, for
example, $Eg(t,x) =
\partial (tAg(t,x))/\partial t$.
\par Therefore, it remains the case of the operator $A$ of the even order $u=2s$.
Take $z=z_0i_0+...+z_{2^v-1}i_{2^v-1}\in {\cal A}_v$, $z_j\in {\bf
R}$. Operators depending on a less set $z_{l_1},...,z_{l_n}$ of
variables can be considered as restrictions of operators by all
variables on spaces of functions constant by variables $z_s$ with
$s\notin \{ l_1,...,l_n \} $.
\par {\bf Theorem.} {\it Let $A=A_u$ be a partial differential operator
of an even order $u=2s$ with locally constant or variable $C^s$ or
$x$-differentiable on $U$ coefficients ${\bf a}_{\alpha }(x)\in
{\cal A}_r$ such that it has the form \par $(4)$ $Af =
c_{u,1}(B_{u,1}f) +...+ c_{u,k}(B_{u,k}f)$, where each
\par $(5)$ $B_{u,p}=B_{u,p,0}+Q_{u-1,p}$ \\ is a partial differential operator
by variables
$x_{m_{u,1}+...+m_{u,p-1}+1}$,...,$x_{m_{u,1}+...+m_{u,p}}$ and of
the order $u$, $m_{u,0}=0$, $c_{u,k}(x)\in {\cal A}_r$ for each $k$,
its principal part \par $(6)$ $B_{u,p,0}= \sum_{|\alpha |=s} {\bf
a}_{p,2\alpha }(x)\partial ^{2\alpha }$ \\ is elliptic with real
coefficients ${\bf a}_{p,2\alpha }(x)\ge 0$, either $0\le r\le 3$
and $f\in C^u(U,{\cal A}_r)$, or $r\ge 4$ and $f\in C^u(U,{\bf R})$.
Then three partial differential operators $\Upsilon ^s$ and
$\Upsilon _1^s$ and $Q$ of orders $s$ and $p$ with $p\le u-1$ with
locally constant or variable $C^s$ or $x$-differentiable
correspondingly on $U$ coefficients with values in ${\cal A}_v$
exist, $r\le v$, such that
\par $(7)$ $Af=\Upsilon ^s(\Upsilon _1^sf) +Qf$.}
\par {\bf Proof.}  Certainly we have $ord Q_{u-1,p}\le u-1$,
$ord (A-A_0) \le u-1$. We choose the following operators:
$$(8)\quad \Upsilon ^s f(x) = \sum_{p=1}^k \sum_{|\alpha |\le s, ~
\alpha _q = 0 \forall q<(m_{u,1}+...+m_{u,p-1}+1) \mbox{ and }
q>(m_{u,1}+...+m_{u,p})} (\partial ^{\alpha } f(x)) [w_p^* \psi _{p,
\alpha }]\mbox{ and}$$
$$(9)\quad \Upsilon ^s_1 f(x) =
\sum_{p=1}^k \sum_{|\alpha |\le s, ~ \alpha _q = 0 \forall
q<(m_{u,1}+...+m_{u,p-1}+1) \mbox{ and } q>(m_{u,1}+...+m_{u,p})}
(\partial ^{\alpha } f(x)) [w_p\psi _{p,\alpha }^*],$$ where
$w_p^2=c_{u,p}$ for all $p$ and ${\psi }_{p,\alpha }^2(x)= - {\bf
a}_{p,2\alpha }(x)$ for each $p$ and $x$, $w_p\in {\cal A}_r$,
${\psi }_{p,\alpha }(x)\in {\cal A}_{r,v}$ and ${\psi }_{p,\alpha
}(x)$ is purely imaginary for ${\bf a}_{p,2\alpha }(x)>0$ for all
$\alpha $ and $x$, $Re (w_p Im (\psi _{p,\alpha }))=0$ for all $p$
and $\alpha $, $Im (x) = (x-x^*)/2$, $v>r$. Here ${\cal A}_{r,v} =
{\cal A}_v/{\cal A}_r$ is the real quotient algebra. The algebra
${\cal A}_{r,v}$ has the generators $i_{j2^r}$, $j=0,...,2^{v-r}-1$.
A natural number $v$ so that $2^{v-r} -1\ge \sum_{p=1}^k
\sum_{q=0}^u {{m_p+q-1}\choose q}$ is sufficient, where ${m\choose
q} = m!/(q!(m-q)!)$ denotes the binomial coefficient,
${{m+q-1}\choose q}$ is the number of different solutions of the
equation $\alpha _1+...+\alpha _m =q$ in non-negative integers
$\alpha _j$. We have either $\partial ^{\alpha + \beta }f\in {\cal
A}_r$ for $0\le r\le 3$ or $\partial ^{\alpha + \beta }f\in {\bf R}$
for $r\ge 4$. Therefore, we can take $\psi _{p,\alpha }(x) \in
i_{2^rq}{\bf R}$, where $q=q(p,\alpha )\ge 1$, $ ~ ~ q(p^1,\alpha
^1)\ne q(p,\alpha )$ when $(p,\alpha )\ne (p^1,\alpha ^1)$.
\par Thus Decomposition $(7)$ is valid due to the following. For $b=
\partial ^{\alpha +\beta }f(z)$ and ${\bf l} = i_{2^rp}$ and $w\in
{\cal A}_r$ one has the identities: \par $(10)$ $(b(w{\bf l}))
(w^*{\bf l}) = ((wb){\bf l})(w^*{\bf l}) = - w(wb) = - w^2b$ and
\par $(11)$ $(((b{\bf l})w^*){\bf l})w = (((bw){\bf l}){\bf l})w = - (bw)w =
- bw^2$ in the considered here cases, since ${\cal A}_r$ is
alternative for $r\le 3$ while ${\bf R}$ is the center of the
Cayley-Dickson algebra (see Formulas $(M1,M2)$ in the introduction).
\par This decomposition of the operator $A_{2s}$ is generally up to
a partial differential operator of order not greater, than $(2s-1)$:
\par $(12)$ $Qf(x) = \sum_{|\alpha |\le s, |\beta | \le s;
\gamma \le \alpha , \epsilon \le \beta , |\gamma + \epsilon
|>0}[\prod _{j=0}^{2^v-1} {{\alpha _j}\choose {\gamma _j}} {{\beta
_j}\choose {\epsilon _j}} ] (\partial ^{\alpha +\beta - \gamma -
\epsilon } f(x))$ \\ $ [(\partial ^{\gamma } {\eta }_{\alpha }(x))
((\partial ^{\epsilon } {\eta }_{\beta }^1(x)]$, \\ where operators
$\Upsilon ^s$ and $\Upsilon ^s_1$ are already written in accordance
with the general form \par $(13)$ $\Upsilon ^sf(x) = \sum_{|\alpha
|\le s} (\partial ^{\alpha }f(x)) \eta _{\alpha }(x)$; \par  $(14)$
$\Upsilon ^s_1 f(x) = \sum_{|\beta |\le s} (\partial ^{\beta }f(x))
\eta _{\beta }^1(x)$.
\par When $A$ in $(3)$ is with constant coefficients, then the
coefficients $w_p$ and $\psi _{p,\alpha }$ for $\Upsilon ^m$ and
$\Upsilon ^m_1$ can also be chosen constant and $Q=0$.
\par {\bf 30. Corollary.} {\it Let suppositions of Theorem 29 be
satisfied. Then a change of variables locally constant or variable
$C^1$ or $x$-differentiable on $U$ correspondingly exists so that
the principal part $A_{2,0}$ of $A_{2}$ becomes with constant
coefficients, when ${\bf a}_{p,2\alpha }>0$ for each $p$, $\alpha $
and $x$.}
\par {\bf 31. Corollary.} {\it If two operators $E=A_{2s}$ and $A=A_{2s-1}$
are related by Equation 29$(3)$, and $A_{2s}$ is presented in
accordance with Formulas 29$(4,5)$, then three operators $\Upsilon
^s$, $\Upsilon ^{s-1}$ and $Q$ of orders $s$, $s-1$ and $2s-2$ exist
so that
\par $(1)$ $A_{2s-1}=\Upsilon ^s\Upsilon ^{s-1} +Q$.}
\par {\bf Proof.} It remains to verify that $ord (Q)\le 2s-2$ in
the case of $A_{2s-1}$, where $Q= \{ \partial (tA_{2s-1})/ \partial
t - \Upsilon ^s\Upsilon _1^s \} |_{t=0}$. Indeed, the form $\lambda
(E)$ corresponding to $E$ is of degree $2s-1$ by $x$ and each
addendum of degree $2s$ in it is of degree not less than $1$ by $t$,
consequently, the product of forms $\lambda (\Upsilon _s) \lambda
(\Upsilon ^s_1)$ corresponding to $\Upsilon ^s$ and $\Upsilon ^s_1$
is also of degree $2s-1$ by $x$ and each addendum of degree $2s$ in
it is of degree not less than $1$ by $t$. But the principal parts of
$\lambda (E)$ and $\lambda (\Upsilon _s) \lambda (\Upsilon ^s_1)$
coincide identically by variables $(t,x)$, hence $ord ( \{ E -
\Upsilon ^s\Upsilon _1^s \} |_{t=0}) \le 2s-2$. Let $a(t,x)$ and
$h(t,x)$ be coefficients from $\Upsilon ^s_1$ and $\Upsilon ^s$.
Using the identities
\par $a(t,x) \partial _t \partial ^{\gamma } tg(x) = a(t,x)
\partial ^{\gamma } g(x)$ and \par $h(t,x) \partial ^{\beta } \partial _t
[a(t,x) \partial ^{\gamma } g(x)] = h(t,x)
\partial ^{\beta } [(\partial _t a(t,x))
\partial ^{\gamma } g(x)] $ \\ for any functions
$g(x)\in C^{2s-1}$ and $a(t,x)\in C^s$, $ord [(h(t,x) \partial
^{\beta }), (a(t,x) \partial ^{\gamma })]|_{t=0}\le 2s-2$, where
$\partial _t =\partial /\partial t$, $|\beta |\le s-1$, $|\gamma
|\le s$, $[A,B] := AB-BA$ denotes the commutator of two operators,
we reduce $(\Upsilon ^s\Upsilon ^s_1 + Q_1)|_{t=0}$ from Formula
29$(7)$ to the form prescribes by equation $(1)$.
\par \par {\bf 32.} We consider operators of the form:
\par $(1)$ $(\Upsilon ^k + \beta I_r) f(z) := \{ \sum_{0<|\alpha |\le k} (\partial ^{\alpha }
f(z) {\eta }_{\alpha }(z) \} + f(z)\beta (z)$, \\ with $\eta
_{\alpha }(z)\in {\cal A}_v$, $\alpha = (\alpha _0,...,\alpha
_{2^r-1})$, $0\le \alpha _k\in {\bf N}$ for each $k$, $|\alpha
|=\alpha _0+...+\alpha _{2^r-1}$, $\beta I_r f(z) := f(z) \beta $,
\par $\partial ^{\alpha }f(z) := \partial ^{|\alpha |} f(z)/\partial
z_0^{\alpha _0}...\partial z_{2^r-1}^{\alpha _{2^r-1}}$, $2\le r \le
v<\infty $, $\beta (z)\in {\cal A}_v$, $z_0,...,z_{2^r-1}\in {\bf
R}$, $z=z_0i_0+...+z_{2^r-1} i_{2^r-1}$.
\par {\bf Proposition.} {\it The operator $(\Upsilon ^k+\beta
)^*(\Upsilon ^k+\beta )$ is elliptic on the space $C^{2k}({\bf
R}^{2^r},{\cal A}_v)$.}
\par {\bf Proof.} We establish the identity \par $(2)$
$(ay)z^* + (az)y^* = a 2 Re (yz^*)$ \\ for any $a, y, z\in {\cal
A}_v$. It is sufficient to prove Equality $(2)$ for any three basic
generators of the Cayley-Dickson algebra ${\cal A}_v$, since the
real field ${\bf R}$ is its center, while the multiplication in
${\cal A}_v$ is distributive $(a+y)z=az+yz$ and $((\alpha a) (\beta
y)) (\gamma z^*) = (\alpha \beta \gamma ) ((ay)z^*)$ for all $\alpha
, \beta , \gamma \in {\bf R}$ and $a, y, z \in {\cal A}_v$. If
$a=i_0$, then $(2)$ is evident, since $yz^* + zy^* = yz^* + (yz^*)^*
= 2 Re (yz^*)$. If $y=i_0$, then $(ay)z^*+ (az)y^*= az^* + az = a 2
Re (z)= a 2 Re (yz^*)$. Analogously for $z=i_0$. \par For three
purely imaginary generators $i_p, i_s, i_k$ consider the minimal
Cayley-Dickson algebra $\Phi = alg_{\bf R} (i_p, i_s, i_k)$ over the
real field generated by them. If it is associative, then it is
isomorphic with either the complex field $\bf C$ or the quaternion
skew field $\bf H$, so that $(ay)z^* + (az)y^* = a(yz^*+zy*) = a 2
Re (yz^*)$. \par If the algebra $\Phi $ is isomorphic with the
octonion algebra, then we use Formulas $(M1,M2)$ from the
introduction for either $a, y\in {\bf H}$ and $z={\bf l}$ or $a,
z\in {\bf H}$ and $y={\bf l}$. This gives $(2)$ in all cases, since
the algebra $alg_{\bf R} (i_p,i_s)$ with two basic generators $i_p$
and $i_s$ is always associative. Particularly, if $y=i_s\ne z=i_k$,
$~ s, ~ k \ge 1$, then the result in $(2)$ is zero.
\par Using $(2)$ we get more generally, that
\par $(3)$ $((ay)z^*)b^* + ((az)y^*)b^* = (a 2 Re (yz^*))b^*=
(ab^*) 2 Re (yz^*) ,$ \\ consequently,
\par $(4)$ $((ay)z^*)b^* + ((az)y^*)b^* + ((by)z^*)a^* + ((bz)y^*)a^*=
4 Re (ab^*) Re (yz^*) $ \\
for any Cayley-Dickson numbers $a, b, y, z\in {\cal A}_v$. In view
of Formulas $(1,4)$ the form corresponding to the principal symbol
of the operator $(\Upsilon ^k+\beta )^*(\Upsilon ^k+\beta )$ is with
real coefficients, of degree $2k$ and non-negative definite,
consequently, the operator $(\Upsilon ^k+\beta )^*(\Upsilon ^k+\beta
)$ is elliptic.

\par {\bf 33. Fundamental solutions.}
Let either $Y$ be a real $Y={\cal A}_v$ or complexified $Y=({\cal
A}_v)_{\bf C}$ or quaternionified $Y=({\cal A}_v)_{\bf H}$
Cayley-Dickson algebra (see \S 28). Consider the space ${\cal
B}({\bf R}^n,Y)$ (see \S 19) supplied with a topology in it is given
by the countable family of semi-norms
\par $(1)$ $p_{\alpha , k} (f) := \sup_{x\in {\bf R}^n} |(1+|x|)^k \partial ^{\alpha
}f(x)|$, \\ where $k=0, 1, 2,...$; $\alpha = (\alpha _1,...,\alpha
_n)$, $0\le \alpha _j\in {\bf Z}$. On this space we take the space
${\cal B}'({\bf R}^n,Y)_l$ of all $Y$ valued continuous generalized
functions (functionals) of the form
\par $(2)$ $f=f_0i_0+...+f_{2^v-1}i_{2^v-1}$ and
$g=g_0i_0+...+g_{2^v-1}i_{2^v-1}$, where $f_j$ and $g_j\in {\cal
B}'({\bf R}^n,Y)$,  with restrictions on ${\cal B}({\bf R}^n,{\bf
R})$ being real or ${\bf C}_{\bf i}$ or ${\bf H}_{J,K,L}$ -valued
generalized functions $f_0,...,f_{2^v-1}, g_0,...,g_{2^v-1}$
respectively. Let $\phi = \phi _0i_0+...+\phi _{2^v-1}i_{2^v-1}$
with $\phi _0,...,\phi _{2^v-1}\in {\cal B}({\bf R}^n,{\bf R})$,
then
\par $(3)$ $[f,\phi ) = \sum_{k,j=0}^{2^v-1} [f_j,\phi _k) i_ki_j$.
We define their convolution  as
\par $(4)$ $[f*g,\phi ) = \sum_{j,k=0}^{2^v-1} ([f_j*g_k,\phi )
i_j)i_k$ for each $\phi \in {\cal B}({\bf R}^n,Y)$. As usually
\par $(5)$ $(f*g)(x) = f(x-y)*g(y) = f(y)*g(x-y)$ \\ for all $x, y \in {\bf R}^n$
due to $(4)$, since the latter equality $(5)$ is satisfied for each
pair $f_j$ and $g_k$. Thus a solution of the equation
\par $(6)$ $(\Upsilon ^s +\beta )f = g $ in ${\cal B}({\bf R}^n,Y)$ or in
the space ${\cal B}'({\bf R}^n,Y)_l$ is: \par $(7)$  $f = {\cal
E}_{\Upsilon ^s+\beta }*g$, where ${\cal E}_{\Upsilon ^s +\beta }$
denotes a fundamental solution of the equation \par $(8)$ $(\Upsilon
^s +\beta ){\cal E}_{\Upsilon +\beta }=\delta $, $(\delta ,\phi
)=\phi (0)$. The fundamental solution of the equation
\par $(9)$ $A_0 {\cal V} = \delta $ with $A_0 = (\Upsilon ^s +\beta ) (\Upsilon ^{s_1}_1+\beta _1)$ \\
using Equalities 32$(2-4)$ can be written as the convolution
$$(10)\quad {\cal V} =: {\cal V}_{A_0} = {\cal E}_{\Upsilon ^s +\beta
} * {\cal E}_{\Upsilon ^{s_1}_1+\beta _1}.$$
\par More generally we can consider the equation
\par $(11)$ $A f=g$ with $A=A_0 + (\Upsilon _2+\beta _2)$, \\ where
$A_0=(\Upsilon +\beta ) (\Upsilon _1 +\beta _1)$, $\Upsilon , ~
\Upsilon _1, ~ \Upsilon _2$ are operators of orders $s$, $s_1$ and
$s_2$ respectively given by 32$(1)$ with $z$-differentiable
coefficients. For $\Upsilon _2+\beta _2=0$ this equation was solved
above. Suppose now, that the operator $\Upsilon _2+\beta _2$ is
non-zero.
\par To solve Equation $(11)$ on a domain $U$ one can write it as the system:
\par $(12)$  $(\Upsilon _1+\beta _1)f =  g_1$, $(\Upsilon +\beta )g_1 = g - (\Upsilon _2+\beta _2)f$. \\
Find at first a fundamental solution ${\cal V}_A$ of Equation $(11)$
for $g=\delta $. We have: \par $(13)$ $f={\cal E}_{\Upsilon _1+\beta
_1}*g_1 = {\cal E}_{\Upsilon _2+\beta _2}*(g- (\Upsilon  +\beta
)g_1)$, consequently, \par $(13.1)$ ${\cal E}_{\Upsilon _1+\beta
_1}*g_1 + {\cal E}_{\Upsilon _2+\beta _2}*((\Upsilon  +\beta )g_1) =
{\cal E}_{\Upsilon _2+\beta _2}*g$. \\ In accordance with $(3-5)$
and 32$(1)$ the identity is satisfied: $[{\cal E}_{\Upsilon _2+\beta
_2}*((\Upsilon  +\beta )g_1),\phi _0) = [ (\Upsilon + \beta ) ({\cal
E}_{\Upsilon _2+\beta _2}*g_1), \phi _0)$. Thus $(13)$ is equivalent
to
\par $(14)$ ${\cal E}_{\Upsilon _1+\beta _1}*g_1 + (\Upsilon +\beta ) ({\cal
E}_{\Upsilon _2+\beta _2}*g_1) = {\cal E}_{\Upsilon _2+\beta _2}$ \\
for $g=\delta $, since ${\cal E}_{\Upsilon _2+\beta _2}*\delta =
{\cal E}_{\Upsilon _2+\beta _2}$. \par We consider the Fourier
transform $F$ by real variables with the generator ${\bf i}$
commuting with $i_j$ for each $j=0,...,2^v-1$ such that
\par $(F1)$ $(Fg)(y) = \int_{{\bf R}^n} e^{- {\bf i}(y,x)}
g(x)dx_1...dx_n$ \\ for any $g\in L^1({\bf R}^n,{\cal A}_v)$, i.e.
$\int_{{\bf R}^n}|g(x)|dx_1...dx_n<\infty $, where $g: {\bf R}^n\to
Y$ is an integrable function, $(y,x)=x_1y_1+...+x_ny_n$,
$x=(x_1,...,x_n)\in {\bf R}^n$, $x_j\in {\bf R}$ for every $j$. The
inverse Fourier transform is:
\par $(F2)$ $(F^{-1}g)(y) = (2\pi )^{-n} \int_{{\bf R}^n} e^{ {\bf i}(y,x)}
g(x)dx_1...dx_n$.
\par For a generalized function $f$ from the space
${\cal B}'({\bf R}^n,Y)_l$ its Fourier transform is defined by the
formula
\par $(F3)$ $(Ff,\phi ) = (f,F\phi )$, $~ ~ (F^{-1}f,\phi ) =
(f,F^{-1}\phi )$. \par In view of $(2-5)$ the Fourier transform of
$(14)$ gives:
\par $(15)$ $ [F({\cal E}_{\Upsilon _1+\beta _1})] [F(g_1)]+
\sum_{j=0}^{2^v-1} [F((\Upsilon +\beta )_j {\cal
E}_{\Upsilon _2+\beta _2})]  [F(g_1)] i_j =F({\cal E}_{\Upsilon _2+\beta _2})$ \\
for $g=\delta $.  With generators $i_0,...,i_{2^v-1}, i_0{\bf
i},...,i_{2^v-1}{\bf i}$ the latter equation gives the linear system
of $2^{v+1}$ equations over the real field, or $2^{v+2}$ equations
when $Y=({\cal A}_v)_{\bf H}$. From it $F(g_1)$ and using the
inverse transform $F^{-1}$ a generalized function $g_1$ can be
found, since $F(g) = F(g_0)i_0+...+F(g_{2^v-1})i_{2^v-1}$ and
$F^{-1}(g) = F^{-1}(g_0)i_0+...+F^{-1}(g_{2^v-1})i_{2^v-1}$ (see
also the Fourier transform of real and complex generalized functions
in \cite{gelshil,vladumf}). Then \par $(16)$ ${\cal V}_A = {\cal
E}_{\Upsilon _1+\beta _1}*g_1$ and $f={\cal V}_A*g$ gives the
solution of $(11)$, where $g_1$ was calculated from $(15)$. \par Let
$\pi ^v_r : ({\cal A}_v)_{\bf H}\to ({\cal A}_r)_{\bf H}$ be the
$\bf R$-linear projection operator defined as the sum of projection
operators $\pi _0+...+\pi _{2^r-1}$, where $\pi _j:  ({\cal
A}_v)_{\bf H}\to {\bf H}i_j$, \par $(17)$ $\pi _j(h)=h_ji_j$, $ ~ h
= \sum_{j=0}^{2^v-1} h_ji_j$, $h_j\in {\bf H}_{J,K,L}$, that gives
the corresponding restrictions when $h_j\in {\bf C}_{\bf i}$ or
$h_j\in {\bf R}$ for $j=0,...,2^r-1$. Indeed, Formulas 2$(5,6)$ have
the natural extension on $({\cal A}_v)_{\bf H}$, since the
generators $J, ~ K$ and $L$ commute with $i_j$ for each $j$.
\par Finally, the restriction from the domain in ${\cal A}_v$ onto the initial domain
of real variables in the real shadow and the extraction of $\pi ^v
_r\circ f\in {\cal A}_r$ with the help of Formulas 2$(5,6)$ gives
the reduction of a solution from ${\cal A}_v$ to ${\cal A}_r$.  \par
Theorems 29, Proposition 32 and Corollaries 30, 31 together with
formulas of this section provide the algorithm for subsequent
resolution of partial differential equations for $s, s-1,...,2$,
because principal parts of operators $A_2$ on the final step are
with constant coefficients. A residue term $Q$ of the first order
can be integrated along a path using a non-commutative line
integration over the Cayley-Dickson algebra \cite{ludoyst,ludfov}.

\par {\bf 34. Multiparameter transforms of generalized functions.}
\par If $\phi \in {\cal B}({\bf R}^n,Y)$ and $g\in {\cal B}'({\bf R}^n,Y)_l$
(see \S \S 19 and 33) we put
\par $(1)$ $\sum_{j=0}^{2^v-1} [{\cal F}^n(g_j;u;p;\zeta ),\phi )i_j := \sum_{j=0}^{2^v-1}
[g_j,{\cal F}^n(\phi ;u;p;\zeta ))i_j$ or shortly
\par $(2)$ $\sum_{j=0}^{2^v-1}[g_je^{-u(p;t;\zeta)},\phi )i_j = \sum_{j=0}^{2^v-1} [g_j,
\phi e^{-u(p;t;\zeta)}) i_j$. \\ If the support $supp (g)$ of $g$ is
contained in a domain $U$, then it is sufficient to take a base
function $\phi $ with the restriction $\phi |_U \in {\cal B}(U,Y)$
and any $\phi |_{{\bf R}^n\setminus U}\in C^{\infty }$.
\par {\bf 34.1. Remark.} It is possible to use Theorem 29,
Corollaries 30 and 31, Proposition 32 and \S 33 for solutions of
definite differential equations with variable coefficients. For this
purpose one can present an operator $A$ as the composition
$A=\Upsilon \Upsilon _1 + Q$, where $ord (A) = ord (\Upsilon ) + ord
(\Upsilon _1)$, $ord (Q) \le ord (A) -1$, $\Upsilon $ and $\Upsilon
_1$ are operators with variable coefficients, $\Upsilon ^*\Upsilon $
and $\Upsilon ^*_1\Upsilon _1$ are elliptic operators with constant
coefficients of their principal symbols at least. Then use Formulas
33$(1-16)$ to find fundamental solutions ${\cal E}_{\Upsilon }$,
${\cal E}_{\Upsilon _1}$ and ${\cal E}_A$ or iterate this procedure
(see also \S 35). A generalization of Feynman's formula over the
Cayley-Dickson algebras for the second order partial differential
operators with the first order addendum $Q$ with variable
coefficients from \cite{ludmmas09} also can be used.

\par {\bf 35. Examples.} Let \par $(1)$ $Af(t)=\sum_{j=1}^n (\partial ^2f(t)/\partial t_j^2)c_j$
\\ be the operator with constant coefficients $c_j\in {\cal A}_r$,
$|c_j|=1$, by the variables $t_1,...,t_n$, $n\ge 2$. We suppose that
$c_j$ are such that the minimal subalgebra $alg_{\bf R}(c_j,c_k)$
containing $c_j$ and $c_k$ is alternative for each $j$ and $k$ and
$|(...(c_1^{1/2}c_2^{1/2})...)c_n^{1/2}|=1$. Since
\par $(2)$ $\partial f(t)/\partial t_j = \sum_{k=1}^n (\partial
f(t(s))/\partial s_k) (\partial s_k/\partial t_j) = \sum_{k=1}^j
\partial f(t(s))/\partial s_k$, the operator $A$ takes the form
\par $(3)$ $Af=\sum_{j=1}^n (\sum_{1\le k, b \le j} (\partial ^2
f(t(s))/\partial s_k\partial s_b)) c_j,$ \\ where $s_j=t_j+...+t_n$
for each $j$. Therefore, by Theorem 12 and Formulas 25$(SO)$ and
28$(6)$ we get:
\par $(4)$ ${\cal F}^n (Af;u;p;\zeta ) = \sum_{j=1}^n
\{ [{\sf R}_{e_j}(p)]^2 F^n_u (p;\zeta ) \} c_j $ for $u(p,t;\zeta
)$ either in ${\cal A}_r$ spherical or ${\cal A}_r$ Cartesian
coordinates with the corresponding operators ${\sf R}_{e_j}(p)$
(see also Formulas 25$(1.1,1.2)$). \\
On the other hand,
\par $(5)$ ${\cal F}^n(\delta ;u;p;\zeta ) =e^{-u(p,0;\zeta
)}=e^{-u(0,0;\zeta )}$ in accordance with Formula 20$(2)$. The delta
function $\delta (t)$ is invariant relative to any invertible linear
operator $C: {\bf R}^n\to {\bf R}^n$ with the determinant $|\det
(C)|=1$, since
$$\int_{{\bf R}^n} \delta (Cx) \phi (x)dx = \int_{{\bf R}^n} \delta
(y) \phi (C^{-1}y) |\det (C)| dy = \phi (C^{-1}0) = \phi (0).$$ Thus
\par $(5)$ ${\cal F}^n (C(Af);u;p;\zeta ) = {\cal F}^n (Af;u;p;\zeta
)$ \\ for any Fundamental solution $f$, where $Cg(t) := g(Ct)$,
$Af=\delta $.  If $C: {\bf R}^n\to {\bf R}^n$ is an invertible
linear operator and $\xi =Ct$, $q=Cp$, $\zeta ' =C\zeta $, then
$t=C^{-1}\xi $, $p=C^{-1}q$ and $\zeta  =C^{-1} \zeta '$. In the
multiparameter noncommutative transform ${\cal F}^n$ there are the
corresponding variables $(t_j,p_j,\zeta _j)$. This is accomplished
in particular for the operator $C(t_1,...,t_n)=(s_1,...,s_n)$. The
operator $C^{-1}$ transforms the right side of Formula $(4)$, when
it is written in the ${\cal A}_r$ spherical coordinates, into
$\sum_{j=1}^n \{ (p_0 + q_j{\sf S}_{e_j})^2 F^n_u (q;\zeta ) \} c_j
.$ The Cayley-Dickson number $q=q_0+q_1i_1+...+q_ni_n$ can be
written as $q=q_0+q_MM$, where $|M|=1$, $M$ is a purely imaginary
Cayley-Dickson number, $q_M\in {\bf R}$, $q_MM=q_1i_1+...+q_ni_n$,
since $q_0=Re (q)$. After a suitable automorphism $\theta : {\cal
A}_r\to {\cal A}_r$ we can take $\theta (q) = q_0+ q_Mi_1$, so that
$\theta (x)=x$ for any real number. The functions $[\sum_{j=1}^n
q_j^2 {\sf S}_{e_j}^2 c_j ]$ and $[\sum_{j=1}^n p_j^2 {\sf
S}_{e_j}^2 c_j ]$ are even by each variable $q_j$ and $p_j$
respectively. Therefore, we deduce in accordance with $(5)$ and
2$(3,4)$ and Corollary 6.1 with parameters $p_0=0$ and $\zeta =0$
and $c_j\in \{ -1, 1 \} $ for each $j$  that
\par $(6)$ $({\cal F}^n)^{-1} (1/[\sum_{j=1}^n
\{ \sum_{1\le k, b \le j} p_k{\sf S}_{e_k} p_b{\sf S}_{e_b}  \} c_j
];u;y;\zeta ) = - [g,e^{ N ([y],[q])})$ \\ in the ${\cal A}_r$
spherical coordinates, where $g=1/[\sum_{j=1}^n  q_j^2 c_j ]$, or
\par $(6.1)$ $({\cal F}^n)^{-1} (1/[\sum_{j=1}^n
\{  p_j^2{\sf S}_{e_j}^2  \} c_j ];u;y;\zeta ) = - [g,e^{ N
([y],[p])})$ \\ in the ${\cal A}_r$ Cartesian coordinates, where
$g=1/[\sum_{j=1}^n  p_j^2 c_j ]$, $N=y/|y|$ for $y\ne 0$, $N=i_1$
for $y=0$, $y=y_1i_1+...+y_ni_n\in {\cal A}_r$,
$[y]=(y_1,...,y_n)\in {\bf R}^n$, $([y],[q])= \sum_{j=1}^n y_jq_j$,
since ${\sf S}^2_{e_k} \cos (\phi + \zeta _k) = \cos (\phi + \zeta
_k + \pi ) = - \cos (\phi + \zeta _k)$ and ${\sf S}^2_{e_k} \sin
(\phi + \zeta _k) = \sin (\phi + \zeta _k + \pi ) = - \sin (\phi +
\zeta _k)$ for each $k$.
\par Particularly, we take $c_j=1$ for each $j=1,...,k_+$ and $c_j=-1$ for any
$j=k_+ + 1,...,n$, where $1\le k_+\le n$.  Thus the inverse Laplace
transform for $q_0=0$ and $\zeta =0$ in accordance with Formulas
2$(1-4)$ reduces to
\par $(7)$ $({\cal F}^n)^{-1} (1/[\sum_{j=1}^n
\{ \sum_{1\le k, b \le j} p_k{\sf S}_{e_k} p_b{\sf S}_{e_b}  \} c_j
];u;y;\zeta ) = $
\par $(2\pi )^{-n} \int_{{\bf R}^n} \exp (  {\bf i}
(q_1y_1+...+q_ny_n) )(1/ [\sum_{j=1}^{k_+} q_j^2 - \sum_{j=k_+ +1}^n
q_j^2])dq_1...dq_n$ \\ in the ${\cal A}_r$ spherical coordinates and
\par $(7.1)$ $({\cal F}^n)^{-1} (1/[\sum_{j=1}^n  p_j^2 {\sf
S}_{e_j}^2 c_j ];u;y;\zeta ) = $
\par $(2\pi )^{-n} \int_{{\bf R}^n} \exp (  {\bf i}
(p_1y_1+...+p_ny_n) )(1/ [\sum_{j=1}^{k_+} p_j^2 - \sum_{j=k_+ +1}^n
p_j^2])dp_1...dp_n$ \\ in the ${\cal A}_r$ Cartesian coordinates,\\
since for any even function its cosine Fourier transform coincides
with the Fourier transform.
\par The inverse Fourier transform $(F^{-1}g)(x)=(2\pi )^{-n}(Fg)(-x)=:\Psi _n$
of the functions $g=1/(\sum_{j=1}^n z_j^2)$ for $n\ge 3$ and ${\cal
P}(1/(\sum_{j=1}^2 z_j^2))$ for $n=2$ in the class of the
generalized functions is known (see \cite{gelshil} and \S \S 9.7 and
11.8 \cite{vladumf}) and gives
\par $(8)$ $\Psi _n(z_1,...,z_n) = C_n (\sum_{j=1}^n z_j^2)^{1-n/2}$
for $3\le n$, where $C_n = - 1/[(n-2) \sigma _n]$, $\sigma _n = 4\pi
^{n/2}/\Gamma ((n/2)-1)$ denotes the surface of the unit sphere in
${\bf R}^n$, $\Gamma (x)$ denotes Euler's gamma-function, while
\par $(9)$ $\Psi _2(z_1,z_2) = C_2 \ln (\sum_{j=1}^2 z_j^2)$ for
$n=2$, where $C_2= 1/(4\pi )$. \\ Thus the technique of \S 2 over
the Cayley-Dickson algebra has permitted to get the solution of the
Laplace operator.
\par For the function \par $(10)$ $P(x) = \sum_{j=1}^{k_+} x_j^2 -
\sum_{j=k_+ +1}^n x_j^2$ with $1\le k_+ <n$ the generalized
functions $(P(x)+{\bf i}0)^{\lambda }$ and $(P(x)-{\bf i}0)^{\lambda
}$ are defined for any $\lambda \in {\bf C} = {\bf R}\oplus {\bf
i}{\bf R}$ (see Chapter 3 in \cite{gelshil}). The function
$P^{\lambda }$ has the cone surface $P(z_1,...,z_n)=0$ of zeros, so
that for the correct definition of generalized functions
corresponding to $P^{\lambda }$ the generalized functions
\par \par $(11)$ $(P(x)+ c{\bf i}0)^{\lambda }=\lim_{0< c\epsilon ,
\epsilon \to 0 } (P(x)^2 + \epsilon ^2)^{\lambda /2} \exp ({\bf i}
\lambda arg (P(x) + {\bf i}c\epsilon ))$ \\ with either $c=-1$ or
$c=1$ were introduced. Therefore, the identity
\par $(12)$ $F(\Psi _{k_+,n-k_+})(x) = - (\sum_{j=1}^{k_+}x_j^2-
\sum_{j=k_+ +1}^n x_j^2) [F(\Psi _{k_+,n-k_+})(x)]^2$ or \par $(13)$
$F(\Psi ) = - 1/(P(x) + c {\bf i} 0)$ follows, where $c=-1$ or
$c=1$.
\par The inverse Fourier transform in the class of the
generalized functions is: \par $(14)$ $F^{-1}((P(x)+c{\bf
i}0)^{\lambda })(z_1,...,z_n) = \exp (- \pi c(n-k_+){\bf i}/2)
2^{2\lambda +n} \pi ^{n/2} \Gamma (\lambda +n/2)(Q(z_1,...,z_n) -
c{\bf i}0)^{- \lambda - n/2)}/[\Gamma (-\lambda )|D|^{1/2}]$ \\ for
each $\lambda \in {\bf C}$ and $n\ge 3$ (see \S IV.2.6
\cite{gelshil}), where $D=\det (g_{j,k})$ denotes a discriminant of
the quadratic form $P(x)=\sum_{j,k=1}^n g_{j,k}x_jx_k$, while $Q(y)=
\sum_{j,k=1}^n g^{j,k}x_jx_k$ is the dual quadratic form so that
$\sum_{k=1}^n g^{j,k}g_{k,l}=\delta ^j_l$ for all $j, l$; $\delta
^j_l=1$ for $j=l$ and $\delta ^j_l=0$ for $j\ne l$.  In the
particular case of $n=2$ the inverse Fourier transform is given by
the formula:
\par $(15)$ $F^{-1}((P(x)+c{\bf i}0)^{-1})(z_1,z_2) = - 4^{-1}|D|^{-1/2} \exp (-
\pi c(n- k_+){\bf i}/2) \ln (Q(z_1,...,z_n) - c{\bf i}0).$ Making
the inverse Fourier transform $F^{-1}$ of the function $- 1/(P(x) +
{\bf i} 0)$ in this particular case of $\lambda =-1$ we get two
complex conjugated fundamental solutions
\par $(16)$ $\Psi _{k_+,n-k_+} (z_1,...,z_n) = - \exp (\pi c(n-k_+)
{\bf i}/2)\Gamma ((n/2) -1) (Q(z_1,...,z_n) + c{\bf
i}0)^{1-(n/2)}/(4\pi ^{n/2})$ for $3\le n$ and $1\le k_+<n$, while
\par $(17)$ $\Psi _{1,1}(z_1,z_2) = 4^{-1} \exp (\pi c(n-k_+){\bf i}/2)
\ln (Q(z_1,z_2) + c{\bf i}0) $ for $n=2$, where either $c= 1$ or
$c=-1$.
\par Generally for the operator $A$ given by Formula $(1)$
we get $P(x) = P_0(x) + P_i(x)$, where $P_0(x)=\sum_{j=1}^n x_j^2Re
(c_j)$ and $P_i(x)= \sum_{j=1}^n x_j^2 Im(c_j)$ are the real and
imaginary parts of $P$, $Im (z) = z- Re (z)$ for any Cayley-Dickson
number $z$. Take ${\bf l} =i_{2^r}$ and consider the form
$P(x)+\epsilon c{\bf l}$ with $\epsilon \ne 0$ and either $c=1$ or
$c=-1$, then $P(x)+\epsilon c{\bf l}\ne 0$ for each $x\in {\bf
R}^n$. We put \par $(18)$ $(P(x)+ c{\bf l}0)^{\lambda }=\lim_{0<
c\epsilon , \epsilon \to 0 } (P(x)^2 + \epsilon ^2)^{\lambda /2}
\exp ({\bf i} \lambda Arg (P(x) + {\bf l}c\epsilon )).$ Consider
$\lambda \in {\bf R}$, the generalized function $(P(x)^2 + \epsilon
^2)^{\lambda /2} \exp ({\bf i} \lambda Arg (P(x) + {\bf l}c\epsilon
))$ is non-degenerate and for it the Fourier transform is defined.
The limit $\lim_{0< c\epsilon , \epsilon \to 0 }$ gives by our
definition the Fourier transform of $(P(x)+ c{\bf l}0)^{\lambda }$.
Since \par $(19)$ $c_j (\beta _j+ \sum_{1\le k\le n, k\ne j}
c_j^{-1}c_k \beta _k) = \sum_{j=1}^nc_j\beta _j$ \\ for all $\beta
_j\in \bf R$ and any $1\le j \le n$ in accordance with the
conditions imposed on $c_j$ at the beginning of this section and
${\bf i}N_j=N_j{\bf i}$ for each $j$, the Fourier transform with the
generator $\bf i$ can be accomplished subsequently by each variable
using Identity $(19)$. The transform $x_j\mapsto (c_j)^{1/2}x_j$ is
diagonal and $|(...((c_1^{1/2}c_2^{1/2})...)c_n^{1/2}|=1$, so we can
put $|D|=1$. \par Each Cayley-Dickson number can be presented in the
polar form $z=|z|e^{\phi M}$, $\phi \in {\bf R}$, $|\phi |\le \pi $,
$M$ is a purely imaginary Cayley-Dickson number $|M|=1$, $Arg (z) =
(\phi + 2\pi k)M$ has the countable number of values, $k\in {\bf Z}$
(see \S 3 in \cite{ludoyst,ludfov}). Therefore, we choose the branch
$z^{1/2} = |z|^{1/2} \exp ( (Arg z)/2)$, $|z|^{1/2}>0$ for $z\ne 0$,
with $|Arg (z)|\le \pi $, $Arg (M) = M\pi /2$ for each purely
imaginary $M$ with $|M|=1$.
\par We treat the iterated integral as in \S 6, i.e. with the same order of brackets.
Taking initially $c_j\in {\bf R}$ and considering the complex
analytic extension of formulas given above in each complex plane
${\bf R}\oplus N_j{\bf R}$ by $c_j$ for each $j$ by induction from
$1$ to $n$, when $c_j$ is not real in the operator $A$, $Im (c_j)\in
{\bf R}N_j$, we get the fundamental solutions for $A$ with the form
$(P(x)+ c{\bf l}0)^{\lambda }$ instead of $(P(x)+ c{\bf
i}0)^{\lambda }$ with multipliers
$(...(c_1^{c/2}c_2^{c/2})...)c_n^{c/2}$ instead of $\exp (\pi
c(n-k_+) {\bf i}/2)$ as above and putting $|D|=1$. Thus \par $(20)$
$\Psi (z_1,...,z_n) = -  \Gamma ((n/2) -1) (P^*(z_1,...,z_n) - c{\bf
l}0)^{1-(n/2)}[(...(c_1^{c/2}c_2^{c/2})...)c_n^{c/2}]^*/(4\pi
^{n/2})$ for $3\le n$, while
\par $(21)$ $\Psi (z_1,z_2) = 4^{-1} [c_1^{c/2} c_2^{c/2}]^* Ln
(P^*(z_1,z_2) - c{\bf l}0)$ for $n=2$, \\ since $c_j^*= c_j^{-1}$
for $|c_j|=1$, $y_jq_j=y_j(c_j^{c/2})^*q_jc_j^{1/2}$, while \par
$(...(dc_1^{c/2}q_1 dc_2^{c/2}q_2)...)dc_n^{c/2}q_n] = dq_1 ...dq_n
[(...(c_1^{c/2}c_2^{c/2})...)c_n^{c/2}]$ and \par
$|(...(c_1^{c/2}c_2^{c/2})...)c_n^{c/2}|=1$.

\par {\bf 36. Partial differential equations
with polynomial real coefficients.} Let \par $(1)$ $A= \sum_{|\alpha
|\le m} a_{\alpha }(q) \partial ^{\alpha }_q$, \par $a_{\alpha }(q)
= \sum_{\beta } a_{\alpha , \beta } q^{\beta }$, $q^{\beta } :=
q_1^{\beta _1}...q_n^{\beta _n}$, $a_{\alpha , \beta }$ and $f$ have
values as in \S 28, and $Af$ be an original. Using the transform in
the ${\cal A}_r$ Cartesian coordinates we take $q_j=t_j$ for each
$j$, while using the transform in ${\cal A}_r$ spherical coordinates
we choose $q_j=s_j(t)$ for each $j$. Then
\par $(2)$ ${\cal F}^n(Af;u;p;\zeta ) = \sum_{\beta } (-1)^{|\beta
|} {\sf S}_{\beta }(p) \partial ^{\beta }_p$ \par $ [\sum_{\beta }
a_{\alpha ,\beta } ( [p_0+p_1{\sf S}_{e_1}]^{\alpha _1} p_2^{\alpha
_2} {\sf S}^{\alpha _2}_{e_2}...p_n^{\alpha _n} {\sf S}^{\alpha
_n}_{e_n})] F^n(p;\zeta ) = G^n(p;\zeta )$ \\ in the ${\cal A}_r$
spherical coordinates and
\par $(2.1)$ ${\cal F}^n(Af;u;p;\zeta ) = \sum_{\beta } (-1)^{|\beta
|} {\sf S}_{\beta }(p) \partial ^{\beta }_p$ \par $ (\sum_{\beta }
a_{\alpha ,\beta } [p_0+p_1{\sf S}_{e_1}]^{\alpha _1} [p_0+p_2{\sf
S}_{e_2}]^{\alpha _2}...[p_0+p_n {\sf S}_{e_n}]^{\alpha _n})
F^n(p;\zeta ) = G^n(p;\zeta )$ \\ in the ${\cal A}_r$ Cartesian
coordinates (see Theorems 12 and 13 above). It may happen that the
second differential equation is simpler than the initial one:
\par $(3)$ $Af=g$. \\ For example, when coefficients depend only on
one variable $t_n$, then the second differential equation is
ordinary and linear.
\par {\bf 37. Noncommutative transforms of products and convolutions of
functions in the ${\cal A}_r$ spherical coordinates.}
\par For any Cayley-Dickson number $z=z_0i_0+...+z_{2^r-1}i_{2^r-1}$
we consider projections \par $(1)$ $\theta _j(z)=z_j$, $z_j\in {\bf
R}$ or ${\bf C}_{\bf i}$ or ${\bf H}_{J,K,L}$, $j=0,...,2^r-1$, $\theta _j(z) = \pi _j(z)i_j^*$, \\
given by Formulas 2$(5,6)$ and 33$(17)$. We define the following
operators
\par $(2)$ ${\cal R}_{\alpha ,j} (F^n(p;\zeta )) :=
F^n(p_0,(-1)^{\alpha _1} p_1,...,(-1)^{\alpha _{j+1-\delta _{j,n}}}
p_{j+1-\delta _{j,n}}, p_{j+2- \delta _{j,n}},$ \\ $...,p_n; \zeta
_0, (-1)^{\alpha _1} \zeta _1 +\pi \alpha _1/2,...,(-1)^{\alpha
_{j+1-\delta _{j,n}}} \zeta _{j+1- \delta _{j,n}} +\pi \alpha _{j+1-
\delta _{j,n}} /2,
\zeta _{j+2- \delta _{j,n}},...,\zeta _n)$ \\
on images $F^n$, $2^{r-1}\le n \le 2^r-1$, $j=0,...,n$. For $\alpha
_j$ and $\beta _j\in \{ 0, 1\} $ their sum $\alpha _j+\beta _j$ is
considered by $(mod ~ 2)$, i.e. in the ring ${\bf Z}_2={\bf
Z}/(2{\bf Z})$, for two vectors $\alpha $ and $\beta \in \{ 0, 1 \}
^{2^r-1}$ their sum is considered componentwise in ${\bf Z}_2$. Let
$$(3)\quad {\cal F}^n(f;u;p;\zeta ) = \sum_{j=0}^n \sum_{k=0}^{2^r-1}\theta _j
({\cal F}^n(\theta _k(f);u;p;\zeta ))i_ki_j,$$ also $F^n_j(p;\zeta )
:= \sum_{k=0}^{2^r-1}\theta _j ({\cal F}^n(\theta _k(f);u;p;\zeta
))i_k$ for an original $f$, where $u(p,t;\zeta )$ is given by
Formulas 2$(1,2,2.1)$. If $f$ is real or ${\bf C}_{\bf i}$ or ${\bf
H}_{J,K,L}$ -valued, then $F^n_j=\theta _j(F^n)$.
\par {\bf Theorem.} {\it If $f$ and $g$ are two originals, then
\par $(4)$ ${\cal F}^n(fg;u;p;\zeta ) = \sum_{j=0}^n \sum_{\alpha ,
\beta \in \{ 0, 1 \} ^n} (-1)^{\alpha _{j+1}(1-\delta
_{j+1,n})}({\cal R}_{\alpha ,j} (F^n_j(p-q_0;\zeta -\eta ))*({\cal
R}_{\beta ,j} (G^n_j(p+q_0-p_0;\eta ))i_j,$ \par $(4.1)$ ${\cal
F}^n(f*g;u;p;\zeta ) = \sum_{j=0}^n \sum_{\alpha , \beta \in \{ 0, 1
\} ^n} (-1)^{\alpha _{j+1}(1-\delta _{j+1,n})}({\cal R}_{\alpha ,j}
(F^n_j(p;\zeta -\eta )) ({\cal R}_{\beta ,j} (G^n_j(p;\eta ))i_j,$
\\ whenever ${\cal F}^n(fg)$, ${\cal F}^n(f)$, ${\cal F}^n(g)$ exist,
where $1\le n \le 2^r-1$, $2\le r$; $\alpha _k + \beta _k =1 ~ (mod
~ 2)$ for $k\le j$ or $k=j+1=n$, $\alpha _k + \beta _k =0 ~ (mod ~
2)$ for $k=j+1<n$ and $\alpha _k=\beta _k=0$ for $k>j+1$ in the
$j$-th addendum on the right of Formulas $(4,4.1)$; the convolution
is by $(p_1,...,p_n)$ in $(4)$, at the same time $q_0\in {\bf R}$
and $\eta \in {\cal A}_r$ are fixed.}
\par {\bf Proof.} The product of two originals can be written in the
form:
\par $(5)$ $f(t)g(t)= \sum_{j=0}^{2^r-1} \sum_{k,l: ~ i_ki_l=i_j}
\theta _k(f(t)) \theta _l(g(t)) i_j$. \\ The functions $\theta
_k(f)$ and $\theta _l(g)$ are real or ${\bf C}_{\bf i}$ or ${\bf
H}_{J,K,L}$ valued respectively. The non-commutative transform of
$fg$ is:
$$(6)\quad {\cal F}^n(fg)(p;\zeta ) = \int_{{\bf R}^n} f(t)g(t)\exp
( - u(p,t;\zeta ))dt =$$  $$ \{ \int_{{\bf R}^n} (f(t)g(t))
e^{-p_0s_1} \cos (p_1s_1+\zeta _1 ) i_0 dt \} +$$ $$ \{
\sum_{j=2}^{n-1} \int_{{\bf R}^n} (f(t)g(t)) e^{-p_0s_1} \sin
(p_1s_1+\zeta _1) ... \sin (p_{j-1}s_{j-1}+\zeta _{j-1}) \cos
(p_js_j+\zeta _j ) i_{j-1} dt \} +$$
$$\int_{{\bf R}^n} (f(t)g(t)) e^{-p_0s_1}
\sin (p_1s_1+\zeta _1) ... \sin (p_ns_n+\zeta _n) i_n dt.$$ On the
other hand,
$$(7)\quad \int_{{\bf R}^n} f(t)g(t) e^{-p_0s_1 + {\bf i} \sum_{j=1}^k (p_js_j+\zeta _j)\gamma _j}
dt =$$ $$\int_{{\bf R}^n}(\int_{{\bf R}^n} f(t) e^{-(p_0-q_0)s_1 +
{\bf i}\sum_{j=1}^k ((p_j-q_j)s_j + \zeta _j - \eta _j)\gamma _j}
dt) (\int_{{\bf R}^n} g(t) e^{-q_0s_1 + {\bf i}\sum_{j=1}^k
(q_js_j+\eta _j)\gamma _j}dt)dq,$$ where $k=1,2,...,n$, $ ~ \gamma
_j\in \{ -1, 1 \} $. Therefore, using Euler's formula $e^{{\bf i}
\phi } = \cos (\phi ) + {\bf i} \sin (\phi )$  and the trigonometric
formulas $\cos (\phi + \psi ) = \cos (\phi ) \cos (\psi ) - \sin
(\phi )\sin (\psi )$, $\quad \sin (\phi + \psi ) = \sin (\phi ) \cos
(\psi ) + \cos (\phi )\sin (\psi )$ for all $\phi , ~ \psi \in {\bf
R}$, and Formulas $(6,7)$, we deduce expressions for $\theta
_j({\cal F}^n(fg))$. We get the integration by $q_1,...,q_n$, which
gives convolutions by the $p_1,...,p_n$ variables. Here $q_0\in {\bf
R}$ and $\eta \in {\cal A}_r$ are any marked numbers. Thus from
Formulas $(5-7)$ and 2$(1,2,2.1,4)$ we deduce Formula $(4)$.
\par Moreover, one certainly has
$$(8)\quad \int_{{\bf R}^n} (f*g)(t) e^{-p_0s_1 + {\bf i} \sum_{j=1}^k (p_js_j+\zeta _j)\gamma _j}
dt =$$ $$(\int_{{\bf R}^n} f(t) e^{- p_0s_1 + {\bf i}\sum_{j=1}^k
(p_js_j +\zeta _j - \eta _j)\gamma _j} dt) (\int_{{\bf R}^n} g(t)
e^{-p_0s_1 + {\bf i}\sum_{j=1}^k (p_js_j+\eta _j)\gamma _j}dt)$$ for
each $1\le k\le n$, $\gamma _j \in \{ -1, 1 \} $, since $s_j(t) =
s_j(t-\tau ) + s_j(\tau )$ for all $j=1,...,n$ and $t, ~ \tau \in
{\bf R}^n$. Thus from Relations $(6,8)$ and 2$(1,2,2.1,4)$ and
Euler's formula one deduces expressions for $\theta _j({\cal
F}^n(f*g))$ and Formula $(4.1)$.
\par {\bf 38. Moving boundary problem.}
\par Let us consider a boundary problem
\par $(1)$ $Af=g$ in the half-space $t_n\ge \phi (t_n)$, where $\phi (0)=0$ and
$\phi (t_n)<t_n$ for each $0\le t_n\in {\bf R}.$ Suppose that the
function $t_n-\phi (t_n) =: \psi (t_n)$ is differentiable and
bijective. For example, if $0<v<1$ and $\phi (t_n)=vt_n$, then the
boundary is moving with the speed $v$. Make the change of variables
$y_n = \psi (t_n)$, $y_1 = t_1,$...,$y_{n-1} = t_{n-1}$, then $t_n =
\psi ^{-1}(y_n)$ and $dt_n = ds_n = (dt_n/dy_n) dy_n$ and due to
Theorem 25 we infer that
$$(2)\quad {\cal F}^n (\sum_{|\alpha |\le m} {\bf b}_{\alpha } \partial ^{\alpha }_s \chi
_{y_n\ge 0} f(t);p;\zeta ) =   \sum_{|\alpha |\le m, 0\le q_n\le
\alpha _n-1} {\bf b}_{\alpha }(\delta _{0,{\alpha _n}} -1)$$  $$
(p_0+{\sf S}_{e_1}p_1)^{\alpha _1}p_2^{\alpha _2}...p_{n-1}^{\alpha
_{n-1}}p_n^{\alpha _n-q_n-1} {\sf S}_{\alpha - \alpha _1e_1- (q_n+1)
e_n}{\cal F}^{n-1,y^n}(\partial ^{q_n}_{t_n}w(y),u(p,(y^n);\zeta
);p;\zeta )$$
$$ + \sum_{|\alpha |\le m} {\bf b}_{\alpha } (p_0+{\sf
S}_{e_1}p_1)^{\alpha _1}p_2^{\alpha _2}...p_n^{\alpha _n} {\sf
S}_{\alpha - \alpha _1e_1} {\cal F}^n( \chi _{y_n\ge 0} (y)
w(y);p;\zeta ) = G^n(p;\zeta )$$ in the ${\cal A}_r$ spherical
coordinates and
$$(2.1)\quad {\cal F}^n (\sum_{|\alpha |\le m} {\bf a}_{\alpha } \partial ^{\alpha }_t \chi
_{y_n\ge 0} f(t);p;\zeta ) =  \sum_{|\alpha |\le m, 0\le q_n\le
\alpha _n-1} {\bf a}_{\alpha }(\delta _{0,{\alpha _n}} -1)$$
$$ (p_0+{\sf S}_{e_1}p_1)^{\alpha _1}(p_0+p_2{\sf S}_{e_2})^{\alpha
_2}...(p_0+p_{n-1}{\sf S}_{e_{n-1}})^{\alpha _{n-1}} (p_0+p_n {\sf
S}_{e_n})^{\alpha _n - q_n-1}$$ $$ {\cal F}^{n-1,y^n}(\partial
^{q_n}_{t_n}w(y),u(p,(y^n);\zeta );p;\zeta )$$
$$ + \sum_{|\alpha |\le m} {\bf a}_{\alpha } (p_0+{\sf
S}_{e_1}p_1)^{\alpha _1}(p_0+p_2{\sf S}_{e_2})^{\alpha
_2}...(p_0+p_n{\sf S}_{e_n})^{\alpha _n}  {\cal F}^n( \chi _{y_n\ge
0} (y) w(y);p;\zeta ) = G^n(p;\zeta )$$ in the ${\cal A}_r$
Cartesian coordinates, where $w(y) := f(t(y)) (dt_n/dy_n)$.
\par Expressing ${\cal F}^n( \chi _{y_n\ge 0} (y) w(y);p;\zeta )$ through
$G^n(p;\zeta )$ and the boundary terms ${\cal F}^{n-1,y^n}(\partial
^{q_n}_{t_n}w(y),u(p,(y^n);\zeta );p;\zeta )$ as in \S 28.3  and
making the inverse transform 8$(4)$ or 8.1$(1)$, or using the
integral kernel $\xi $ as in \S 28.5, one gets a solution $w(y)$ or
$f(t)=w(y(t)) (dy_n(t_n)/dt_n)$.
\par {\bf 39. Partial differential equations with discontinuous
coefficients.}
\par Consider a domain $U$ and its subdomains $U\supset U_1\supset ...\supset U_k$
satisfying Conditions 28$(D1,D4,i-vii)$ so that coefficients of an
operator $A$ (see 28$(2)$) are constant on $Int (U_k)$ and on $V_1=
U\setminus Int (U_1)$, $V_2=U_1\setminus Int
(U_2)$,...,$V_k=U_{k-1}\setminus Int (U_k)$ and are allowed to be
discontinuous at the common borders $\partial V_j\cap
\partial U_j$ for each $j=1,...,k$. Each function $f\chi _{U_j}$ is an original
on $U$ or a generalized function with the support $supp (f\chi
_{U_j})\subset U_j$ if $f$ is an original or a generalized function
on $U$. Choose operators $A^j$ with constant coefficients on $U^j$
and $A^j|_{Int (V_j)}=0$, where $U^0=U$, so that $A|_{U_k}=A^k$,...,
$A|_{U_j}=A^j+...+A^k$,..., $A|_U=A^0+...+A^k$. Therefore, in the
class of originals or generalized functions on $U$ the problem (see
28$(1,2)$) can be written as
\par $(1)$ $Af=g$, or
\par $(2)$ $A^0f\chi _{V_1}=g\chi _{V_1}$,...,$A^{k-1}f\chi_{V_k}=g\chi
_{V_k}$, $A^kf\chi _{U_k}=g\chi _{U_k},$ \\ since $\chi
_{V_1}+...+\chi _{V_k}+\chi _{U_k}=\chi _U$. Thus the equivalent
problem is: \par $(3)$ $A^0f^0=g^0$, $A^1f^1=g^1$,...,$A^kf^k=g^k$ \\
with $f^k=f\chi _{U_k}$, $g^k=g\chi _{U_k}$, also $f^j=f\chi
_{V_{j+1}}$, $g^j=g\chi _{V_{j+1}}$ for each $j=0,...,k-1$. On
$\partial U$ take the boundary condition in accordance with
28$(5.1)$. With any boundary conditions in the class of originals or
generalized functions on additional borders $\partial U_j\setminus
\partial U$ given in accordance with 28$(5.1)$ a solution $f^j$ on
$U^j$ exists, when the corresponding condition 8$(3)$ is satisfied
(see Theorems 8 and 28.1). \par Each problem $A^jf^j=g^j$ can be
considered on $U_j$, since $supp(g^j)\subset U_j$. Extend $f^j$ by
zero on $U\setminus V_j$ for each $0\le j \le k-1$. When the right
side of 28$(6)$ is non-trivial, then $f^j$ is non-trivial. If
$f^{j-1}$ is calculated, then the boundary conditions on $\partial
U^j\setminus \partial U$ can be chosen in accordance with values of
$f^{j-1}$ and its corresponding derivatives $(\partial ^{\beta }
f^{j-1} /\partial \nu ^{\beta })|_{(\partial  U^j\setminus \partial
U)}$ for some $\beta < ord (A^j)$ in accordance with the operator
$A^j$ and the boundary conditions 28$(5.1)$ on the boundary
$\partial U^j\setminus \partial U$. Having found $f^j$ for each
$j=0,...,k$ one gets the solution $f=f^0+...+f^k$  on $U$ of Problem
$(1)$ with the boundary conditions 28$(5.1)$ on $\partial U$.
\par {\bf 40. Remark.} The multiparameter noncommutative transform
over the Cayley-Dickson algebras presented above is the natural
generalization of the usual complex one-parameter Laplace transform.
It opens new opportunities for solving partial differential
equations of different types.
\par It may happen that Theorem 13 is simpler to use, than Theorem 21 for
partial differential equations with real variables. Theorem 13 has
an advantage that it can be simpler used for partial differential
equations of complex and hyper-complex variables, because each pair
$(p_l + p_ji_l^*i_j)$ for $l\ne j$ is the complex variable. In these
variants boundary conditions may be for $F^k(p;\zeta )$ on a
hyperplane $Re (p) = a$ in ${\cal A}_r$.
\par As it was seen above the appearing integrals are by
multidimensional domains. For their calculations the Fubini's
theorem, residues, Jordan Lemma and tables of known integrals also
can be used. Generally in computational mathematics integrals are
easier to calculate, than to solve partial differential equations
numerically. As a rule iterations of algorithms for integrals
converge faster, than iterations of numerical methods for partial
differential equations.
\par Functions with octonion values may be used to resolve
systems of partial differential equations. Using conjugations of
Cayley-Dickson numbers one gets the transition between operators
with coefficients either on the left or on the right of partial
derivatives: $[(\partial ^{\alpha } f(x))c_{\alpha }]^* = c_{\alpha
}^*(\partial ^{\alpha } f(x))^*$, particularly, $(\partial ^{\alpha
} f(x))^* = \partial ^{\alpha } f^*(x)$ for $x\in {\bf R}^n$,
$\partial ^{\alpha }=\partial ^{\alpha }_x$.
\par Using of Formulas 2$(5,6)$ gives variables $t_j=z_j$ for $z\in {\cal
A}_r$. So one can consider a class of super-differentiable originals
$f(z)$, $z\in V\subset {\cal A}_r$. In the class of piecewise on
open subsets super-differentiable originals $f(z)$, $z\in V\subset
{\cal A}_r$, with $t_j=z_j$ for each $j=1,...,n$, $n=2^r-1$, in the
fixed $z$-representations we get the noncommutative transform for
$f(z)\chi _V(z)$ relative to the Cayley-Dickson variable $z\in {\cal
A}_r$. Therefore, the results given above transfer on this variant
also.
\par Theorem 17 also opens new opportunities to investigate and
solve certain types of nonlinear partial differential equations
using previous results on spectral theory of functions of operators
\cite{luspraaca,lujmsalop}. For example, analytic functions $q(z)$
in Theorem 17 permit to consider nonlinear operators $q(\sigma )$,
where $\sigma f(z) := \sum_{j=0}^{2^r-1} (\partial f(z)/\partial
z_j)i_j$. It is planned to study in the next paper.
\par Partial differential equations with periodic $g$
and $f$ with vector period corresponding to $Q^n$ may be considered
also. Certainly others classes of smoothness, for example, Sobolev's
or generalized functions can also be considered. It is planned in a
next paper to consider this and also problems with boundary
conditions as well as with non-constant coefficients in more
details. \par The technique described above permits to consider
partial differential equations of different types and write their
solutions in integral forms. If appearing integrals can be
calculated in elementary or special of generalized functions, then
this gives the explicit formulas in terms of known functions. In
conjunction with the line integration over the Cayley-Dickson
algebras it permits to solve some types of non linear partial
differential equations. The multiparameter Laplace transform over
the Cayley-Dickson algebras takes into account the boundary
conditions. It naturally means the treatment of systems of partial
differential equations due to the multidimensionality of the
Cayley-Dickson algebras.

\par Department of Applied Mathematics,
\par Moscow State Technical University MIREA,
av. Vernadsky 78,
\par Moscow, Russia
\par e-mail: sludkowski@mail.ru
\end{document}